\documentclass[11pt,oneside,reqno]{amsart}
\usepackage{epsfig,psfrag,subfigure}
\usepackage{amsmath,amssymb,amsfonts,latexsym, mathrsfs,euscript}
\usepackage{graphicx,graphics,cite}
\usepackage{wrapfig}
\usepackage[usenames]{color}
\usepackage{geometry}                
\usepackage[table]{xcolor}

\usepackage{fullpage}
\usepackage{amsthm}

\makeatletter \@addtoreset{figure}{section}
\def\thefigure{\thesection.\@arabic\c@figure} \def\fps@figure{h, t}
\@addtoreset{table}{section}
\def\thetable{\thesection.\@arabic\c@table}\def\fps@table{h, t}
\@addtoreset{equation}{section}

\makeatother

\theoremstyle{plain}

\theoremstyle{remark}
\newtheorem{note}{Note}

\theoremstyle{definition}
\newtheorem{example}{Example}

\newcommand{\beq}{\begin{equation}}
\newcommand{\eeq}{\end{equation}}
\newcommand{\beqa}{\begin{eqnarray}}

\newcommand{\eeqa}{\end{eqnarray}}

\newcounter{nfig}

\newcommand{\me}{\ensuremath{\mathrm{e}}}
\newcommand{\mi}{\ensuremath{\mathrm{i}}}
\newcommand{\dif}{\ensuremath{\mathrm{d}}}
\DeclareMathOperator{\re}{Re}
\DeclareMathOperator{\im}{Im}
\renewcommand{\Re}{\re}
\renewcommand{\Im}{\im}

\def\erf{\text{erf}}

\begin{document}

\title[Analysis of compact DRB FD schemes]{Analysis and development of compact finite difference schemes with optimized numerical dispersion relations}
\author{Yi-Hung Kuo}
\email{yhkuo@g.ucla.edu}
\author{Long Lee}
\email{llee@uwyo.edu}
\author{Gregory Lyng}
\email{glyng@uwyo.edu}
\address{Department of Mathematics, University of Wyoming, Laramie, 82071, USA}
\date{\today}

\maketitle

\begin{abstract}
Finite difference approximation, in addition to Taylor truncation errors, introduces numerical dispersion-and-dissipation errors into numerical solutions of partial differential equations. We analyze a class of finite difference schemes which are designed to
minimize these errors (at the expense of formal order of accuracy), and
we analyze the interplay between the Taylor truncation errors and the dispersion-and-dissipation errors during mesh refinement.
In particular, 
we study the numerical dispersion relation of the fully discretized non-dispersive transport equation in one and two space dimensions. We derive the numerical phase error and the $L^2$-norm error of the solution in terms of the dispersion-and-dissipation error. Based on our analysis, we investigate the error dynamics among various optimized compact schemes and the unoptimized higher-order generalized Pad\'e compact schemes, taking into account four important factors, namely, (i) error tolerance, (ii) computer memory capacity, (iii) resolvable wavenumber, and (iv) CPU/GPU time. The dynamics shed light on the principles of designing suitable optimized compact schemes for a given problem. Using these principles as guidelines, we then propose an optimized scheme that prescribes the numerical dispersion relation before finding the corresponding discretization. This approach produces smaller numerical dispersion-and-dissipation errors for linear and nonlinear problems, compared with the unoptimized higher-order compact schemes and other optimized schemes developed in the literature. 
Finally, we discuss the difficulty of developing an optimized composite boundary scheme for problems with non-trivial boundary conditions. We propose a composite scheme that introduces a buffer zone to connect an optimized interior scheme and an unoptimized boundary scheme. Our numerical experiments show that this strategy produces small $L^2$-norm error when a wave packet passes through the non-periodic boundary.
\end{abstract}


\section{Introduction}\label{sec:intro}
\subsection{Background}\label{ssec:back}
Among the virtues possessed by finite difference (FD) schemes as tools for solving partial differential equations (PDEs) are ease of implementation and flexibility with regard to boundary conditions. Many problems of interest (e.g., in aeroacoustics and turbulent fluid flows) feature wide ranges of temporal and spatial scales, and the construction of accurate (approximate) solutions of these problems over all scales, especially within the constraint of limited computational power, remains a substantial challenge. Thus, in the last several decades, different criteria---in addition to the traditional truncation error of Taylor series expansion---have been used to analyze and design FD schemes which are somehow ``optimal.''  
Trefethen \cite{t82} was among the pioneers who recognized that errors due to FD approximation are not random perturbations, but a systematic superposition of dispersion-and-dissipation errors of various orders associated with the approximation.  
He analyzed the dispersion relations and numerical group velocities of various fully discretized finite difference schemes solving the non-dispersive simple transport equation, and he illustrated that local truncation error is only one factor contributing to the inaccuracy of a given scheme. 
%
%
Indeed, much subsequent work has shown that reducing or minimizing the dispersion-and-dissipation errors allows one to develop a FD scheme that, in principle, can make the most of a coarse-grid calculation when further refinement of the mesh would be too expensive \cite{lele92,tw93,ht94}. 

We refer to this class of algorithms as dispersion-relation-based (DRB) schemes. For a typical DRB scheme, a process is employed to determine a set of spatial or temporal differencing coefficients for the FD stencils so that the numerical approximations for the derivatives produce the least numerical dispersion and/or dissipation errors at the expense of formal order of accuracy. It is known that when considering a refinement process for which the mesh size continues to decrease, higher order schemes are superior asymptotically. However, for a given mesh size, there may exist a DRB scheme of lower formal order-of-accuracy with superior accuracy on that grid for an appropriate family of initial data \cite{ht94}. While there is no DRB scheme that is superior to its higher-order counterpart on all meshes, we give here an analysis of the interplay between the Taylor truncation error and the dispersion-and-dissipation error when refining meshes for FD schemes. 
This analysis sheds light on how to choose/design a proper scheme for a given problem, whether an un-optimized higher-order scheme or an optimized DRB scheme, based on the problem at hand. 


\subsection{DRB schemes: a little history}\label{ssec:drb}
Generalized Pad\'e compact (implicit) FD schemes with spectral-like resolution were introduced---in the course of establishing improved representation of a range of scales in the evaluation of spatial derivatives---by Lele in his seminal paper \cite{lele92}.
Lele introduced a general framework of finding the coefficients of a compact FD scheme approximating the first or the second spatial derivative, based on the truncation error of the Taylor series expansion, and then he described a compact DRB scheme. 
Lele's approach proceeds in two steps. First, the differencing coefficients are chosen to meet the desired formal order of accuracy, and, second, the remaining coefficients---which could be used to generate a higher-order-accurate scheme---are found by equations deduced from matching the so-called effective wavenumber and the exact wavenumber at some chosen values. The resulting FD schemes have a lower formal order of accuracy, but they have much better resolution at high effective wavenumber than that of  corresponding higher-order schemes. Although the choice of matching locations was ad hoc, and no attempt was made to optimize the choice, the result turns out to be remarkably good. 

At roughly the same time, Tam and Webb \cite{tw93} proposed a technique for minimizing the dispersion-and-dissipation error for spatial and temporal FD discretization. They considered explicit FD schemes to approximate spatial derivatives. Unlike Lele's scheme, after some of the differencing coefficients are chosen to meet the desired formal order of accuracy, the remaining coefficients are determined using the equations obtained by minimizing an integral of the difference between the effective and the exact wavenumber over the range of interest. A similar minimization procedure was applied to the error of effective angular frequency to find the coefficients of a chosen explicit multi-step time integrator. The authors called the resulting scheme a \emph{dispersion-relation-preserving (DRP) FD scheme}. In a variation on this theme, Haras and Ta'asan \cite{ht94} introduced a DRB scheme, in which the coefficients of compact FD schemes were found by minimizing a ``global'' truncation error using a constrained minimization problem that incorporated initial conditions 
into the kernel of the objective function.  Similarly, Kim and Lee \cite{kl96} proposed an optimized compact scheme by introducing a weight function in the kernel of the objective function to streamline the minimization process. 

There are now many other DRB schemes developed using different optimization strategies; these strategies are usually based on the requirements of the considered application. We do not give an exhaustive list here of existing DRB schemes; indeed, after more than two decades of development, there are now so many different forms of DRB methods that we could hardly compile such a list.
However, to give a sense of the variety of DRB methods, we do mention that examples of existing optimized DRB schemes include upwind-type schemes  \cite{pirozzi09, sw08, zyhs07, zc98},  multi-scale and multi-grid-size schemes \cite{bbmb07,tk03}, DRB schemes with non-uniform girds for general geometry \cite{cl01,sw08}, composite boundary schemes \cite{j07, j11,kl97,kim07},  optimized prefactored compact schemes \cite{az03}, and optimized DRB schemes that couple space-time discretization \cite{rbsl06, srb11}.  In addition to the development of optimized DRB schemes, there are also comparison studies of the performance of optimized DRB schemes \cite{ht94,ps02, p06, zigg00}, applications of DRB schemes in fluid flow and non-linear waves \cite{cs09, cls09, skcl11, srb11}, and error dynamics of optimized DRB schemes \cite{sds07}.

\subsection{Numerical dispersion-and-dissipation errors}
Suppose that a PDE admits 
a solution of the form
\begin{equation}\label{eq:sol}
u(x,t)=\me^{\mi(kx-\omega t)}\,,
\end{equation}
where $\mi=\sqrt{-1}$, $k$ is the wavenumber, and $\omega$ is the angular frequency. 
We assume that for each real wavenumber, there exists a corresponding frequency $\omega(k)$ determined by the PDE in question. Evidently, the above solution travels at a speed $c(k) = \omega(k) / k$.
The frequency can thus be written in terms of the wavenumber as $\omega(k) = c(k) k$; this is called the dispersion relation of the PDE. The group velocity of the solution (\ref{eq:sol}) is the rate of change of $\omega(k)$ with respect to $k$, i.e., $c_{g}(k) =\dif\omega(k) / \dif k$. 
If we consider the non-dispersive linear transport equation, 
\beq\label{eq:oww}
u_t+cu_x=0\,,
\eeq
it is a simple matter to see that the phase speed $c(k)$ is simply constant, $c(k)\equiv c$, and the dispersion relation is just a linear function of $k$, namely, $\omega(k) = c k$. In this case the group velocity is also constant, $c_{g}(k)=c$. Thus, with the advantage of constant phase velocity and group velocity, the PDE \eqref{eq:oww} serves as a convenient platform for analyzing numerical dispersion-and-dissipation errors, and thus for assessing the performance of a FD scheme.

Dispersion-and-dissipation properties of FD schemes depend on their numerical phase and group velocities which depend on their effective (numerical) wavenumber and angular frequency. The effective wavenumber can be obtained by taking the Fourier transform of the spatial FD stencil. To see this, suppose that the spatial derivative of a function $f$ at some point $x$ on a uniform grid has been approximated using an explicit FD stencil with $N$ nodes to the left of $x$ and $M$ nodes to the right of $x$ by
\begin{equation}\label{eq:sten}
\frac{\partial f}{\partial x}(x) \simeq \frac{1}{\Delta x}\sum_{j=-N}^{M}a_j f(x+j \Delta x)\,.
\end{equation}
Taking the Fourier transform of \eqref{eq:sten}, we find
\begin{equation}\label{eq:k}
\mi k\hat{f} \simeq \left(\frac{1}{\Delta x}\sum_{j=-N}^{M}a_j\me^{\mi jk\Delta x}\right)\hat{f}\,,
\end{equation}
where $\hat{f}$ denotes the Fourier transform of $f$. Comparing the left- and right-hand sides of \eqref{eq:k}, we define the \emph{effective wavenumber} for the FD scheme in \eqref{eq:sten} by
\begin{equation}\label{eq:kstar}
k^* = \frac{-\mi}{\Delta x}\sum_{j=-N}^{M}a_j\me^{\mi jk\Delta x}\,,
\end{equation}
so that the harmonic wave $u^*(x,t)=\me^{\mi(kx-\omega^* t)} $, where $\omega^* = c k^*$ (if the time integration is exact), is an approximation to the solution of the one-dimensional scalar transport equation.
We note that the approximation and the exact solution agree initially, that is $u^*(x,0)=u(x, 0) =\me^{\mi kx}$.  However, at any finite time $t=T$, the approximate solution can be written as
\begin{equation}\label{eq:errors}
u^*(x, T) 
= \me^{\mi(kx-\omega^* T)} 
= \me^{\mi(kx-ck^* T + ckT  - ckT)}  
= \underbrace{\me^{\mi c(k-\Re k^*)T}}_{\text{dispersion}}
\cdot \underbrace{\me^{-c\Im k^*T}}_{\text{dissipation}}
\cdot \underbrace{\me^{\mi(kx-\omega T)}}_{\text{exact solution}}\,,
\end{equation}
where $k^*=\Re k^*+\mi\Im k^*$. It is clear from the above back-of-the-envelope calculation that when symmetric central difference FD schemes ($N=M$ and $a_{-j} = -a_{j}$) are considered, $k^*$ is real and the dissipation factor is equal to 1. On the contrary, if an unsymmetric stencil is used, the effective wavenumber $k^*$ is complex for each real $k$.  In this case, $c\Im k^*$ must be positive to prevent $u^*$ from growing exponentially in time, which leads to numerical instability. Using \eqref{eq:kstar}, it is straightforward to show that, for the simplest case, the condition $c\Im k^*>0$ implies either $c a_{-1} < 0$ for $c>0$, or  $c a_{1}  < 0$ for $c<0$. This finding is consistent with the stability requirement for upwinding schemes; this condition is based on the characteristic direction of the transport equation. In general, upwind numerical schemes have built-in numerical damping. This damping rate can be calculated precisely from \eqref{eq:kstar}, \eqref{eq:errors} once the stencil size and stencil coefficients are known \cite{tam06}. Equation \eqref{eq:errors} also shows that dispersion is inevitable for FD schemes, even though the waves supported by the original PDE are non-dispersive.

Temporal discretization introduces numerical dispersion-and-dissipation in the same way.  Independently of the nature of the spatial discretization, Tam and Webb \cite{tw93, tam06} investigated the dispersion-and-dissipation errors introduced by explicit multi-step time-marching schemes and presented an optimization process to minimize the errors.  The  low-dissipation and low-dispersion Runge-Kutta (LDDRK) method developed by Hu et al.\ \cite{hhm96}  minimized the temporal errors of one-step, multi-stage RK methods by quantifying the errors as functions of $c\Delta t k^*$. The analysis was later extended to analyze and regulate dispersion-and-dissipation effects of the optimized RK methods for spatially discretized equations \cite{ad11, bp09, p06, sh98}.

Since there is no obvious benefit to optimizing temporal and spatial discretization at the same time, in this paper, we focus on optimizing the spatial dispersion-and-dissipation error, but incorporate the phase error of the time integrator into our overall error quantification. We show that this error can be controlled to a range that is much smaller than the spatial dispersion-and-dissipation error.  We also discuss the error of numerical group velocity and its relation with the spatial dispersion-and-dissipation error. 
Our error analysis begins with the analysis of compact FD schemes for the transport equation with constant coefficients and with periodic boundary conditions.  
We show that, with some modifications, the analysis can be extended to problems with non-constant coefficients as well.  We also discuss the difficulties of extending the analysis to problems with non-periodic boundary conditions.  Based on our analysis, we develop an optimized compact DRB scheme in periodic domains and a composite boundary scheme for general boundary conditions. 
%
We show that our algorithms give spectral-like resolution. We also show that the proposed algorithms provide accurate results for a nonlinear problem and a higher-dimensional problem.  
Furthermore, our analysis and numerical experiments confirm that under the same computational cost and feasible computer power the optimized compact schemes obtained by sacrificing formal order of accuracy to reduce the dispersion-and-dissipation errors are advantageous for solutions with high modified wavenumber in nature.


\subsection{Plan} 
To aid the reader, we outline the contents of the remainder of this paper as follows. 
In  \S\ref{sec:analysis}, using the constant-coefficient transport equation as a model, we introduce the error analysis for optimized compact DRB schemes in a periodic domain. We give a quantitative description of numerical phase error and $L^2$-norm error. We use examples to illustrate the effects of different optimization strategies.  
In \S\ref{sec:DRB-design}, based on the error dynamics of DRB schemes, we propose a compact DRB scheme that defines the numerical effective modified wavenumber \emph{before} finding the corresponding discretization. This approach allows us to enlarge the number of degrees of freedom for minimizing the dispersion-and-dissipation error, and we use all degrees of freedom to minimize the dispersion-and-dissipation errors (rather than eliminating the Taylor truncation errors). 
We extend the analysis and the algorithm to the transport equation with spatially various coefficients in  \S\ref{sec:non-constant}. We show that the aliasing error for the non-constant coefficient problems is at least as important as the dispersion-and-dissipation error. 
In \S\ref{sec:nonlinear} and \S\ref{sec:2D}, we treat both a nonlinear and higher-dimensional examples. Specifically, we analyze the performance of compact DRB schemes for the nonlinear Hopf equation at a time before wave breaking when the spectrum swiftly moves to the large-wavenumber regime. We find that the benefits of DRB schemes carry over to this nonlinear example. 
We then extend our error analysis to higher-dimensional case, and we illustrate the advantages of our proposed scheme by means of a two-dimensional test problem, Zalesak's rotating disk. 
We also discuss the difficulties of extending the minimization process to problems with non-periodic boundary conditions and propose a hybrid composite compact DRB scheme as an alternative in \S\ref{sec:non-periodic}. Finally, we make some concluding remarks in \S\ref{sec:conclusion}.

\section{Analysis of compact DRB schemes}\label{sec:analysis}

\subsection{Beyond Taylor truncation error}
After spatial discretization, a given time-evolution PDE becomes a semi-discretized system of ordinary differential equations (ODEs) that may be written as 
\begin{equation}\label{2_1}
\frac{\dif}{\dif t}\textbf{U}=F(\textbf{U})\,,
\end{equation}
where $\textbf{U}$ is a vector-valued function of $t$.  We suppose that the function $F$ satisfies a Lipschitz condition 
\begin{equation}\label{2_2}
\Vert F(\textbf{U})-F(\textbf{V})\Vert\leq\mathscr{L}\Vert\textbf{U}-\textbf{V}\Vert\,,
\end{equation}
where  $\mathscr{L}$ is the Lipschitz constant. This constant plays a role in the uniqueness and stability of the ODE system and in the error estimate of its numerical solution.

If we numerically integrate the ODE system \eqref{2_1}---without imposing any stability constraints on the mesh sizes $\Delta x$ and $\Delta t$--- by means of some temporal scheme from $t=0$ to $t=T$, the $L^{\infty}$-norm error associated with the numerical solution at time $T$, denoted by $e(T)$, is bounded by
\begin{equation}\label{2_3}
e(T)\leq T \me^{\mathscr{L}T}\left(\Delta x^s+\Delta t^p\right),
\end{equation}
where $s$ and $p$ are the orders of local truncation errors for spatial and temporal discretization, respectively.
It is clear that even for modest sizes of $\mathscr{L}$ and $T$, this bound can be incredibly large due to the exponential factor $\me^{\mathscr{L}T}$. Indeed, it is known that this bound is not some theoretical  overestimate; the error of such numerical solutions can be unacceptably large.

After imposing suitable stability constraints, we obtain a more useful error estimate:
\begin{equation}\label{2_4}
e(T)\leq C_T\left(\Delta x^s+\Delta t^p\right),
\end{equation}
where $C_T$ is a constant depending on $T$.  The new error bound (\ref{2_4}) gives a much more practical relationship between mesh sizes and corresponding numerical errors. Nevertheless, the magnitude of $C_T$ can still be so large that in practice obtaining accurate numerical solutions is at the expense of computational cost (i.e., small $\Delta x$ or $\Delta t$ or both). 
More specifically, the cost enters in two ways: (1) small $\Delta x$ means more memory is required, and (2) small $\Delta x$ and $\Delta t$ result in a large CPU/GPU time. While the latter issue may be addressed by adopting parallel computing techniques, memory capacity is and will continue to be a problem for large scale simulations. In order to better utilize the available computer power, a more specific relationship between the adopted mesh sizes and the smallest scale of the resolvable features must be established. Optimized DRB schemes are developed exactly for attacking this issue.

In order to design a FD scheme that is capable of providing accurate numerical solutions with a smaller memory requirement and less CPU/GPU time, various optimization methodologies, based on the idea of preserving the dispersion relation, have been proposed. For a numerical solution, achieving the desired accuracy means that the associated error is less than a prescribed error tolerance. The question we must ask is whether sacrificing formal order of accuracy in order to preserve the dispersion relation is more advantageous than keeping the higher-order truncation errors without optimization, especially if we assume that all methods have the same computational cost. 
To answer this question properly, we must assess four factors: (1) error tolerance, (2) maximal available memory capacity, (3) the highest wavenumber we want to resolve, and (4) the actual CPU/GPU time for a simulation. If any of these four factors is not properly addressed, any conclusion in a comparison study may be misleading and subject to question. In the following sections, we provide several examples to illustrate these points.  

We remark that for most (if not all) optimized DRB schemes, the optimization procedures and the underlying analysis are developed based on Fourier analysis. Consequently, the $L^2$-norm becomes a natural choice for measuring the errors. As a spatial-average, the $L^2$ error provides a better framework for problems originating from data assimilation and certain inverse problems. For the error estimates derived based on the Taylor truncation error (for higher-order, unoptimized schemes), we refer the reader to the text of LeVeque \cite{LeVeque07}.

\subsection{Non-dispersive scalar transport equation}
Consider the one-dimensional scalar transport equation
\begin{equation}\label{2_5}
u_t+cu_x=0\,,\quad 0\leq x\leq\ell\,,
\end{equation}
with periodic boundary conditions. Unless stated otherwise, we assume that $c>0$ in the following discussion. For a prescribed initial condition $u(x,0)=u_0(x)$, the true solution at time $T$ is known to be $u(x,T)=u_0(x-cT)$. The initial profile simply travels to the right at a constant speed $c$ without being deformed. Suppose that the solution can also be represented in terms of the Fourier series
\begin{equation}\label{2_7}
u(x,t)=\sum_{k'=-\infty}^{\infty}\widehat{u}(k',t)\me^{2\pi \mi k'x/\ell}
\end{equation}
where
\begin{equation}\label{2_6}
\widehat{u}(k',t)\equiv\dfrac{1}{\ell}\int_0^{\ell} u(x,t)\me^{-2\pi \mi k'x/\ell}\,\dif x,\quad k'\in\mathbb{Z}\,,
\end{equation}
denotes the $k'$-th Fourier coefficient.
Substituting the solution ansatz \eqref{2_7} into the PDE and applying the initial condition, we solve for $\widehat{u}(k',t)$, and we obtain
\begin{equation}\label{2_8}
\widehat{u}(k',t)=\widehat{u}_0(k')\me^{-2\pi \mi ck't/\ell}\,,\quad k'\in\mathbb{Z}\,,
\end{equation}
where $\widehat{u}_0(k')$ is the $k'$-th Fourier coefficient of $u_0(x)$. 

\begin{note}
Evidently, the above procedure is equivalent to taking the Fourier transform of \eqref{2_5} and its initial condition, and then solving the resulting ODE initial-value problem (IVP).  The solution is exactly 
\begin{equation}\label{2_9}
u(x,t)=\sum_{k'=-\infty}^{\infty}\widehat{u}_0(k')\me^{2\pi \mi k'(x-ct)/\ell}\,.
\end{equation}
We remark that each component $\me^{2\pi \mi k'x/\ell}$ in the solution ansatz (\ref{2_7}) is an eigenfunction of the eigenvalue problem associated with the differential operator $\dif/\dif x$ in the periodic domain, so that the solution representation \eqref{2_9} is indeed an eigenfunction expansion. 
\end{note}

When the spatial derivative is approximated by a FD scheme, the solution of the resulting semi-discretized ODE (also an IVP) can be represented by a Fourier series as well (shown later, cf. \eqref{2_38}). The difference between these two Fourier series is that instead of  $\me^{2\pi \mi k'(x-ct)/\ell}$, the eigenfunction of the approximate solution is  $\me^{2\pi \mi(kx_j-\omega^*t)/\ell}$, which suggests that the wave travels at an approximate speed $c^*\equiv\omega^*/k$, rather than at the true speed $c$. Our analysis in the following sections aims to use Fourier representations to quantify the numerical errors associated with compact FD schemes. This analysis parallels the classical von Neumann analysis, but we derive new error quantification criteria to aid the future development of more accurate FD schemes.


\subsection{Compact FD discretization}\label{sec:compact}

Consider an $\ell$-periodic function $u(x)$ on a uniform mesh $x_j=j\Delta x$, $\Delta x=\ell/N$, $1\leq j\leq N$. Let $U_j=u(x_j)$. A generalized Pad\'e compact FD scheme \cite{lele92} to approximate the derivative $u'(x_j)$ is 
\begin{equation}\label{2_10}
\beta U'_{j-2}+\alpha U'_{j-1}+U'_{j}+\alpha U'_{j+1}+\beta U'_{j+2}=a\dfrac{U_{j+1}-U_{j-1}}{2\Delta x}+b\dfrac{U_{j+2}-U_{j-2}}{4\Delta x}+c\dfrac{U_{j+3}-U_{j-3}}{6\Delta x}\,,
\end{equation}
where $U'_j\approx u'(x_j)$. The coefficients $a$, $b$, $c$ and $\alpha$, $\beta$ are determined by matching the desired order. The order of the approximation (\ref{2_10}) is determined by the following constraints:
\begin{subequations}\label{eq:padeconstraints}
\begin{align}
a+b+c & = 1+2\alpha+2\beta \,\,\,\, \qquad\qquad\text{(second order)}\\
a+2^2b+3^2c &= 2\cdot 3(\alpha+2^2\beta) \qquad\qquad\text{(fourth order)}\\
a+2^4b+3^4c &= 2\cdot 5(\alpha+2^4\beta) \qquad\qquad\text{(sixth order)}\\
a+2^6b+3^6c &= 2\cdot 7(\alpha+2^6\beta) \qquad\qquad\text{(eighth order)}\\
a+2^8b+3^8c &= 2\cdot 9(\alpha+2^8\beta) \qquad\qquad\text{(tenth order)}
\end{align}
\end{subequations}
We introduce the matrix-vector notations 
$$
\textbf{U}=\left(U_1,\:U_2,\cdots,\:U_N\right)^T
\quad
\text{and} 
\quad
\textbf{U}'=\left(U'_1,\:U'_2,\cdots,\:U'_N\right)^T\,,
$$ 
so that \eqref{2_10} can be expressed compactly as
\begin{equation}\label{2_16}
\textbf{L}\textbf{U}'=\textbf{R}\textbf{U}\,,
\end{equation}
where $\textbf{L}$ and $\textbf{R}$ are cyclic pentadiagonal and cyclic heptadiagonal matrices, respectively, whose nonzero entries can be computed from  \eqref{2_10}. We rewrite \eqref{2_16} as 
\[
\textbf{U}'=\textbf{L}^{-1}\textbf{R}\textbf{U}\,,
\]
and we let $\textbf{D}\equiv\textbf{L}^{-1}\textbf{R}$. This matrix $\textbf{D}$ is the discrete analogue of the differential operator $\dif/\dif x$. 

\begin{note}
We remark that the classical centered difference formula is a special case of the compact FD scheme with $\alpha=\beta=0$. For nonzero $\alpha$ and $\beta$, it is necessary to solve the pentadiagonal (or tridiagonal if $\beta=0$) system to obtain $\textbf{U}'$. However, the higher computational cost is compensated by higher order of accuracy and smaller truncation error. Furthermore, compact schemes require less stencil points for achieving a desired order of accuracy compared with explicit schemes. This feature may be advantageous when dealing with stability issues.
\end{note}

Due to the spatial homogeneity of coefficients and periodicity, the eigenvectors of $\textbf{D}$ have components 
\begin{equation}\label{2_17}
U_j=\me^{2\pi \mi k'x_j/\ell}\,,\quad j=1,\ldots,N\,,
\end{equation}
where $k'$ ranges from $-N/2+1$ to $N/2$. For simplicity, we assume that $N$ is even; a similar formulation can be derived for odd $N$. For each particular $k'$, $k=(2\pi/\ell)k'$ is the angular wavenumber, or simply wavenumber, and the components of the corresponding eigenvector can be expressed as
\begin{equation}\label{2_18}
U_j=\me^{\mi kx_j}\,,\quad j=1,\ldots,N\,.
\end{equation}
For a FD scheme, the corresponding $U_j'$ in \eqref{2_16} for the above eigenvector can be represented as 
\begin{equation}\label{2_19}
U'_j=\mi k^*\me^{\mi kx_j},\quad j=1,\ldots,N\,.
\end{equation}
Substituting \eqref{2_18} and (\ref{2_19}) into \eqref{2_10}, we find an expression for the effective (numerical) wavenumber $k^*$ 
\begin{equation}\label{2_20}
\kappa^*=\dfrac{a\sin(\kappa)+(b/2)\sin(2\kappa)+(c/3)\sin(3\kappa)}{1+2\alpha\cos(\kappa)+2\beta\cos(2\kappa)}\,.
\end{equation}
Here $\kappa\equiv k\Delta x$ and $\kappa^*\equiv k^*\Delta x$ are known as the modified wavenumber and modified effective wavenumber, respectively. We refer to \eqref{2_20} as the {\em resolution characteristic} of the spatial discretization.
The maximum of $\vert k\vert$ is specified by the Nyquist limit $\pi/\Delta x$, and is evidently inversely proportional to the spatial mesh size $\Delta x$. The modified wavenumber $\kappa$, however, lies in the fixed range $\left[-\pi,\pi\right]$. The relationship between $\kappa^*$ and $\kappa$ is fully determined by the coefficients of the spatial discretization and is independent of $\Delta x$ and the domain size $\ell$.  We note that because $\kappa^*$ is an odd function in $\kappa$, it suffices to consider only nonnegative $\kappa$. 
 
The best numerical approximation to $u(x)=\me^{\mi kx}$ a numerical method can achieve is given by the spectral method, for which $U_j'=\mi k^*\me^{\mi kx_j}$, where $k^{*}\equiv k$. For FD schemes, the error of the approximation of the first derivative is quantified by the difference $\vert k^*-k\vert$ (or, equivalently, $\vert\kappa^*-\kappa\vert$).
For example, for a FD scheme whose derivation is solely based on Taylor truncation error, it is straightforward to show, using Taylor expansion, that for the approximation of spatial derivative of order $s$, the relation between $\kappa^*$ and $\kappa$ is
\begin{equation}\label{2_21}
\kappa^* = \kappa + c_{s+1} \kappa^{s+1} + O(\kappa^{s+2}),
\end{equation}
for some constant $c_{s+1}$. This indicates that $\kappa^*$ is only a good approximation to $\kappa$ within a relatively small region. The idea of DRB schemes is to extend the region of good approximation of $\kappa$ to one that is as large as possible.

Returning now to our model problem, the transport equation \eqref{2_5}, discretization in space yields an ODE system
\begin{equation}\label{2_22}
\textbf{U}_t+c\textbf{D}\textbf{U}=0.
\end{equation}
Recalling the Discrete Fourier Transform (DFT) in space, 
\begin{equation}\label{2_24}
\widehat{\textbf{U}}(k',t)\equiv\dfrac{1}{N}\sum_{j=1}^{N}U_j(t)\me^{-2\pi \mi k'x_j/\ell}\,,
\end{equation}
we note that the solution of \eqref{2_22} on the grid point $j$ at some time $T$ is 
\begin{equation}\label{2_26}
U_j(T)=\sum_{k'=-N/2+1}^{N/2}\widehat{\textbf{U}}(k',0)\me^{\mi(kx_j-ck^*T)}\,,\quad1\leq j\leq N\,.
\end{equation}
Suppose the initial condition is specified on a mesh by $U_j(0)=u_0(x_j)$, $1\leq j\leq N$, where $u_0(x)$ is a function defined on the continuous interval $[0,\ell]$, then we can represent $\widehat{\textbf{U}}(k',0)$ by
\begin{equation}\label{2_27}
\widehat{\textbf{U}}(k',0)=\dfrac{1}{N}\sum_{\alpha\equiv k'\pmod N}\widehat{u}_0(\alpha)\,,\quad
-\dfrac{N}{2}+1\leq k'\leq \dfrac{N}{2},
\end{equation}
where  $\widehat{u}_0(\alpha)$ denotes the Fourier coefficients of $u_0(x)$. Since $N<\infty$ is finite, there always exists aliasing error (where $\widehat{\textbf{U}}(\alpha,0)$ for $\alpha > N/2$, or $\alpha<-N/2$ does not exist and is represented by  $\widehat{\textbf{U}}(\alpha',0)$, for some $\alpha'$,  $-N/2<\alpha'<N/2$). We assume that in our numerical experiments $N$ is large enough (or $\Delta x$ is small enough), so that the aliasing error only occurs in the high wavenumber that are beyond the range in which we are interested.

In the following sections, to illustrate numerical phase and group velocity errors due to the time integrator, we solve the above semi-discretized system by an explicit $M$-stage, $p$-th order RK method. The numerical solution of the this time integrator approximates the matrix exponential of the exact solution by a polynomial.  We express the numerical phase and group velocity errors in terms of the spatial dispersion-and-dissipation error plus a higher-order term, and show that this higher-order term can made small compared with the spatial dispersion-and-dissipation error. This approach provides a framework of analysis that can be generalized to linear systems with non-constant coefficients, systems with non-periodic boundary conditions,  and to multi-dimensional systems.

\subsection{Numerical solution of a RK temporal integration}\label{sec:RK}
The exact solution of \eqref{2_22} at $t=T+\Delta t$ can be written as
\begin{equation}\label{2_28}
\textbf{U}(T+\Delta t)=\me^{-c\Delta t\textbf{D}}\textbf{U}(T)\,.
\end{equation}
where the matrix exponential is defined to be
\begin{equation}\label{2_29}
\me^{-c\Delta t\textbf{D}}\equiv\sum_{m=0}^{\infty}\dfrac{1}{m!}(-c\Delta t\textbf{D})^m.
\end{equation}
We write $U_j^n$ for the FD approximation to $u(x_j,t_n)$, $x_j=j\Delta x$, $t_n=n\Delta t$, and $\textbf{U}^n=\left(U_1^n,\:U_2^n,\cdots,\:U_N^n\right)^T$. An explicit $M$-stage RK integrator that advances the solution from $t_n$ to $t_{n+1}$ is
\begin{equation}\label{2_30}
\begin{split}
\textbf{K}_1&=-c\Delta t\textbf{D}\textbf{U}^n\,,\\
\textbf{K}_i&=-c\Delta t\textbf{D}(\textbf{U}^n + w_{i-1}\textbf{K}_{i-1}),\quad 2\leq i\leq M\,,\\
\textbf{U}^{n+1}&=\textbf{U}^n + w_M\textbf{K}_M\,.
\end{split}
\end{equation}
Effectively, an $M$-stage $p$-th order (or $(p+1)$-th order in terms of local truncation error) RK integrator like \eqref{2_30} approximates the matrix exponential $\me^{-c\Delta t\textbf{D}}$ by a polynomial of degree $M$:
\begin{equation}\label{2_31}
\mathcal{P}(-c\Delta t\textbf{D})\equiv\sum_{m=0}^M a_m(-c\Delta t\textbf{D})^m\,,
\end{equation}
with the first $p+1$ coefficients subject to conditions $a_m=1/m!$, $0\leq m\leq p$. The $w_i$'s and $a_m$'s are related by
\begin{equation}\label{2_32}
\begin{split}
a_0&\equiv 1\,,\\
a_1&=w_M\,,\\
a_2&=w_M w_{M-1}\,,\\
&\vdots\\
a_M&=w_M w_{M-1}\cdots w_1\,.
\end{split}
\end{equation}
Hence the time step in \eqref{2_30} can be rewritten as
\begin{equation}\label{2_33}
\textbf{U}^{n+1}=\mathcal{P}(-c\Delta t\textbf{D})\textbf{U}^n\,.
\end{equation}
The implementation of (\ref{2_30}) requires only two levels of storage, one for $\textbf{U}^n$ and one for $\textbf{K}_i$ \cite{hhm96}. Note that if the coefficients of the original PDE are spatially inhomogeneous, or more complicated (linear) boundary conditions are imposed, we only need to modify the matrix $\textbf{D}$ in the above formulation. That is, the time advancing scheme (\ref{2_30}) with coefficients $w_i$'s,  determined by \eqref{2_32}, works for all linear cases. We note that Butcher \cite{butcher64} has shown that, for general nonlinear problems, it is impossible to reduce the storage and computational cost of the RK schemes except for certain four-stage schemes.

Taking the DFT in space of \eqref{2_33}, we obtain
\begin{equation}\label{2_34}
\widehat{\textbf{U}}^{n+1}(k')=\mathcal{P}(-\mi ck^*\Delta t)\widehat{\textbf{U}}^n(k')\,,\quad -\dfrac{N}{2}+1\leq k'\leq \dfrac{N}{2}\,.
\end{equation}
We now introduce the numerical temporal frequency
\begin{equation}\label{2_35}
\omega^*\equiv\dfrac{\mi}{\Delta t}\log\mathcal{P}(-\mi ck^*\Delta t)\,.
\end{equation}
Since the polynomial $\mathcal{P}$ approximates the exponential up to order $p$, we have
\begin{equation}\label{2_36}
\log\mathcal{P}(-\mi ck^*\Delta t)=-\mi ck^*\Delta t + O((ck^*\Delta t)^{p+1})\,,
\end{equation}
and
\begin{equation}\label{2_37}
\omega^*=ck^* + O((ck^*)^{p+1}(\Delta t)^p)\,.
\end{equation}
As a result, it is clear that \eqref{2_37} specifies the dispersion relation of the overall FD scheme. With this definition of $\omega^*$, we write \eqref{2_34} as $\widehat{\textbf{U}}^{n+1}(k')=\widehat{\textbf{U}}^n(k')\me^{-\mi\omega^*\Delta t}$, and we can explicitly write down the exact solution associated with the overall FD scheme:
\begin{equation}\label{2_38}
U_j^n=\sum_{k'=-N/2+1}^{N/2}\widehat{\textbf{U}}(k',0)\me^{\mi(kx_j-\omega^*t_n)}\,,
\quad 1\leq j\leq N\,.
\end{equation}
Recall that $k=(2\pi/\ell)k'$, and the relation between $\kappa\equiv k\Delta x$ and $\kappa^*\equiv k^*\Delta x$ is given in \eqref{2_20}. Since
$\omega^*$ depends on $k^*$ and $\Delta t$, it is a function of $k$, $\Delta x$ and $\Delta t$.  

\subsection{Phase error and $L^2$-norm error}
Let $\omega^*\Delta t$ be the effective (numerical) modified temporal frequency. Using \eqref{2_37}, we can express the dispersion relation \eqref{2_35} as a function of $\kappa^*$ times a dimensionless number
\begin{equation}\label{2_39}
\omega^*\Delta t = r\kappa^* + O((r\kappa^*)^{p+1})\,.
\end{equation}
Here $r=c\Delta t/\Delta x$ is the Courant-Friedrichs-Lewy (CFL) number for the transport equation. While $\kappa^*$ only depends on $\kappa$, the quantity $\omega^*\Delta t$ depends on both $\kappa$ and $r$. The behavior of the fully discretized problem is completely characterized by $\omega^*\Delta t$ and $\kappa^*$.

Suppose that we are interested in resolving a wave solution whose length scale is specified by the wavenumber $k$. From \eqref{2_38}, the numerical phase speed $c^*$ associated with this particular wave component is given by
\begin{equation}\label{2_40}
c^*=\dfrac{\omega^*}{k}\\
=\dfrac{\Delta x}{\Delta t}\dfrac{\omega^*\Delta t}{k\Delta x}
=\dfrac{c}{r}\dfrac{r\kappa^*+O((r\kappa^*)^{p+1})}{\kappa}
=c\dfrac{\kappa^*}{\kappa}\left(1+O((r\kappa^*)^{p})\right)\,,
\end{equation}
and hence the numerical phase error $\theta_e$ at time $T=n\Delta t$ is 
\begin{equation}\label{2_41}
\theta_e=\omega^*T-ckT
=(c^*-c)kT
=ckT\dfrac{\kappa^*-\kappa}{\kappa} + ckT\dfrac{\kappa^*}{\kappa}O((r\kappa^*)^{p})\,.
\end{equation}
The first term in the last expression of $\theta_e$ is mainly due to the spatial discretization, whereas the second term is  largely from the temporal integration. The latter term can be effectively suppressed by adopting a small CFL number (equivalently, small $\Delta t$ compared with $\Delta x$). It is interesting to see that in terms of phase error, the error introduced by spatial discretization has a first-order effect, while the error caused by the temporal integration is negligible if a sufficiently small $r$ is used.  Therefore, when making the comparison of phase-accuracy for various spatial discretizations, we must make sure that the CFL number is sufficiently small so that the comparison is not influenced by the error introduced by the temporal scheme. 

To see the relation between the overall numerical error and the phase error, we compare the numerical solution (\ref{2_38}) with the true solution. The error at time $T=n\Delta t$, denoted by $e(T)$, in the $L^2$-norm is
\begin{equation}\label{2_43}
\begin{split}
\Vert e(T)\Vert_{L^2}& = \sqrt{\sum_{j=1}^N\vert u(x_j,T)-U_j(T)\vert^2\Delta x}
=\sqrt{\ell\sum_{k'=-N/2+1}^{N/2}\vert\widehat{U}(k',0)\vert^2\vert \me^{\mi\theta_e}-1\vert^2}\\
&\approx\sqrt{\ell\sum_{k'=-N/2+1}^{N/2}\vert\widehat{U}(k',0)\vert^2 |\theta_e|^2}\,.
\end{split}
\end{equation}
The second equality holds thanks to \eqref{2_27} and Parseval's identity. The last approximation assumes that $|\theta_e|$ is small so that $\me^{\mi\theta_e}\approx 1+\mi\theta_e$. Equation \eqref{2_43} shows that numerical error in the $L^2$-norm depends on the phase error $\theta_e$. 
Now,  from \eqref{2_41}, in the limiting case $r\rightarrow 0^+$, we have simply
\begin{equation}\label{2_44}
\theta_e\approx ckT\dfrac{\kappa^*-\kappa}{\kappa}\,,
\end{equation}
where $k=\dfrac{2\pi}{\ell}k'$, $\kappa=k\Delta x$. 
Substituting \eqref{2_44} into the last expression in \eqref{2_43}, we see that
\begin{multline}\label{2_45}
\Vert e(T)\Vert_{L^2}\approx cT\sqrt{\ell\sum_{k'=-N/2+1}^{N/2}\vert\widehat{U}(k',0)\vert^2k^2\left|\frac{\kappa^*-\kappa}{\kappa}\right|^2}\\
=\dfrac{cT}{\Delta x}\sqrt{\ell\sum_{k'=-N/2+1}^{N/2}\vert\widehat{U}(k',0)\vert^2\vert\kappa^*-\kappa\vert^2}\,,
\end{multline}
provided that the phase error $|\theta_e|$ is sufficiently small for those $k'$ with large nonzero $\widehat{U}(k',0)$. Equation \eqref{2_45} indicates that a good understanding of the resolution characteristic of the spatial discretization, or, more precisely, the error $\vert\kappa^*-\kappa\vert$, is a prerequisite for controlling the numerical error associated with these FD schemes. Furthermore, \eqref{2_45} suggests that it makes good sense to design a spatial discretization based on minimizing $\vert\kappa^*-\kappa\vert$.  
\begin{note}[A word of caution]
Equation \eqref{2_43} indicates that the quantitative description of the $L^2$-norm error is determined by the numerical phase error $\theta_e$ which is, in turn, related to the numerical phase speed $c^*$. However,  we remind the reader that the derivation is obtained from a simple model, namely the non-dispersive transport equation \eqref{2_5} with periodic boundary conditions. For more general problems, the relation between the $L^2$-norm numerical error and the phase error requires further investigation and may be more involved. 
\end{note}

\subsection{Numerical group velocity error} The numerical group velocity $c_g^*$ of a fully discretized PDE by a FD scheme is given by
\begin{equation}\label{2_42}
c_g^*=\dfrac{\partial\omega^*}{\partial k}\\
=\dfrac{\Delta x}{\Delta t}\dfrac{\partial\omega^*\Delta t}{\partial\kappa}\\
=c\dfrac{\partial\kappa^*}{\kappa}\left(1+O((r\kappa^*)^{p})\right).
\end{equation}
One potentially useful property of the numerical group velocity, pointed out by Tam \cite{tam06}, is that numerical group velocity can provide some level of quantitative prediction of the position of the numerical oscillations near sharp gradients or jump discontinuities in a solution wave packet.
As in the analysis of the numerical phase error, the numerical group velocity error associated with $c_g^*$ can be controlled by two terms. One is from spatial discretization, and the other is from temporal integration. For the transport equation \eqref{2_5}, the group velocity is exactly $c$. We remark that in the limiting case $r\rightarrow 0^+$ (or the time integration is exact), $c_g^*=c\partial\kappa^*/\partial\kappa$. This was studied previously by Tam \cite{tam06} and others \cite{sds07, t82}. However, the utility of the information provided by $c_g^*$ appears to be limited in the context of the development of accurate numerical FD algorithms. The basic reason stems from the fact that $\kappa^*$ is a function of $\kappa$. By the fundamental theorem of calculus, therefore, we expect to find an estimate for $\kappa^*$ from knowledge of $\partial\kappa^*/\partial\kappa$. The converse is generally not true. 
However, in practice, we often see that a DRB scheme derived from minimizing  $\vert\kappa^*-\kappa\vert$ also gives a good approximation for the numerical group velocity (or $\partial\kappa^*/\partial\kappa$). This can be explained by the fact that if $\kappa^*$ is a smooth function in $\kappa\in[-\pi,\pi]$, with good a priori knowledge of the signs of $\kappa^*-\kappa$ and $\partial^2\kappa^*/\partial\kappa^2$, and a bound of $\partial^2\kappa^*/\partial\kappa^2$, it is possible to derive an estimate for $\partial\kappa^*/\partial\kappa$ in terms of $\vert\kappa^*-\kappa\vert$. On the other hand, even when we can obtain such an estimate, it is often not general and very much relies on a priori knowledge. Therefore, we will not pursue this point further here. 

\subsection{Resolution characteristic of spatial discretization}
From our previous discussion, it is clear that the error $\vert\kappa^*-\kappa\vert$ plays a crucial role in the error analysis. A detailed knowledge of the behavior of $\kappa^*$ and $\kappa^*-\kappa$, as functions of $\kappa$, can help us better understand the performance of FD schemes from the viewpoint of the dispersion-and-dissipation error. 

We begin with Equation \eqref{2_21}. In order to quantify the range of modified wavenumber $\kappa$ within which the dispersion-and-dissipation error, $\kappa^*-\kappa$, is dominated by the leading-order term $c_{s+1}\kappa^{s+1}$, we define $\gamma>0$ to be the largest number so that if the modified wavenumber $\kappa$ satisfies  $0<\kappa\leq\gamma$, we have
\begin{equation}\label{2_47}
\left|\dfrac{\kappa^*-\kappa-c_{s+1}\kappa^{s+1}}{c_{s+1}\kappa^{s+1}}\right|
\leq\epsilon_{\gamma}\,,
\end{equation}
where $\epsilon_{\gamma}\ll1$ is a chosen threshold measuring the tolerance of the difference between the leading-order error and overall error.  The interval of interest for $\kappa$ is $0\leq\kappa\leq\pi$, and our numerical experiments suggest that $\gamma$ is often $O(1)$. Usually, smaller values of $\epsilon_{\gamma}$ imply a smaller $\gamma$ for a fixed spatial discretization.

\begin{figure}[ht]
\begin{centering}
\includegraphics[scale=0.6]{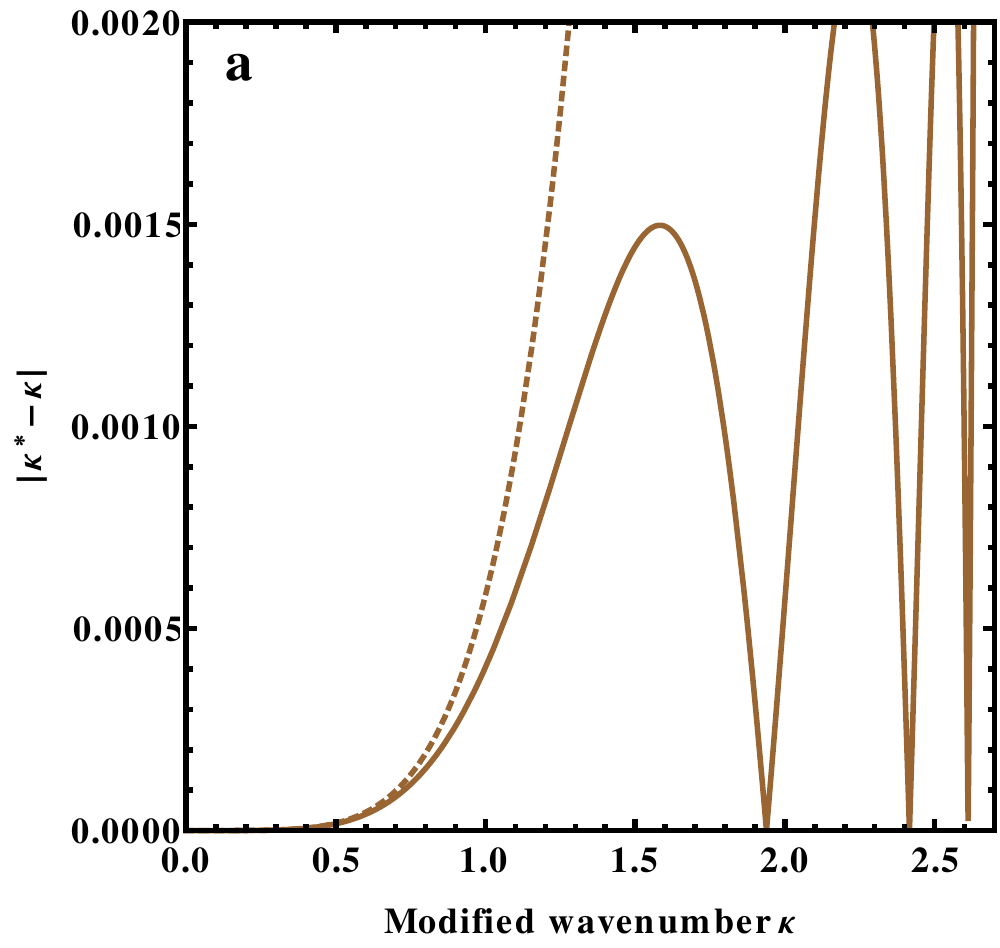}\hskip 0.58in
\includegraphics[scale=0.6]{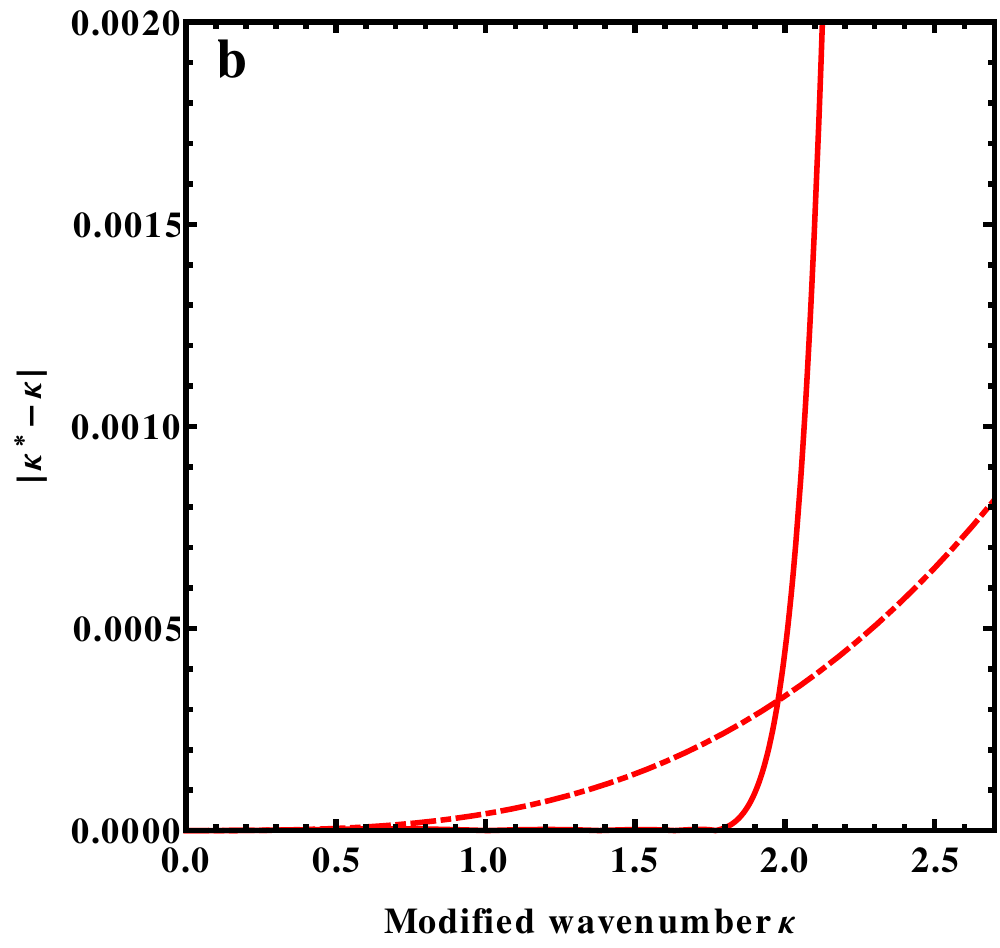}\\
\includegraphics[scale=0.6]{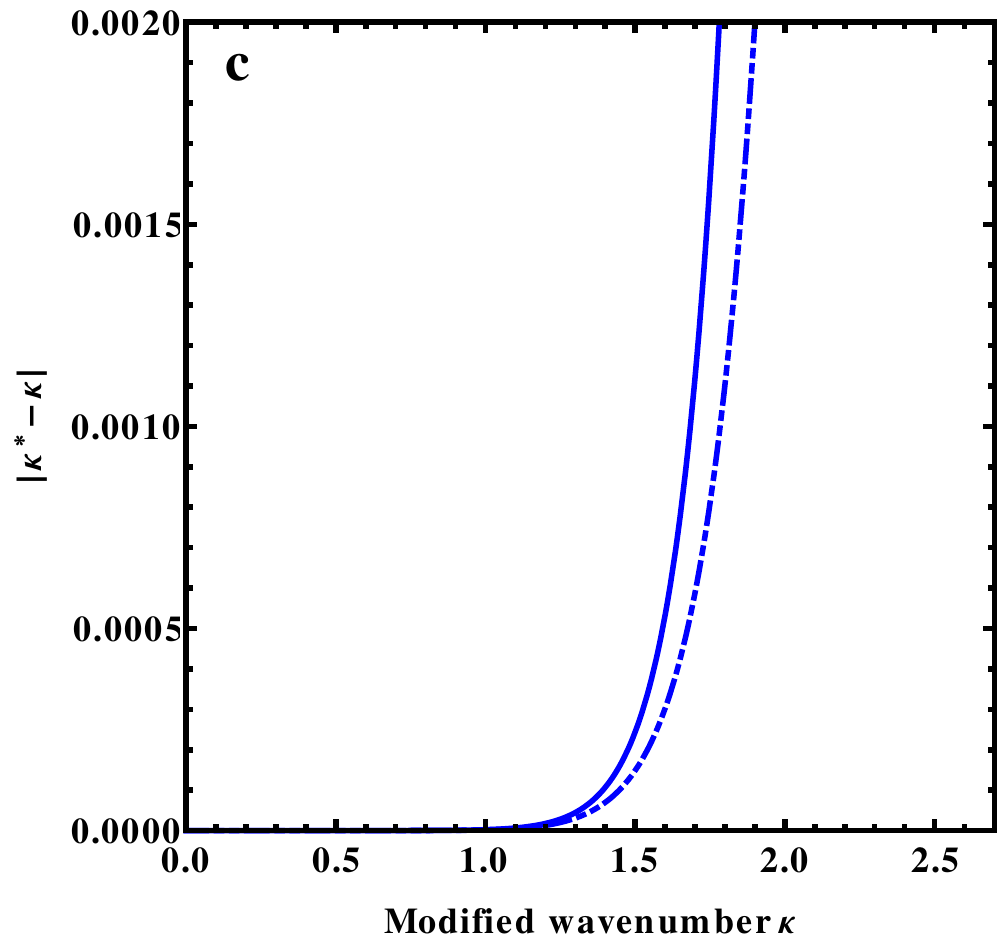}\hskip 0.6in
\includegraphics[scale=0.6]{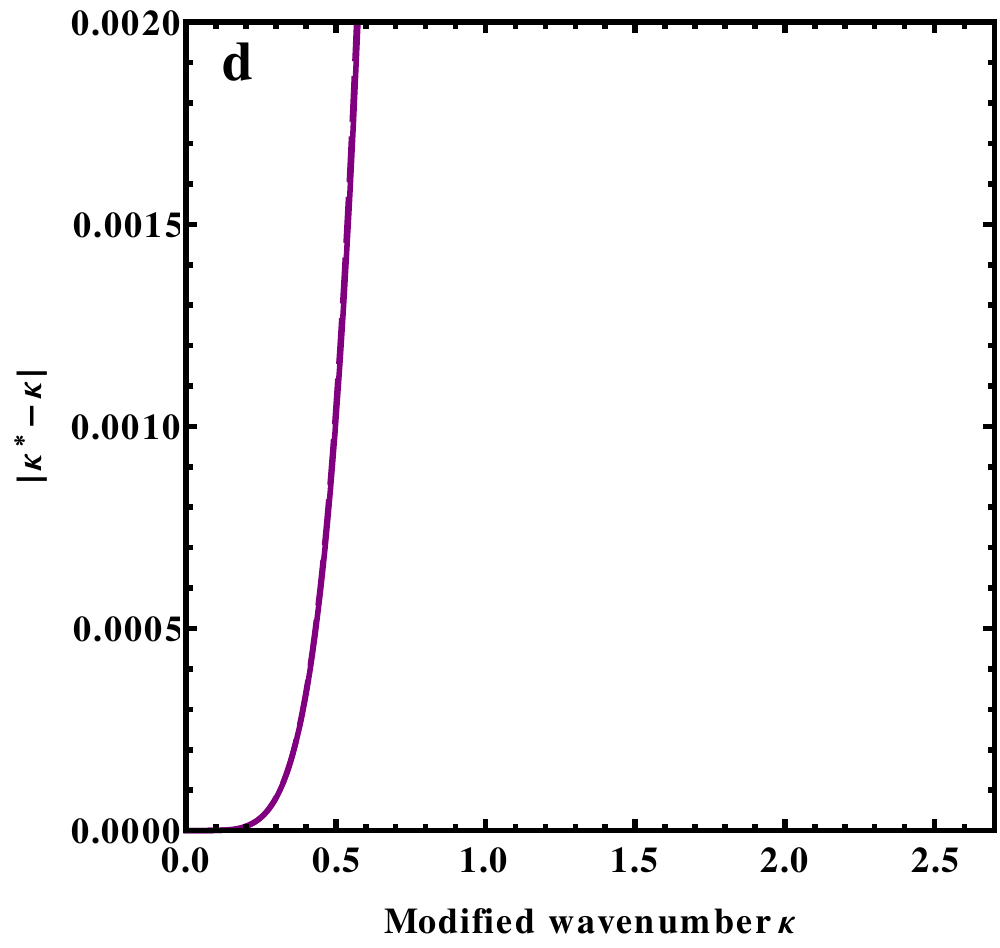}
\end{centering}
\caption{The solid lines are $\vert\kappa^*-\kappa\vert$ and the dashed lines are the leading order error $\vert c_{s+1}\kappa^{s+1}\vert$. The horizontal axis is the modified wavenumber $\kappa$. (a) The fourth-order optimized pentadiagonal scheme by Kim \cite{kim07} (referred to as ``KIM4TH''); (b) The optimized second-order pentadiagonal scheme developed in this paper (referred to as ``OPT2ND1P8''); (c) The tenth-order unoptimized pentadiagonal compact scheme in \cite{lele92} (referred to as ``UNOPT10TH''); (d) The explicit fourth-order central differencing scheme (referred to as ``CD4'').}
\label{Fig:2_2_2}
\end{figure}

When $s\ge 1$ and $\kappa$ is sufficiently small, the leading-order error is a good approximation to the dispersion-and-dissipation error. However this does not imply that the error itself is small, because of the factor $c_{s+1}$. For example, in Figure  \ref{Fig:2_2_2}(d), we show that for the explicit fourth-order central differencing scheme (referred to as ``CD4''), the leading-order error is a good approximation to the dispersion-and-dissipation error. The error, however, is larger than that of all other schemes in comparison for $\kappa\le 0.5$

Suppose $\kappa$ lies in $[0,\gamma]$ (the dispersion-and-dissipation error is well approximated by the leading-order term). Since $\kappa=k\Delta x$, decreasing $\kappa$ (equivalent to decreasing $\Delta x$ for a fixed $k$) will decrease the dispersion-and-dissipation error. Hence for a given FD scheme, if we can determine $\gamma$, then we can determine a suitable mesh size to limit the dispersion error.  In particular for problems whose solutions are dominated by a particular wavenumber $k$.

For unoptimized FD schemes derived solely from Taylor series expansion, the value of $\gamma$ completely characterizes its resolution properties. For such schemes, $|\kappa^*-\kappa|$ is a monotonic increasing function of $\kappa$, as shown in Figure \ref{Fig:2_2_2}(c) \& (d) for the tenth-order unoptimized pentadiagonal compact scheme (referred to as ``UNOPT10TH'') in \cite{lele92}  and CD4, respectively.  Consequently, for this type of schemes decreasing  $\kappa$ is the only way to reduce the dispersion-dissipation error $\vert\kappa^*-\kappa\vert$.
%

Under a more general setting, let us consider an initial condition whose spectrum consists of a wide range of wavenumber, say $0\leq k\leq k_{max}$. Suppose that for each $k$, the dispersion-dissipation error satisfies $\vert\kappa^*-\kappa\vert\leq\epsilon_{\kappa}$, where $\epsilon_{\kappa}>0$ is a prescribed tolerance. Based on our analysis, for this initial condition, if an unoptimized scheme is used to approximate the spatial derivative, the grid size $\Delta x$ must satisfy $\kappa_{max}=k_{max}\Delta x\leq\gamma$ and $\vert c_{s+1}\kappa_{max}^{s+1}\vert\leq\epsilon_{\kappa}$ for some $\gamma$ and $\epsilon_{\kappa}$ in order to obtain a small dispersion-and-dissipation error. In the case of large $k_{max}$ and small $\gamma$, small $\Delta x$ is required. In contrast to the unoptimized FD schemes, the principle of an optimized DRB scheme is to seek a spatial discretization that for a large $\Delta x$ we still have $\vert\kappa^*-\kappa\vert\le \epsilon_{\kappa}$ uniformly for $0<\kappa<\Gamma$, where $\Gamma\gg\gamma$. Without going into the detail, Figure \ref{Fig:2_2_2}(a) \& (b) are $|\kappa^*-\kappa|$ versus $\kappa$ for the fourth-order optimized pentadiagonal scheme by Kim \cite{kim07} (referred to as ``KIM4TH'') and the optimized second-order pentadiagonal scheme (referred to as ``OPT2ND1P8'') described later in this paper, respectively. Compared with the unoptimized schemes in (c) \& (d),  the optimized schemes are able to control $|\kappa^*-\kappa|<\epsilon_{\kappa}=0.0015$  for $0 < \kappa < 1.8$. Especially, OPT2ND1P8 has $\epsilon_{\kappa} < 10^{-5}$ in this range of $\kappa$. It is worth noting that KIM4TH and OPT2ND1P8 are DRB schemes optimized from UNOPT10TH. All these three schemes in principle have the same computational cost. 

\begin{figure}[ht]
\begin{centering}
\includegraphics[scale=0.6]{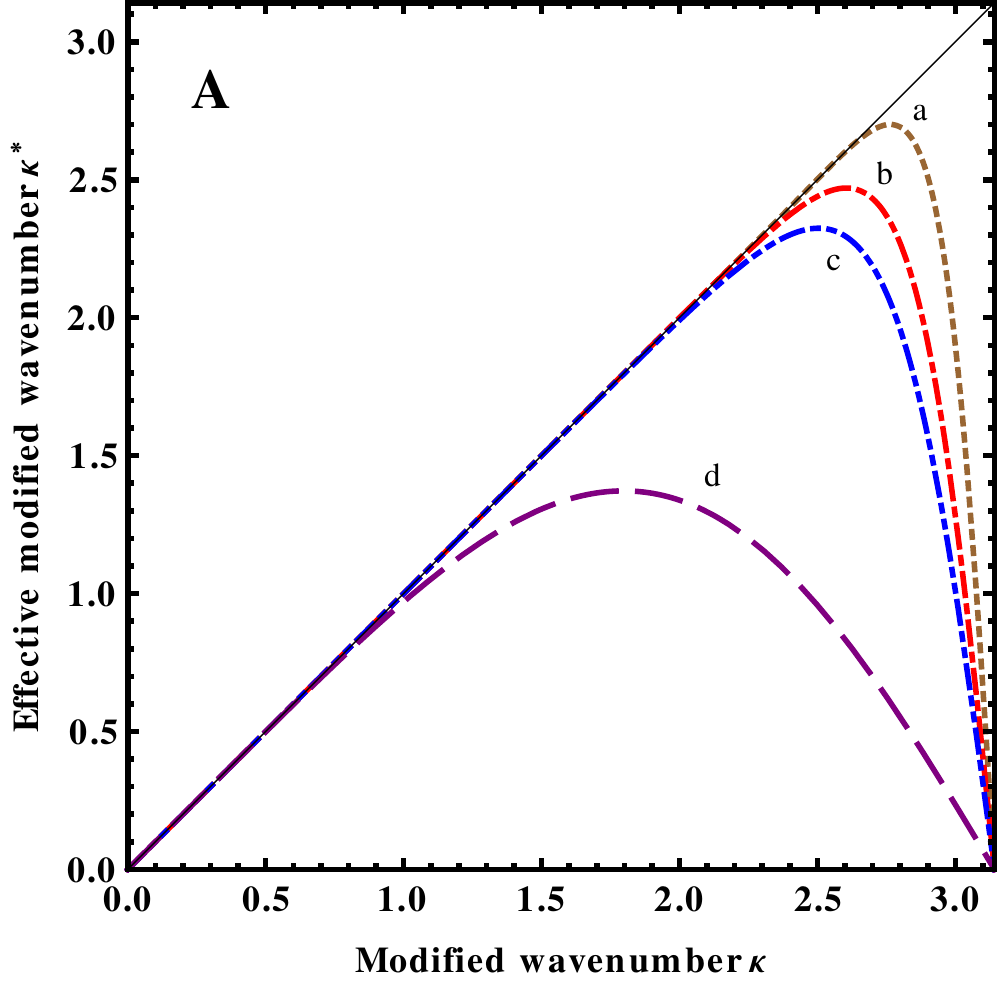}  \hskip 24pt
\includegraphics[scale=0.64]{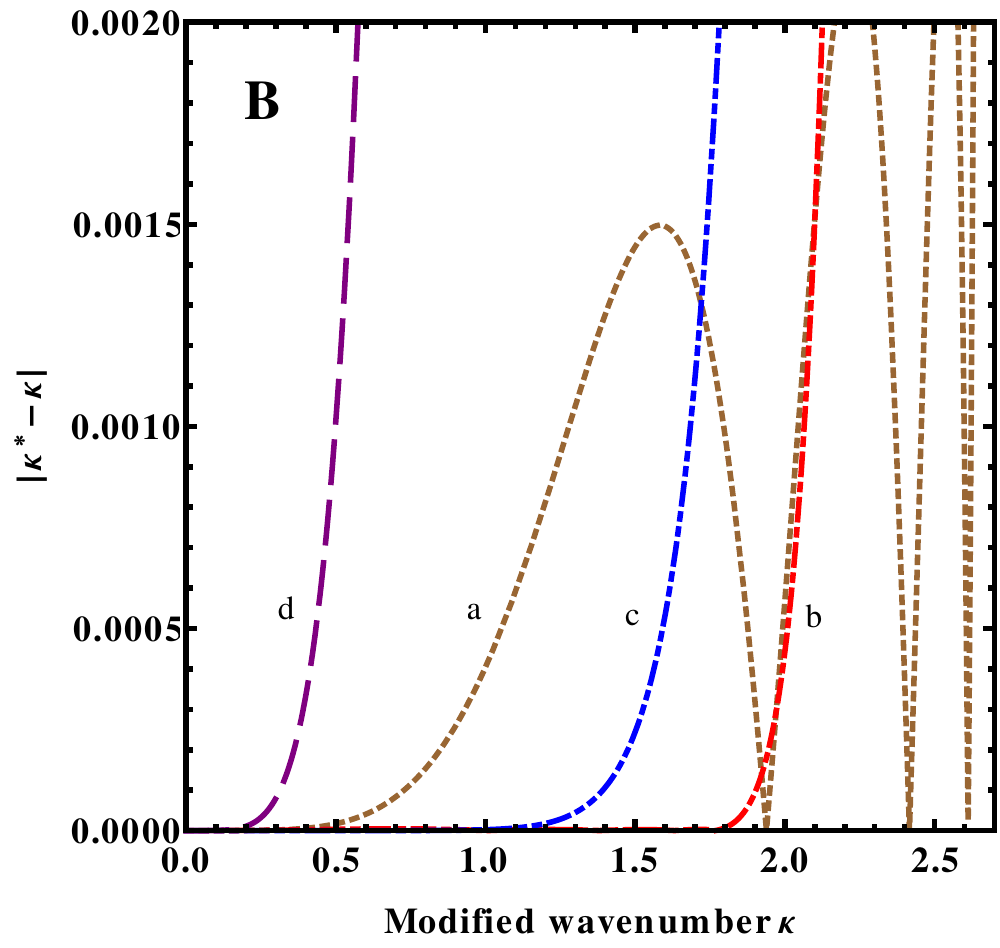}
\end{centering}
\caption{(A) : $\kappa^*$ versus $\kappa$ , and (B): $\vert\kappa^*-\kappa\vert$ versus $\kappa$ for (a) KIM4TH, (b) OPT2ND1P8, (c)UNOPT10TH, (d) CD4.}
\label{Fig:2_2}
\end{figure}

We plot the modified effective wavenumber $\kappa^*$ of these four discretizations in Figure \ref{Fig:2_2}(A), and their dispersion-and-dissipation errors in \ref{Fig:2_2}(B) (a collective plot of Figure \ref{Fig:2_2_2}). Figures \ref{Fig:2_2_2} and \ref{Fig:2_2} show that CD4 has the worst resolution characteristic among the four tested FD schemes, in the sense that within the interested range of wavenumber, CD4 has the worst dispersion-and-dissipation error. It is, however, not so clear from Figure \ref{Fig:2_2} that among the other three schemes which one has the most desirable resolution characteristic. At a glance of Figure \ref{Fig:2_2}(A), the effective wavenumber of KIM4TH seems to be the one closest to the real wavenumber. However, when we zoom in the dispersion-and-dissipation errors of these four schemes in Figure \ref{Fig:2_2}(B), we found that OPT2ND1P8 has the least dispersion-and-dissipation error between $0 < \kappa < 1.8$. Even UNOPT10TH has smaller error than that of 
KIM4TH before $\kappa\approx 1.6$. After $\kappa \approx 1.8$, the error of OPT2ND1P8 increases monotonically, while the error of KIM4TH begins to oscillate. The  question is whether smaller dispersion-and-dissipation error at high modified wavenumber produces smaller $L^2$-norm error for numerical solutions? We will further examine the $L^2$-norm error of the numerical solutions of the transport equation for these four schemes in the next section.


\subsection{Numerical examples}
In this section we examine the performance of several spatial discretization schemes for solving the transport equation with an initial condition traveling as a wave packet. The numerical experiments shed light on (1) the use of analysis derived previously in determining an accurate spatial discretization, and (2) the fundamental difference between the optimized and unoptimized schemes. 

\begin{example}[Initial condition with a localized spectrum]\label{ex:localized}

Consider the initial value problem \cite{bib:clawpack}

\begin{subequations}\label{2_50}
\begin{align}
u_t&+u_x=0,\;\quad 0\leq x\leq 1,\\
u(x,0)&=u_0(x)=\cos(30\pi x)\me^{-80(x-0.5)^2},
\end{align}
\end{subequations}
with periodic boundary conditions. We compute the numerical solutions at $T=1$ and evaluate the resulting $L^2$-norm errors according to the definition given by \eqref{2_43}. Our purpose is to study the relationship between the $L^2$-norm errors and the adopted spatial mesh sizes for both optimized and unoptimized FD schemes.


\begin{figure}[ht]
\begin{centering}
\includegraphics[scale=0.6]{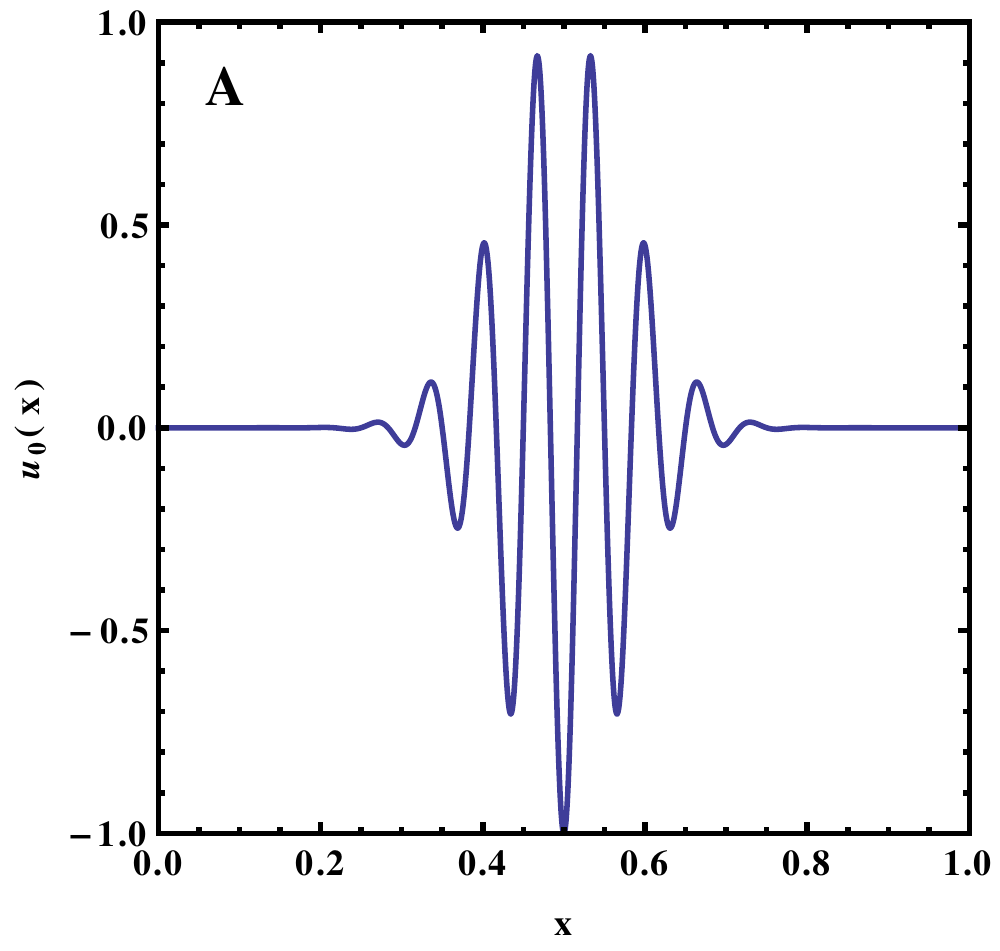} \hskip 36pt
\includegraphics[scale=0.6]{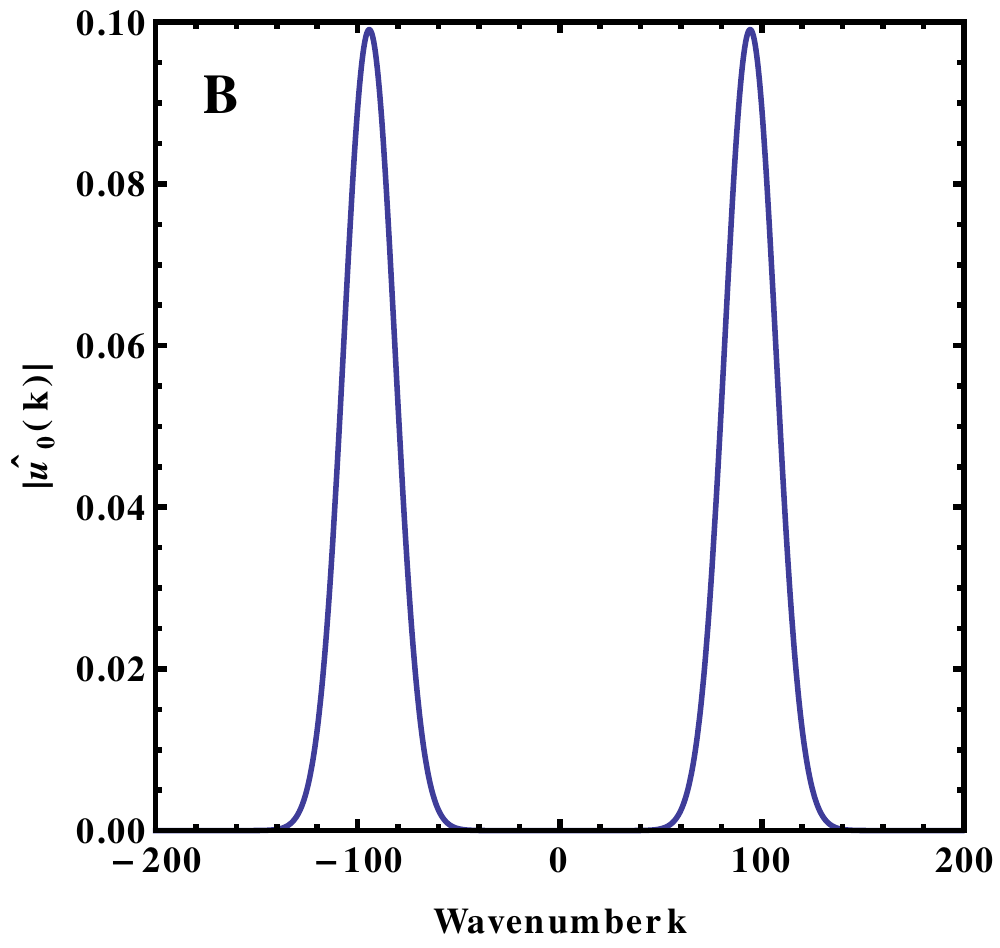}
\end{centering}
\caption{(A): The wave packet of $u_0(x)=\cos(30\pi x)\me^{-80(x-0.5)^2}$. (B): The spectrum of $u_0$. The wavenumber shown is $-200\leq k\leq 200$. The spectrum outside this range is neglected.}
\label{Fig:2_1}
\end{figure}


\begin{figure}[ht]
\begin{centering}
\includegraphics[scale=0.6]{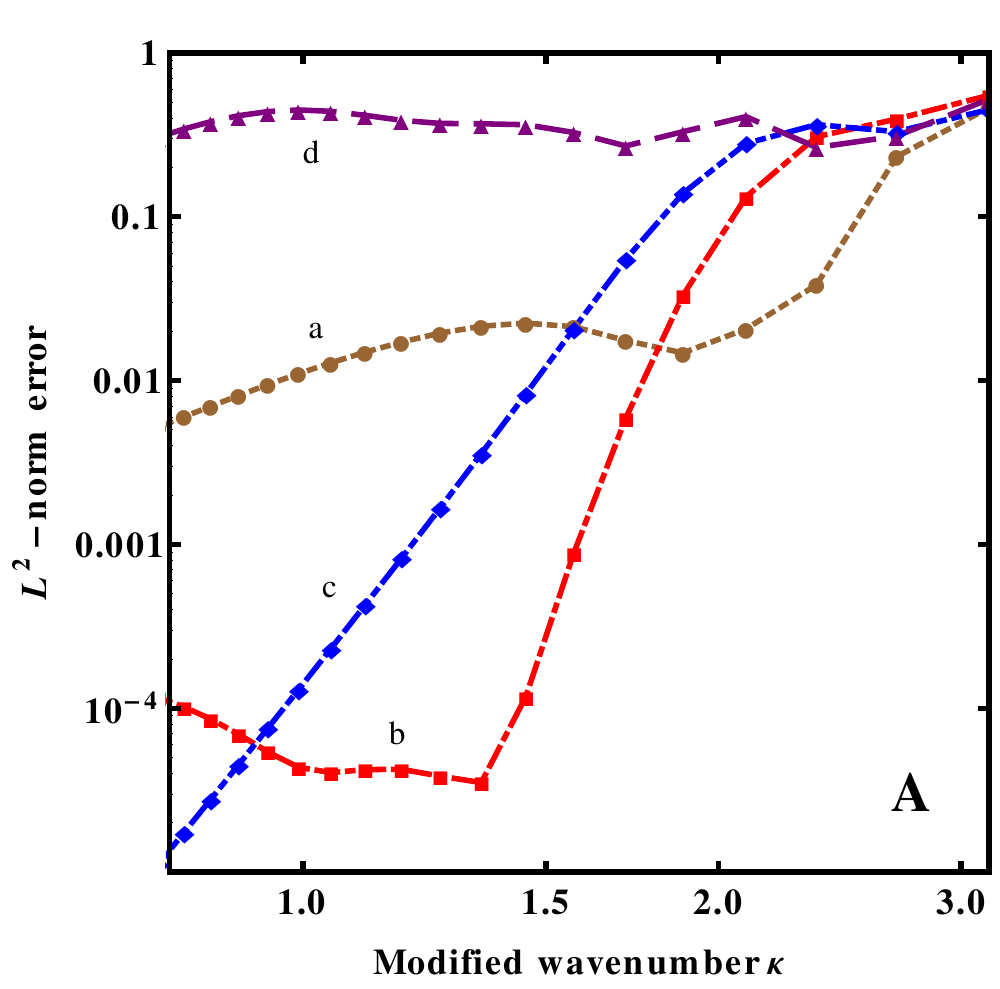} \hskip 0.6in
\includegraphics[scale=0.6]{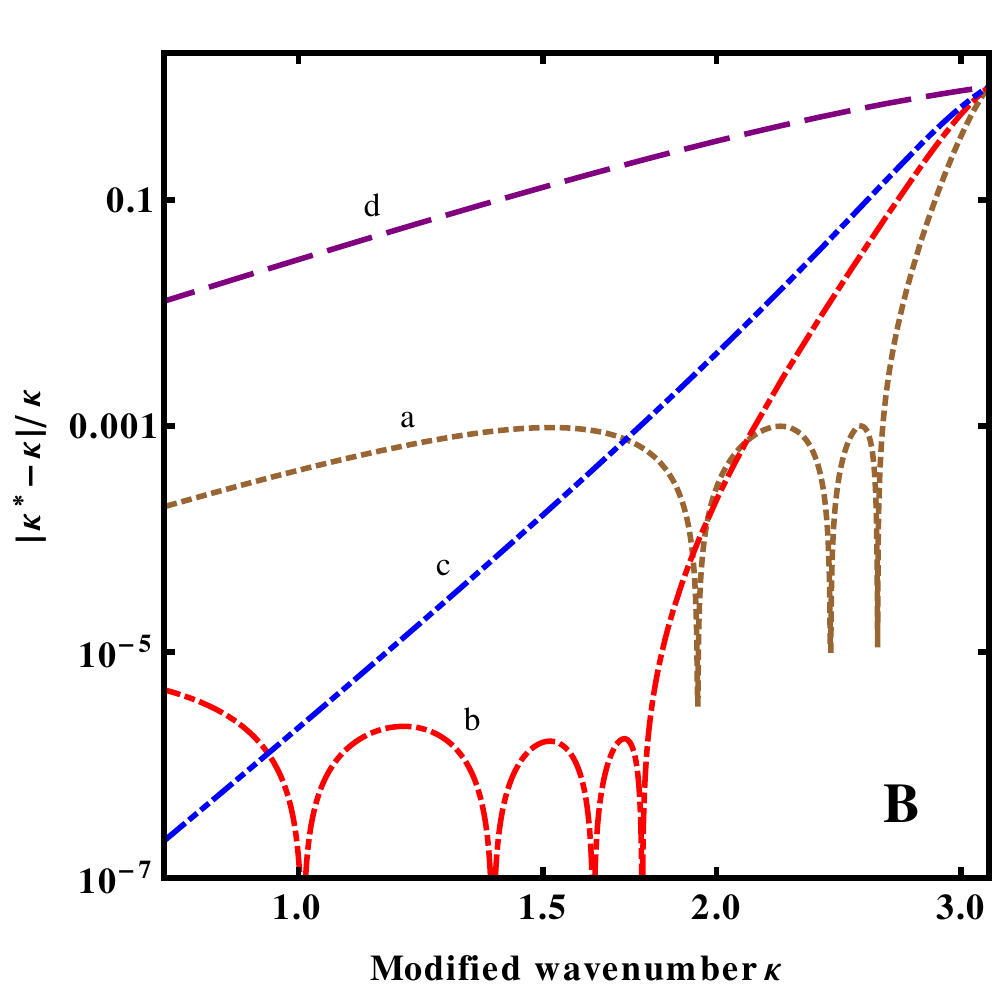}
\end{centering}
\caption{(A): The $L^2$-norm error versus $\kappa$ in log-log scale for \eqref{2_50} at $T=1$. Herein and after, the markers are where the errors are evaluated. (B): the $\vert\kappa^*-\kappa\vert/\kappa$ versus $\kappa$ in log-log scale for (a) KIM4TH, (b) OPT2ND1P8, (c) UNOPT10TH, (d) CD4.  In (A), the discrete modified wavenumber is computed with respect to $k=30\pi$, the primary wavenumber of the wave packet.}
\label{Fig:2_3}
\end{figure}


Due to the exponential factor, $u_0(x)$ can be considered a $1$-periodic function. The profile of $u_0(x)$ and its spectrum $\vert\widehat{u}_0(k)\vert$ are shown in Figure \ref{Fig:2_1}. It is clear that $\vert\widehat{u}_0(k)\vert$ is a binormal distribution centered at $k=\pm30\pi$ with standard deviation $\sigma_k=2\sqrt{80}\approx 17.9$. In practice, since the primary spectrum is a continuous smooth hump centered at $30\pi$ and the bandwidth is narrow,  the summation of $\vert\widehat{U}(k',0)\vert$  in  \eqref{2_43} can be approximated by a scalar multiple of $30\pi$. Hence we can treat $\vert\widehat{u}_0(k)\vert$ as a function localized at $k=\pm 30\pi$.
In our following numerical experiments, we assume this initial condition has single wavenumber $k=30\pi$.

The temporal scheme we adopt in our numerical experiments is the explicit eighth-order Runge-Kutta scheme (referred to as ``RK8''). The actual algorithm of RK8 is specified in \eqref{2_30}, \eqref{2_31}, and \eqref{2_32}, with $M=8$, $a_m=1/m!$ for $1\leq m\leq M=8$. We remark that the eigenvalues for a compact scheme like \eqref{2_10} are purely imaginary. The stability region of RK8 covers a large portion of imaginary axis, and the method is stable with a relatively large CFL number. However, in order to reduce the error generated by the temporal scheme, we choose a fixed CFL number $r=0.1$. Based on our previous analysis, with the setting of this temporal scheme, the dominant error should come from the spatial discretization, and the comparison among various spatial discretizations in principle is not influenced by the temporal integration. 

The spatial discretization schemes examined in our numerical experiments are UNOPT10TH, KIM4TH, OPT2ND1P8, and CD4. The first three are the implicit pentadiagonal compact schemes described in \eqref{2_10}. UNOPT10TH features the highest order-of-accuracy, in terms of Taylor truncation error, for this compact scheme, while CD4 is the only explicit scheme.

We now compute the solution of  \eqref{2_50} with a sequence of $\Delta x$ by using the aforementioned FD schemes to the final time $T=1$ . The $L^2$-norm error and the relative error $\vert\kappa^*-\kappa\vert/\kappa$ for each of the tested scheme are shown in Figure \ref{Fig:2_3}(A) \& (B), respectively. The modified wavenumber is chosen to be $0.8\leq\kappa\leq\pi$ for both (A) \& (B), with respect to the primary wavenumber, $k=30\pi$, of the wave packet.
Figure \ref{Fig:2_3}(A) shows that the $L^2$-norm error of CD4 does not decrease as $\kappa$ decreases, whereas Figure \ref{Fig:2_3}(B) shows that the relative error decreases at a rate of four-order-accuracy (slop of the straight line $\approx 4$). The discrepancy can be explained if we reexamine \eqref{2_45}. 

Recall the approximation of the $L^2$-norm error, in term of the phase error, in \eqref{2_45}.  This approximation is valid, provided the phase error (proportional to the relative error) in \eqref{2_44} is small, so that the exponential can be approximated by its leading-order term. For this specific problem, the relative error for CD4 is larger than $2\times 10^{-3}\approx 0.2\%$ and $|\theta_e|$ is larger than $0.19$. This is evidently too big to make the approximation in \eqref{2_45} valid. This explains why the $L^2$-norm error of CD4 does not decrease as $\kappa$ decreases. Furthermore, phase error depends on the coefficient $ckT$. We expect that we will not observe the linear relation between the $L^2$-norm error and the relative error for CD4, unless much smaller $\Delta x$ or $T$ are applied.  Similar behavior is observed for other unoptimized schemes as well. For example, the $L^2$-norm error of the unoptimized tenth-order scheme, UNOPT10TH, does not decrease at the same rate as the relative error until $\kappa\approx 2$. 

\begin{note}
For a fixed spatial discretization and a fixed $\kappa$, if a smaller mesh size becomes available due to increased computer memory, the scheme has the potential to support higher wavenumbers. However, the phase error depends linearly on the wavenumber (with respect to a fixed $\kappa$), thus the overall error may increase. If this is this case, it is necessary to derive a new optimized scheme that accounts for the new error tolerance.   
\end{note}

For optimized DRB schemes, the situation is similar but more complicated. The fourth-order optimize  scheme KIM4TH sacrifices the accuracy at the lower modified wavenumber to control the relative error to the level below $10^{-3}$ for $0 <\kappa<2.7$. The optimized scheme OPT2ND1P8 developed in this paper, on the other hand, controls the relative error to the level below $10^{-5}$ for $0< \kappa <1.8$. As a result, the $L^2$-norm error of these scheme are small as well.  We see that for  $\kappa\le 1.8$, OPT2ND1P8 has a smaller $L^2$-norm error than that of KIM4TH, and for $\kappa >1.8$, KIM4TH produces the most accurate solution among the four tested schemes. 
\end{example}

\begin{note}
We make two observations about this numerical experiment.
\begin{enumerate}
\item[(a)] The typical practice of plotting the effective modified wavenumber versus modified wavenumber, as in Figure \ref{Fig:2_2}(A), may be misleading. Merely visually matching $\kappa^{*}$ to $\kappa$ at high modified wavenumber at this scale does not guarantee the most appropriate (or accurate) numerical scheme for solving the given problem. 
Indeed, the process of identifying such a scheme will depend on the problem itself and the available computer power.  For example, in experiment \ref{ex:localized}, the primary wavenumber is $k=30\pi$. If the computer memory allows us to refine the grid to $0.7<\kappa=k\Delta x<1.8$, then OPT2ND1P8 may be the most appropriate scheme to obtain the most accurate numerical solution. For the most extreme case, if $\kappa$ can be controlled to the range below 0.7, then the tenth-order unoptimized scheme may be the best. However, small $\Delta x$ usually implies small $\Delta t$ so that large CPU/GPU time is required. We have to keep that in mind as well to design (or select) a proper numerical scheme for our problems. 
\item[(b)]  To design a DRB scheme, ideally, we need to take into account the limitation of our computer power, including the size of memory, CPU/GPU power, and the ability of parallelization, as well as the spectrum of the solution of the  problem. Unfortunately, such information is rarely available.  Kim \cite{kim07}, however, demonstrates a good basic philosophy for designing a DRB scheme in the absence of this information. That is, maintaining the relative error $\vert\kappa^*-\kappa\vert/\kappa\leq10^{-3}$ for $0\leq\kappa\leq\Gamma_r$, while keeping $\Gamma_r$ as large as 2.6. As  a result, the $L^2$-norm error decreases almost right away, after $\kappa$ starts to decrease from the modified wavenumber around 2.7.  In the next section, we will outline the basic principles for designing an optimized DRB scheme. Moreover, we demonstrate these principles by designing new schemes.
\end{enumerate}
\end{note}

\begin{example}[Initial condition with a broad spectrum]\label{ex:broad}
Consider the following initial-value problem 
\begin{subequations}\label{2_52}
\begin{align}
&u_t+u_x=0,\\
&u(x,0)=u_0(x)=\cos(k(x)x),
\end{align}
\end{subequations}
where $k(x)=30\pi\exp(-0.75(x-5)^2)$ and $0\le x\le 10$.  Since $k(x)$ is an exponential function, $u_0(0^{+})=u_0(10^{-})\approx 1$. We place this initial condition in a periodic domain, with the period equal to 10. 
The spectrum $\widehat{u}_0(k)$ has wavenumber in the range of $0\leq\vert k\vert\leq380$.  The distribution of the spectrum is shown in Figure \ref{Fig:2_4}(A). The exponential function emulates the short length of data record,  which results in the so-called spectral leakage (side-lobe). The instantaneous wavenumber (analog to the instantaneous frequency in time series analysis) is defined as $k_c(x)=\dif\theta/\dif x$, where $\theta=k(x)\,x$. Hence locally $\max\{k_c(x)\} = \max\{k'(x)x+k(x)\} \approx 360$, for $x\in [0, 10]$. However, if we look at the spectrum of $u_{0}$, we find that the maximal wavenumber at which the spectrum is significantly nonzero is around 380. 
For this example, we define $k_{\mathrm{max}}=380$.

\begin{figure}[ht]
\begin{centering}
\includegraphics[scale=0.6]{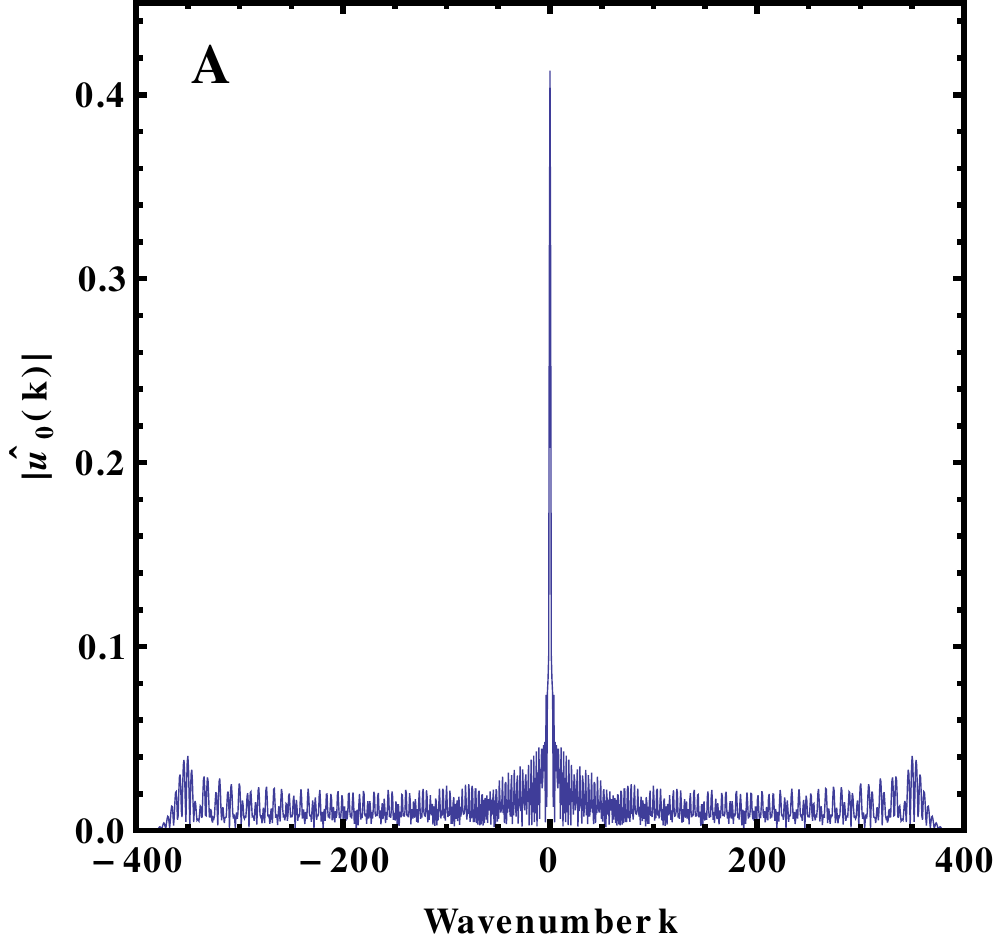} \hskip 0.5in
\includegraphics[scale=0.6]{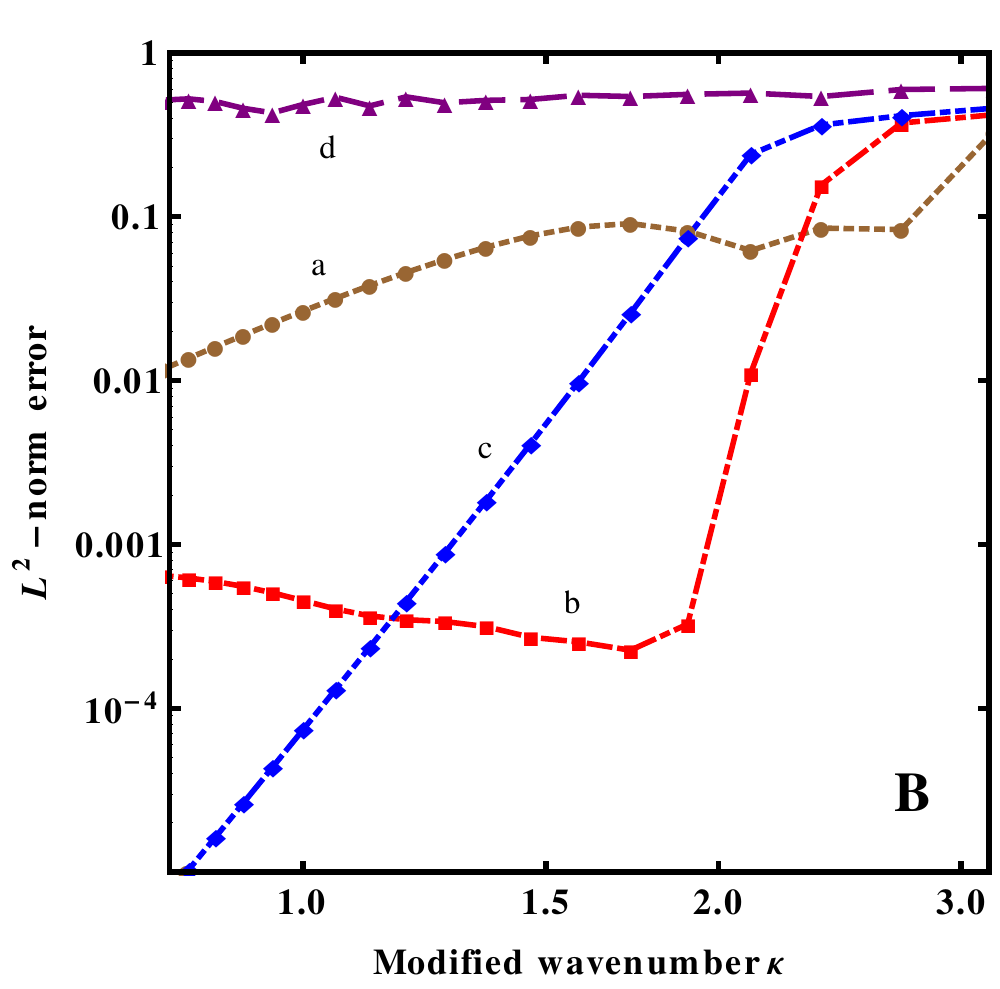}
\end{centering}
\caption{(A): The spectrum $\vert\widehat{u}_0(k)\vert$ for the initial condition $u_0(x)$ in Equation \eqref{2_52}. (B): The $L^2$-norm error versus $\kappa$ at the final time $T=1$ in the log-log scale for (a) KIM4TH, (b) OPT2ND1P8, (c) UNOPT10TH, (d) CD4.  In (B), the discrete modified wavenumber is computed with respect to $k_{\mathrm{max}}=380$.}
\label{Fig:2_4}
\end{figure}

As we did for Example \ref{ex:localized}, we solve the initial-value problem (\ref{2_52}) with various $\Delta x$, and compute the $L^2$-norm error. Figure \ref{Fig:2_4}(B) shows the errors plotted against the modified wavenumber, $\kappa=k_{\mathrm{max}}\Delta x=380\Delta x$, for (a) KIM4TH, (b) OPT2ND1P8, (c) UNOPT10TH, and (d) CD4, respectively. 
Structure-wise, Figure \ref{Fig:2_3}(A) and Figure \ref{Fig:2_4}(B) are virtually identical, in particular, for the distinguishable feature between the optimized and the unoptimized (Taylor series expansion based) schemes.  
We now use the schematic plots (``cartoons'') in Figure \ref{Fig:2_5} to illustrate the fundamental difference between unoptimized and optimized schemes and the principles of designing a ``good'' DRB scheme.


\begin{figure}[tbh]
\begin{centering}
\includegraphics[scale=0.6]{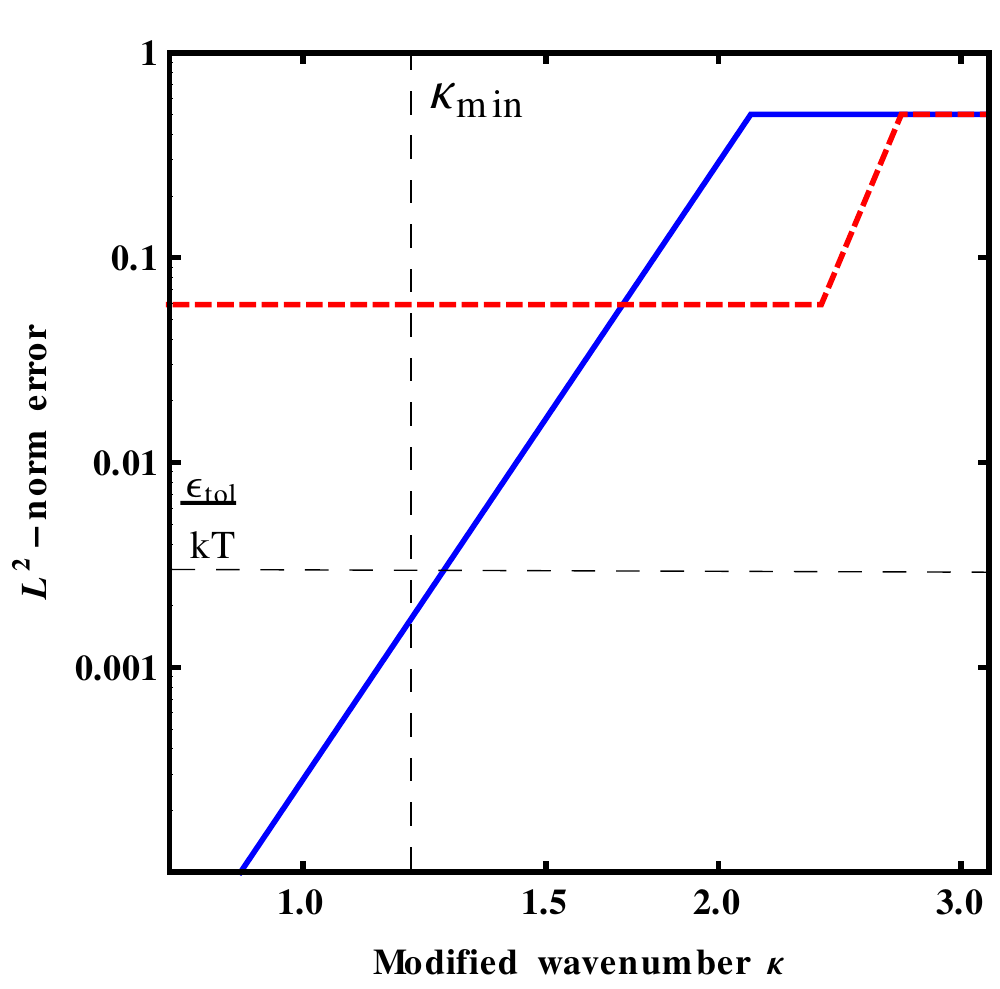} \hskip 0.5in
\includegraphics[scale=0.6]{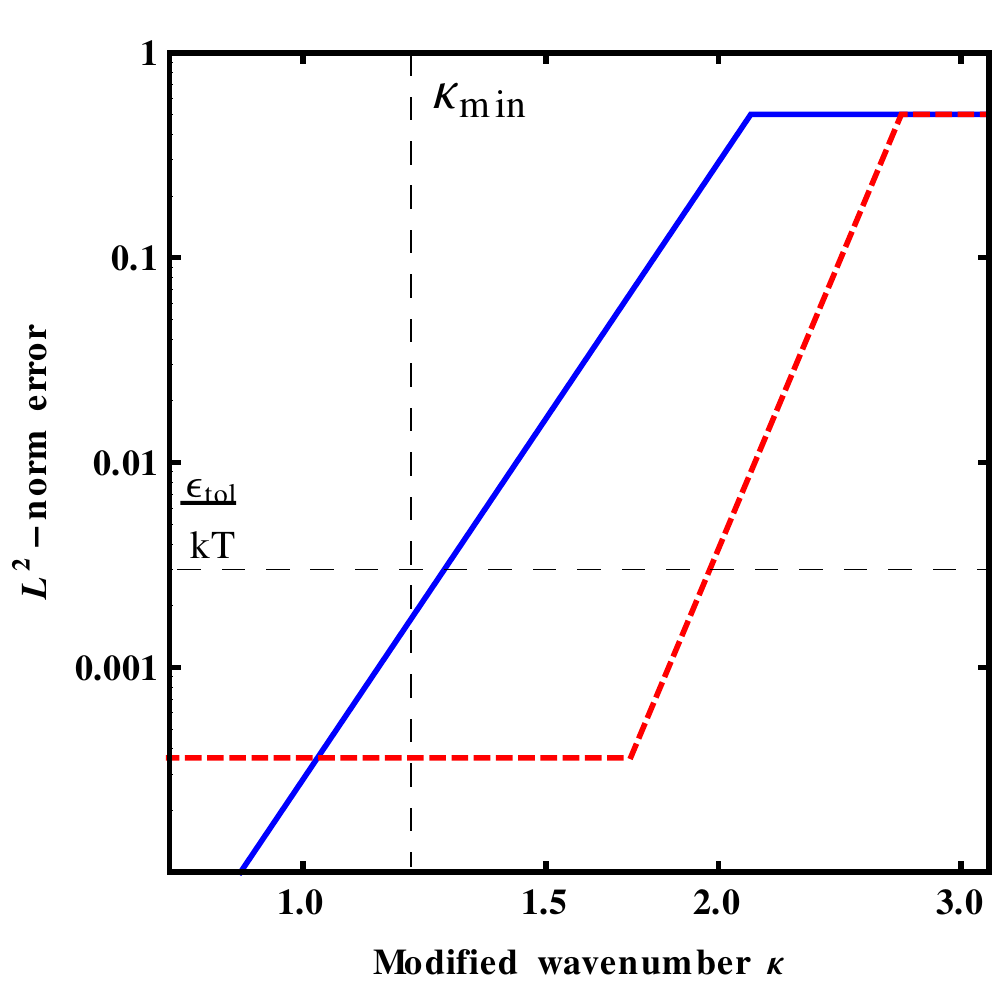}
\end{centering}
\caption{Schematic illustrations of the errors versus modified wavenumber for unoptimized and optimized DRB schemes. The horizontal long-dashed line is the effective maximum error tolerance. The vertical long-dashed line is the minimum modified wavenumber that a computer can resolve. The solid line is the error of an unoptimized scheme, and the short-dashed is the error of an optimized scheme.}
\label{Fig:2_5}
\end{figure}

We define the desired maximum error tolerance and the smallest mesh size allowed by the computer power for a simulation to be $\epsilon_{\mathrm{tol}}$ and $\Delta x_{\mathrm{min}}$, respectively, where $\epsilon_{\mathrm{tol}}$ is specified in some metric norm, such as $L^2$ norm. From \eqref{2_45}, it makes sense to define the \emph{effective error tolerance} as $\epsilon_{\mathrm{tol}}/kT$, where $k$ is a specified wavenumber and $T$ the final time for the simulation. We also define the minimum modified wavenumber to be $\kappa_{\mathrm{min}}=k\Delta x_{\mathrm{min}}$.

Now, in Figure \ref{Fig:2_5}, $\epsilon_{\mathrm{tol}}/kT$ is represented by the horizontal long-dashed line, while $\kappa_{\mathrm{min}}$ is the vertical long-dashed line. The solid line is the $L^2$-norm error of some unoptimized scheme, and the short-dashed line represents that of an optimized DRB scheme. For an unoptimized scheme, it is well-known that for a fixed wavenumber, at the region of large modified wavenumber (coarser grids), the error does not decrease according to the order-of-accuracy of the scheme, until the effective modified wavenumber is reasonably small (finer grids). The scheme will continue to enjoy the decrease of error when refining the grid until the modified wavenumber hits $\kappa_{\mathrm{min}}$ (where the memory of the computer cannot afford to refine the grid anymore).  If this point is below the horizontal line of $\epsilon_{\mathrm{tol}}/kT$ (as shown in Figure \ref{Fig:2_5}(A)), then the unoptimized scheme is appropriate for the problem, otherwise we need to design a better numerical scheme for the problem. 

For a typical optimized DRB scheme, the error decays earlier and faster when refining the grid than its unoptimized counterpart. After that the error may decay much slower or even stop decaying before it starts decaying again at a relatively small modified wavenumber when the effect of the order-of-accuracy takes over. Suppose when an optimized DRB scheme enters the aforementioned ``plateau'' region, the error is bigger than $\epsilon_{\mathrm{tol}}/kT$ and the plateau extends to the region over $\kappa_{\mathrm{min}}$, as shown in Figure \ref{Fig:2_5}(A),  then this optimized DRB scheme is not suitable for the problem. On the other hand, if the error decreases to be lower than  $\epsilon_{\mathrm{tol}}/kT$ before entering the plateau region, while $\kappa > \kappa_{\mathrm{min}}$, as shown in Figure \ref{Fig:2_5}(B), then this optimized DRB scheme could be advantageous for solving the problem.
Consider Figure \ref{Fig:2_4}(B) of Example \ref{ex:broad} as an illustration of the above point of view. Suppose that our desired effective maximum error tolerate is $\epsilon_{\mathrm{tol}}/kT=0.001$ and our computer can only resolve $\kappa_{\mathrm{min}}=1.5$, then OPT2ND1P8 is the only scheme capable of delivering the desired result among the four methods.

\begin{note}
It is known that FD approximation introduces both truncation and dispersion-and-dissipation errors. Traditionally, the development of FD schemes has focused on reducing the Taylor truncation error and associated stability issues. After Lele \cite{lele92} and Tam et al. \cite{twd93} introduced the concept of reducing/minimizing the dispersion-and-dissipation error, a large variety of optimized FD schemes have been proposed. However, the development of these schemes, like the design of KIM4TH described in \cite{kim07}, is (i) built on spatial discretizations that have their origins in Taylor expansion   
 and (ii) aim to balance 
order of accuracy and dispersion-and-dissipation error. The schemes OPT2ND1P5 and OPT2ND 1P8 tilt the balance entirely toward dispersion-and-dissipation error. Thus, in contrast to the traditional perspective of Taylor truncation error, we devote virtually all of our effort to minimizing the dispersion-and-dissipation error.  The constraint of the optimization is merely first-order accuracy. The second-order accuracy of the scheme is achieved automatically when symmetry is imposed. It is worth noting that from the error dynamics described in Figure \ref{Fig:2_5}, if the solution of our problem concentrates at high-wavenumber region, perhaps it is futile to maintain high order accuracy when designing a DRB scheme. It makes sense to focus on reducing the dispersion-and-dissipation error.

\end{note}

\end{example}

\section{Optimized DRB schemes}\label{sec:DRB-design}

In this section, we discuss the principles of optimization process. We then extend these principles to propose a new optimized DRB scheme whose derivation is conceptually different from the prevailing techniques utilized to design DRB schemes. Finally, we show that our new scheme gives the least $L^2$-norm error among the tested algorithms for a large range of $\kappa$.

\subsection{Designing principles of optimized DRB schemes}\label{sec:general}
The scheme OPT2ND1P8 may serve as an exemplar of the typical approach to deriving an optimized DRB scheme. Conventionally, the process begins with a specific form of discretization, such as \eqref{2_10}. The next step is to take the Fourier Transform to determined the effective (numerical) modified wavenumber $\kappa^*$ in terms of the unknown coefficients of the discretization and the real wavenumber $\kappa$, such as \eqref{2_20}. The optimization (minimization) is against the objective (error) function 
\begin{equation}\label{3_2}
E\equiv\int_0^{\Gamma}\vert\kappa^*-\kappa\vert^2W(\kappa)\,\dif\kappa\,,
\end{equation}
where $W(\kappa)$ is a weight function chosen to facilitate the minimization process. 
This process is usually coupled to equations that arise from the elimination of some lower-order Taylor truncation errors; these appear as constraints in the minimization process. Thus, such FD schemes are designed to have a small error between the numerical and real wavenumber for a prescribed error tolerance and to maintain certain order of accuracy.

We remark that a weight function $W(\kappa)\ne 1$ may ease the optimization process \cite{kl96}, but there is no other benefit. Simply put, the minimizer obtained from minimizing the objective function with $W(\kappa)\ne 1$, need not be a minimizer for the objective function
\begin{equation}\label{3_1}
E\equiv\int_0^{\Gamma}\vert\kappa^*-\kappa\vert^2\,\dif \kappa\,.
\end{equation}
In other words, with a weight function $W(\kappa)\ne 1$, instead of minimizing the dispersion-and-dissipation error, we minimize the product of the dispersion-and-dissipation error and a (continuous) function.
Often, the minimizer found by using \eqref{3_2} does not give the minimum dispersion-and-dissipation error.

The scheme OPT2ND1P8 is derived from the discretization \eqref{2_10} and the choice of $\kappa^*$ given in \eqref{2_20}. We retain only the first equation in \eqref{eq:padeconstraints}, and the unknowns $a$, $b$, $c$ together with $\alpha$ and $\beta$ are derived from minimizing the objective function \eqref{3_1}. 
Hence, OPT2ND1P8 is a second-order accurate method in terms of Taylor truncation error, while most of the degrees of freedom are dedicated to minimizing the dispersion-and-dissipation error. 
As our experiments show, the number of degrees of freedom reserved for optimizing the dispersion-and-dissipation error affects the behavior of $L^2$-norm error. That is, using more degrees of freedom to minimize the dispersion-and-dissipation error tends to result in earlier decay at higher modified wavenumber and smaller $L^2$-norm error in the plateau region. 
The trade off is that the error behavior after the plateau region is governed by the lower-order accuracy of the scheme. However, that behavior is usually not so important, because most likely it will happen after $\kappa_{\mathrm{min}}$ shown in Figure \ref{Fig:2_5}. Indeed, for such small $\kappa$, we expect that a high-order unoptimized scheme will offer better performance. 

The value of $\Gamma$, the upper limit of integration,  in the objective function \eqref{3_1} also plays an important role in the design of suitable optimized DRB schemes.
Our numerical experiments suggest that using larger values of $\Gamma$ causes the $L^2$-norm error to start decaying at large modified wavenumber. But this comes at a cost; the error at lower modified wavenumbers (especially the plateau region) is large. In addition, in the case of periodic boundary conditions, $\kappa^*$ is usually a  $2\pi$-periodic function in $\kappa$.  In this scenario, $\vert\kappa^*-\kappa\vert$ is bounded, but is large near $\kappa=2\pi$, and the range of $\vert\kappa^*-\kappa\vert$ may extend several orders of magnitude. Larger values of $\Gamma$ cause difficulty in the numerical optimization process and require high precision to achieve the desired stopping criteria. These issues do not arise if smaller values of $\Gamma$ are used. In our numerical experiments, we choose $\Gamma=1.8$ for OPT2ND1P8, and in the next section we will compare OPT2ND1P8 with OPT2ND1P5 whose $\Gamma$ value is 1.5. 

For the DRB schemes in this paper (and for many in the literature), the objective function is as in \eqref{3_1} with the kernel $|\kappa^*-\kappa|^2$. We view this kernel as the product $|\kappa^*-\kappa|\,|\kappa^*-\kappa|$, where the second factor of $|\kappa^*-\kappa|$ behaves as a ``dynamical'' weight function for the first one.  Consequently, more weight is put on the modified wavenumber where the dispersion-and-dissipation error is large.  We observe that, in principle, there is nothing preventing us from using a different kernel $|\kappa^*-\kappa|^n$, where $n\ne 2$. However, our numerical experiments show that the dynamics of the dispersion-and-dissipation error are very similar for $n=1, 2,$ and $4$. However, in the case $n=1$, the lack of smoothness in the kernel hinders the optimization process. Consequently, the optimization process requires a large working precision and is more time consuming. 

Finally,  the optimization process depends on the form of $\kappa^*$. The function $\kappa^*$ is fixed (and hence the number of degrees of freedom), once the discretization is determined. 
To sidestep this, we now propose a new approach to determining $\kappa^*$. Instead of first determining the discretization and then finding the function $\kappa^*$, we begin by prescribing the function $\kappa^*$ in terms of unknown coefficients first. Then we determine the corresponding discretization that matches the prescribed $\kappa^*$.  In particular, we propose an unknown form of spatial discretization for which the numerical modified wavenumber takes the form
\begin{subequations}\label{3_6}
\begin{align}
\kappa^*&=\dfrac{a\sin(\kappa)+(b/2)\sin(2\kappa)+(c/3)\sin(3\kappa)}{1+2\alpha\cos(\kappa)+2\beta\cos(2\kappa)} + \left(d\sin(\kappa)+\dfrac{e}{2}\sin(2\kappa)+\dfrac{f}{3}\sin(3\kappa)\right) \\
&= \kappa_1^* + \left(d\sin(\kappa)+\dfrac{e}{2}\sin(2\kappa)+\dfrac{f}{3}\sin(3\kappa)\right)\,.
\end{align}
\end{subequations}
The import of the choice in \eqref{3_6} is that we gain more degrees of freedom to optimize the dispersion-and-dissipation error. It remains to find the corresponding discretization. We recall that the first term, $\kappa_1^*$, arises from the discretization, \eqref{2_10}---or equivalently \eqref{2_16}, via the linear system $\textbf{L}\textbf{U}'=\textbf{R}\textbf{U}$.  The numerator of  $\kappa_1^*$, namely the term  
\[
\sin(\kappa)+(b/2)\sin(2\kappa)+(c/3)\sin(3\kappa)\,,
\] 
results from the matrix-vector multiplication $\textbf{R}\textbf{U}$, while the denominator,  $1+2\alpha\cos(\kappa)+2\beta\cos(2\kappa)$, comes from the multiplication of $\textbf{L}^{-1}$ and $\textbf{R}\textbf{U}$ (i.e., solving the linear system for $\textbf{U}'$). 
Since the second part of $\kappa^*$ has the same form as the numerator of $\kappa_1^*$, the numerical derivatives at the grid points are found by 
\beq\label{3_6_1}
\textbf{U}'=\textbf{U}_1'+\textbf{W},
\eeq
where $\textbf{U}_1' = \textbf{L}^{-1}\textbf{R}\textbf{U}$ and $\textbf{W}=\tilde{\textbf{R}}\textbf{U}$, so that the corresponding compact discretization scheme for $\kappa^*$ takes the form of 
\beq\label{3_7}
\textbf{L}\textbf{U}' = \textbf{R}\textbf{U} + \textbf{L}\tilde{\textbf{R}}\textbf{U},
\eeq
for which $U_j'$ is given by
\beq\label{3_7_1}
U_j' = \mi k^*U_j,
\eeq
provided $U_j=\me^{\mi k x_j}$. Here, $\textbf{L}$ and $\textbf{R}$ are the cyclic pentadiagonal and heptadiagonal matrices defined by \eqref{2_10}, and $\tilde{\textbf{R}}$ is a cyclic heptadiagonal matrix that shares exactly the same structure as $\textbf{R}$,  except the non-zero entries are $\pm d/(2\Delta x)$, $\pm e/(4\Delta x)$, and $\pm f/(6\Delta x)$. Therefore the discretization can be expressed explicitly as
\begin{multline}\label{5_14}
\beta U'_{j-2}+\alpha U'_{j-1}+U'_{j}+\alpha U'_{j+1}+\beta U'_{j+2}
\\=  a\dfrac{U_{j+1}-U_{j-1}}{2\Delta x}+b\dfrac{U_{j+2}-U_{j-2}}{4\Delta x}+c\dfrac{U_{j+3}-U_{j-3}}{6\Delta x}
\\+ \beta w_{j-2}+\alpha w_{j-1}+w_{j}+\alpha w_{j+1}+\beta w_{j+2}\,,
\end{multline}
where
\begin{equation}\label{5_15}
w_j=d\dfrac{U_{j+1}-U_{j-1}}{2\Delta x}+e\dfrac{U_{j+2}-U_{j-2}}{4\Delta x}+f\dfrac{U_{j+3}-U_{j-3}}{6\Delta x}.
\end{equation}

We remark that the computational costs for solving the linear systems \eqref{2_16} and \eqref{3_7} are virtually the same; The principal costs is in finding $\textbf{L}^{-1}$. We refer to the compact DRB scheme obtained by using the discretization \eqref{5_14} and the numerical dispersion relation \eqref{3_6} as KLL2ND.

\begin{note}
Equation \eqref{3_6_1} is the corresponding discretization for the prescribed numerical modified wavenumber (\ref{3_6}). In fact, if we look at  \eqref{3_6}, \eqref{3_6_1}, and  \eqref{3_7} together, we see that mathematically $\kappa^*_1$  is associated with the linear transformation ${\bf{L}^{-1}{\bf{R}}}$ and 
\[
\left(d\sin(\kappa)+\dfrac{e}{2}\sin(2\kappa)+\dfrac{f}{3}\sin(3\kappa)\right)
\]
is associated with another linear transformation $\tilde{{\bf{R}}}$. Therefore, in principle, we can prescribe any numerical modified wavenumber $\kappa^{*}$ in the form $\kappa^{*}=\kappa_1^*+\cdots+\kappa_M^*$. 
As long as we know the linear transformation that corresponds to each $\kappa_i^*$, we may recover the corresponding discretization for $\kappa^{*}$.
\end{note}

\begin{note}
The introduction of the spatial discretization (\ref{5_14}) does not make much sense from the view point of Taylor truncation error. The philosophy behind the prescribed numerical modified wavenumber $\kappa^*$ in \eqref{3_6}, which leads, in turn, to the discretization \eqref{5_14}, is described as follows. 
%
%
By proposing \eqref{3_6}, we introduce a conceptually different perspective.  We introduce the resolution characteristic before finding the corresponding discretization. While similar to OPT2ND1P5 and OPT2ND1P8, we devote virtually all of our effort to minimizing the dispersion-and-dissipation error, the extra degrees of freedom gained from the additional unknowns in the prescribed resolution characteristic further reduce the dispersion-dissipation error in the minimization process, compared with OPT2ND1P5 and OPT2ND1P8.  The advantages of this approach are shown in the next round of numerical experiments.
\end{note}

\subsection{More numerical examples}
\begin{figure}[ht]
\begin{centering}
\includegraphics[scale=0.55]{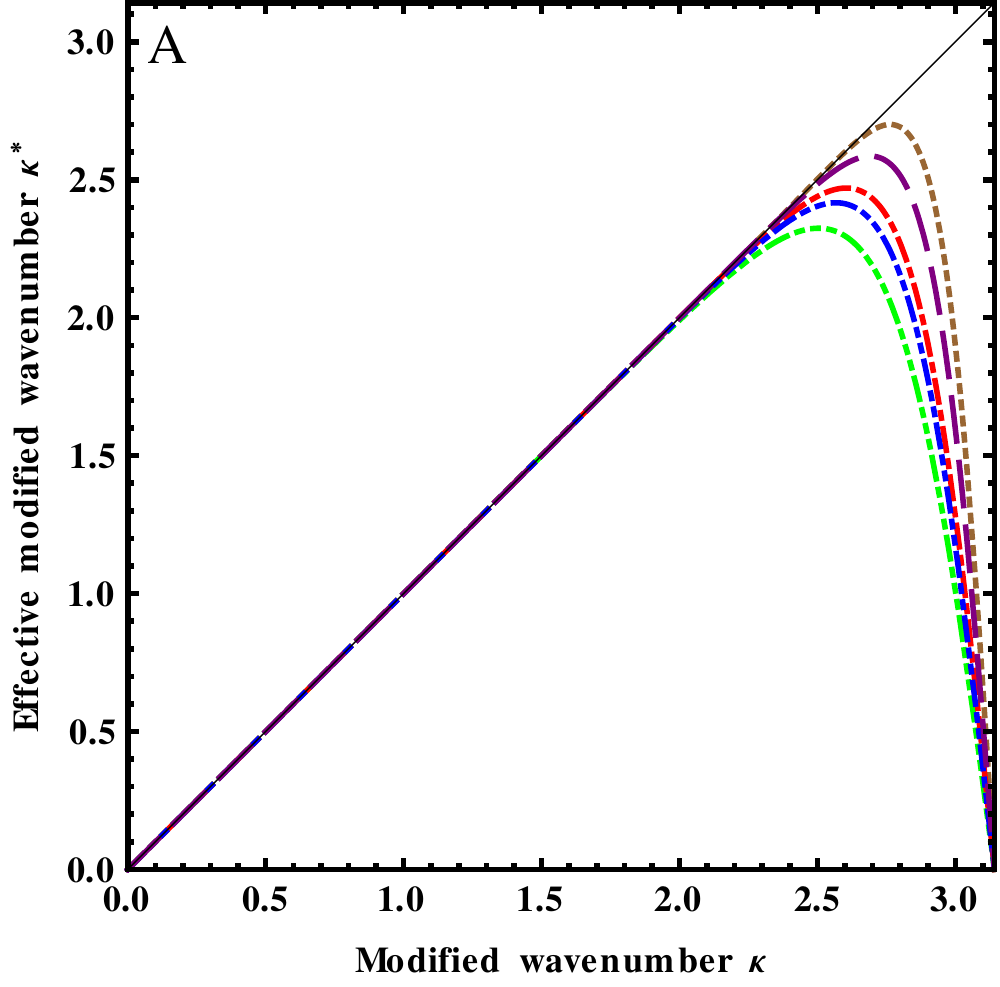}\hskip 0.5in
\includegraphics[scale=0.6]{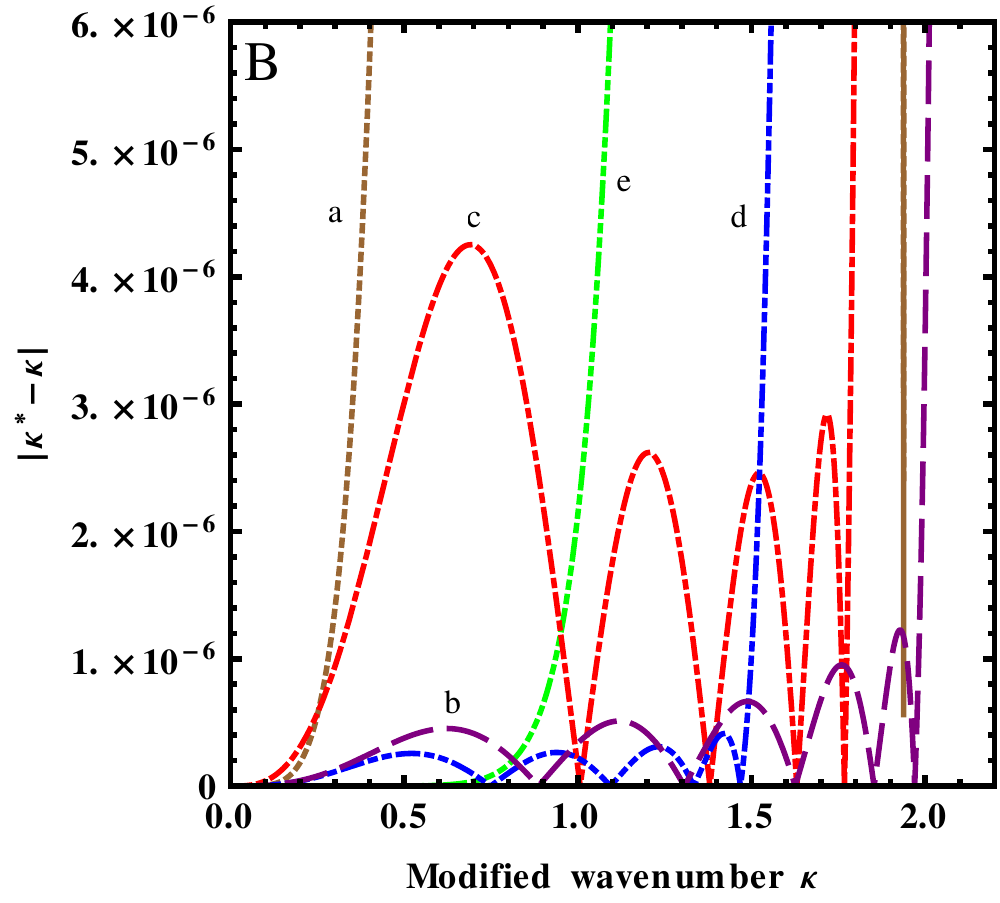}\\ \vskip 0.3in
\includegraphics[scale=0.6]{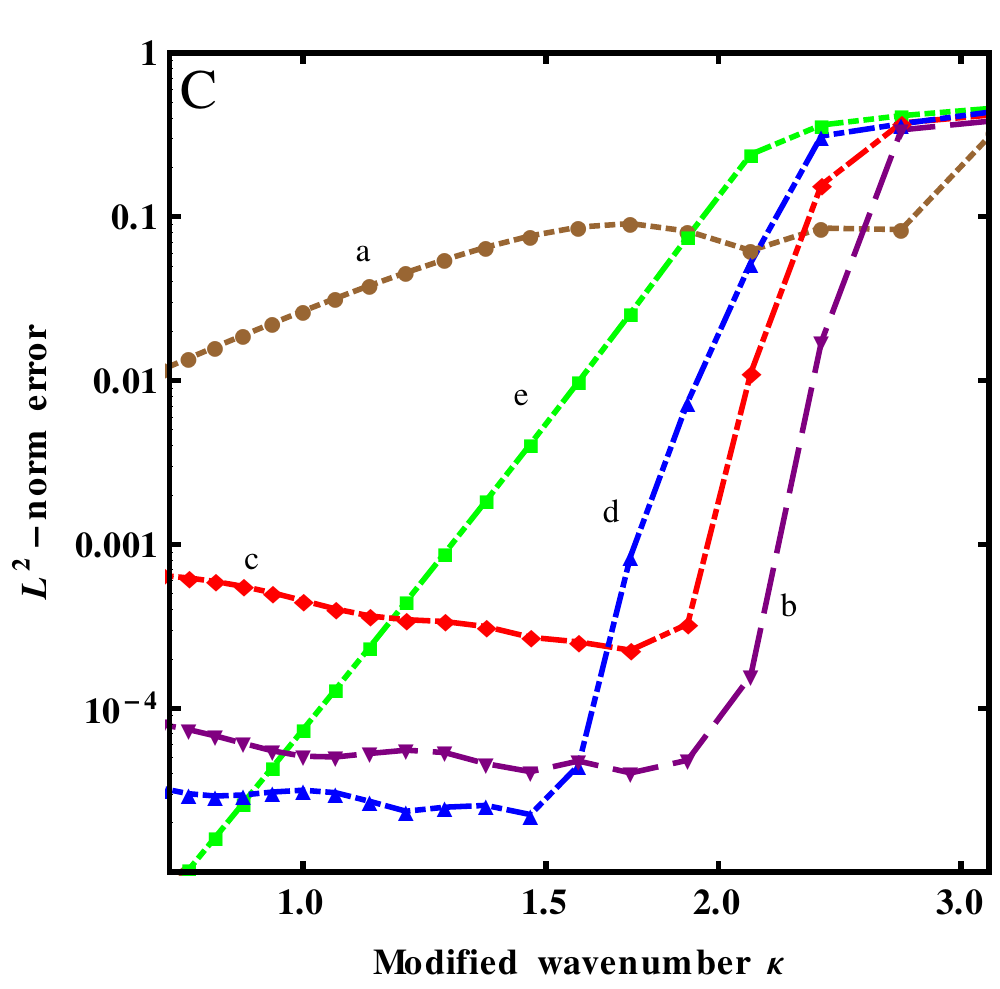}\hskip 0.5in
\includegraphics[scale=0.6]{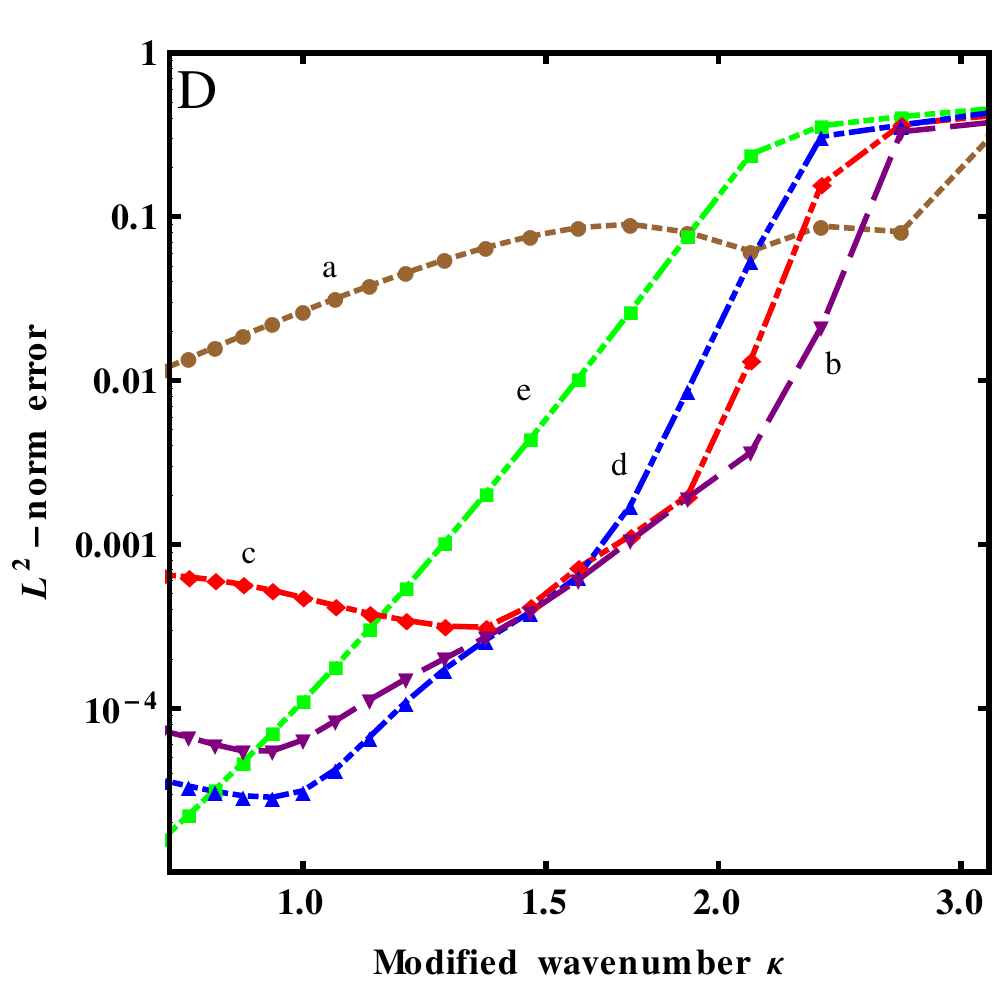}
\end{centering}
\caption{Numerical experiments for the schemes with pentadiagonal spatial discretizations. The time integrator is RK8. (A) Effective modified wavenumber. From top to bottom, (a) KIM4TH, (b) KLL2ND, (c) OPT2ND1P8, (d) OPT2ND1P5, and (e) UNOPT10TH. (B) The dispersion-and-dissipation error.
(C) The $L^2$-norm error with CFL number $r=0.1$. The effective wavenumber is against $k_{max}=380$. (D) The same as (C), but with CFL number $r=0.5$.}
\label{Fig:3_2}
\end{figure}

In this section, we compare several schemes with pentadiagonal spatial discretizations for the initial value problem (\ref{2_52}).  We plot the effective modified wavenumber of the tested schemes in Figure \ref{Fig:3_2}(A).  From top to bottom, the schemes are (a) KIM4TH, (b) KLL2ND, (c) OPT2ND1P8, (d) OPT2ND1P5, and (e) UNOPT10TH. The integration limit ($\Gamma$ value) for KLL2ND is $\Gamma =2.0$. The time integrator for the experiments is RK8.
Figure \ref{Fig:3_2}(B) plots the dispersion-and-dissipation error, $|\kappa^*-\kappa|$, for these schemes. We observe that for  the three DRB schemes developed in this paper, (b) KLL2ND, (c) OPT2ND1P8, and (d) OPT2ND1P5, their dispersion-and-dissipation errors are below $5\times 10^{-6}$ for the modified wavenumber smaller than the $\Gamma$ values used in their objective functions (2.0, 1.8, and 1,5, respectively). Figure \ref{Fig:3_2}(C) plots the $L^2$-norm errors for the tested scheme. The CFL number is $0.1$ for these experiments. It is worth noting that the $L^2$-norm error of KLL2ND starts decreasing around $\kappa\approx 2.7$, and remains small after $\kappa\approx 2.0$. This scheme has the potential to satisfy a larger range of thresholds discussed in Figure \ref{Fig:2_5}, compared with the other tested schemes. The two schemes, OPT2ND1P8, and OPT2ND1P5 have the same tendency as KLL2ND. Namely, the $L^2$-norm error decays and enters the plateau region around the modified wavenumber that coincides with the $\Gamma$ values used in their objective functions (1.8 and 1.5, respectively), thus the $\Gamma$ value determines the modified wavenumber where the $L^2$-norm error stops decaying and enters the plateau region. 

Similar to Figure \ref{Fig:3_2}(C), Figure \ref{Fig:3_2}(D) is a plot of the $L^2$-norm errors, except the CFL number is $r=0.5$.  When comparing Figure \ref{Fig:3_2}(C) \& (D), we see that for the schemes (b) KLL2ND, (c) OPT2ND1P8, (d) OPT2ND1P5, and (e) UNOPT10TH, when $\Delta t$ is not small enough ($r=0.5$),  the error from the time integrator is not negligible, whereas for (a) KIM4TH, the $L^2$-norm errors are very much the same for both $r=0.1$ and $r=0.5$. Looking at \eqref{2_41}, this suggests that for a fixed time-step $\Delta t$, when refining the spatial mesh (keeping $r < 1$), the time-integration error eventually becomes the dominant error for the schemes (b), (c), (d), and (e) ($r$ increases from 0.1 to 0.5 in this example), whereas the spatial-discretization error continues to be the dominant error for the scheme (a).
\begin{figure}[ht]
\begin{centering}
\includegraphics[scale=0.65]{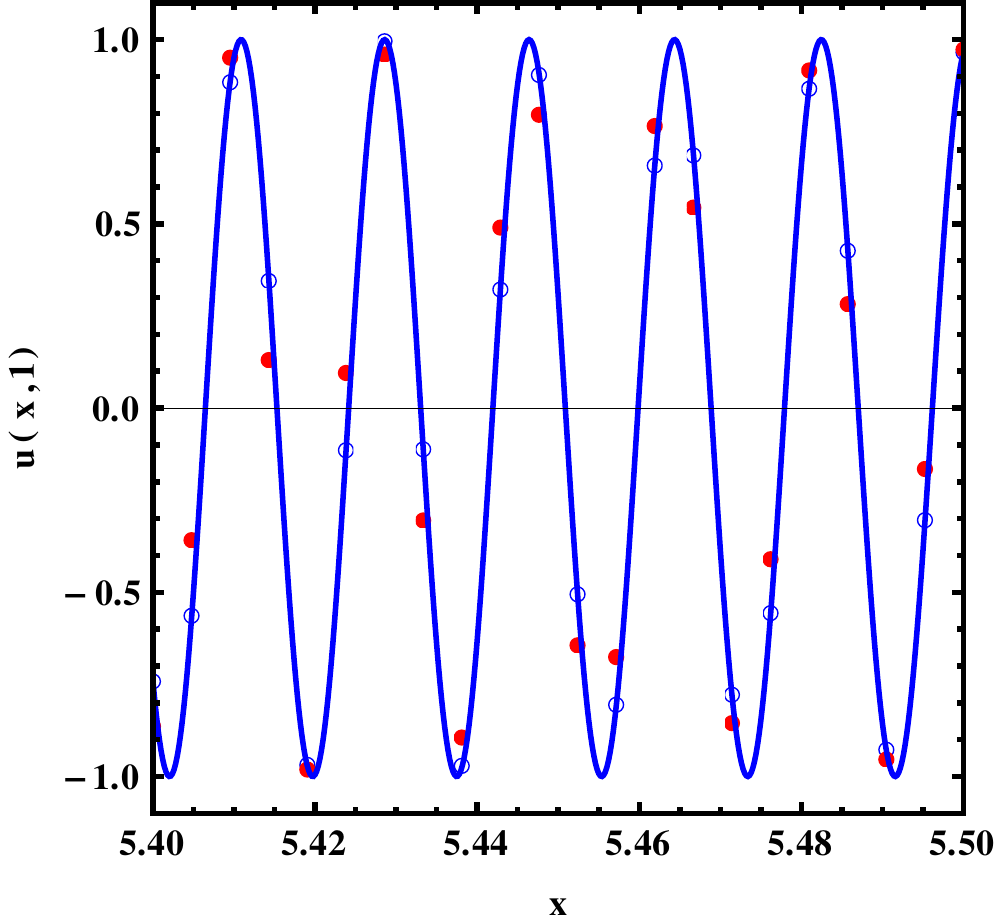}
\end{centering}
\caption{Comparison of numerical solutions of  KLL2ND and UNOPT10TH against the exact solution for the initial value problem (\ref{2_52}) at $T=1$. Solid line is the exact solution. Circle is KLL2ND and filled circle is UNOPT10TH. RK8 is the time integrator. $\Delta x= 1/210$ and CFL number is $r=0.1$.}
\label{Fig:3_3}
\end{figure}

We now revisit the initial value problem (\ref{2_52}). We compare the numerical solutions computed by UNOPT10TH and KLL2ND, respectively, against the exact solution. The mesh size for the calculation is $\Delta x=1/210$. The time integrator is RK8, and the CFL number is $r=0.1$.  Figure \ref{Fig:3_3} shows the solution waveform for $x\in [5.4, 5.5]$ at $T=1$. Since the speed of the transport equation is $c=1$, initially at $T=0$, the interval of $x$ for the current window is $[4.4, 4.5]$. Recall the instantaneous wavenumber $k_c(x)$ defined in the previous section. The instantaneous wavenumber at the midpoint of $[4.4, 4.5]$ is $k_c(4.45)\approx 350.89$. Hence locally the wavelength is $\lambda\approx 2\pi/k_c(4.45)\approx 0.0179$, as shown in Figure \ref{Fig:3_3}. Since FD numerical solution is affected not only by the local wavenumber, in fact the $L^2$-norm error of the FD solution has influence from the effective wavenumber as high as $\kappa_{max}=k_{\mathrm{max}}\Delta x\approx 1.8190$ in this example. If we look at Figure \ref{Fig:3_2}(C), at $\kappa_{\mathrm{max}}$ KLL2ND has much smaller $L^2$-norm error than UNOPT10TH. This reflects to the numerical waveforms in Figure \ref{Fig:3_3} as well. The solid line in Figure \ref{Fig:3_3} is the exact solution. The circle is the solution by using KLL2ND and the filled circle is the solution by using UNOPT10TH. We see that the KLL2ND solution visually has no phase shift, whereas the UNOPT10TH solution produces an obvious phase error.

\section{Non-constant coefficient transport equation}\label{sec:non-constant}

In this section, we extend the error analysis to a case with non-constant coefficients. Consider 
\begin{equation}\label{4_1}
u_t+c(x)u_x=0,\quad 0\leq x\leq\ell\,,
\end{equation}
with periodic boundary conditions. After an application of the Fourier transform, the equation becomes an infinite-dimensional system of ODEs,
\begin{equation}\label{4_2}
\dfrac{\dif}{\dif t}\widehat{u}(k',t) + (\widehat{c}*\widehat{u_x})(k',t)=0\,,\quad k'\in\mathbb{Z}\,,
\end{equation}
where $\widehat{u}(k',t)$ and $\widehat{c}(k')$ are the (spatial) Fourier coefficients of $u(x,t)$ and $c(x)$, respectively. 
The Fourier coefficient $\widehat{u_x}$ of the first derivative $u_x(x,t)$, corresponding to wavenumber $k'$, is given by $\frac{2\pi \mi k'}{\ell} \widehat{u}(k',t)$, and the convolution $(\widehat{c}*\widehat{u_x})(k',t)$ is defined by
\begin{equation}\label{4_3}
(\widehat{c}*\widehat{u_x})(k',t)\equiv\sum_{\alpha=-\infty}^{\infty}\widehat{c}(k'-\alpha)\cdot\frac{2\pi \mi \alpha}{\ell}\widehat{u}(\alpha,t)\,.
\end{equation}
In the previously studied case that $c(x)=c$ for some constant $c$, the corresponding Fourier coefficient is $c$ for $k=0$ and zero otherwise. Thus, in that case, the convolution collapses to a scalar multiplication, and the ODE system \eqref{4_2} becomes fully decoupled. We recall that the solution of that simpler problem is given in \eqref{2_8}. Inspection of that formula shows that the Fourier modes corresponding to different wavenumbers cannot interact or interchange energy with each other. 
Unfortunately, this is no longer true for general $c(x)$. In this case different modes can interact through the convolution, and the initial spectrum can migrate to different wavenumbers.

Considering a mesh grid with the grid points $x_j=j\Delta x$, $\Delta x=\ell/N$, $j=1,\:2,\cdots,\:N$ equally distributed,  we denote the approximation to $u(x_j,t)$ on the grid points by $U_j(t)$ and define $\textbf{U}(t)\equiv(U_1(t),\cdots,\:U_N(t))^T$. If we suppose that the first derivative $u_x$ is approximated by a compact FD scheme, $\textbf{L}\textbf{U}'=\textbf{R}\textbf{U}$, the PDE \eqref{4_1} can be written as a semi-discretized system
\begin{equation}\label{4_4}
\dfrac{\dif}{\dif t}\textbf{U}(t)+\textbf{CDU}=0\,,
\end{equation}
where $\textbf{D}=\textbf{L}^{-1}\textbf{R}$ and $\textbf{C}$ is a diagonal matrix whose $j$-th diagonal element is given by $c(x_j)$. Thus, in terms of the matrix exponential, the solution of \eqref{4_4} can be written as
\begin{equation}\label{4_5}
\textbf{U}(t)=\me^{-t\textbf{CD}}\textbf{U}(0)\,.
\end{equation}
Now, we may solve \eqref{4_4} by the RK method described in section \ref{sec:RK}, and all of the analysis in  \S\ref{sec:analysis} is still valid. The only difference is that the term $c\textbf{D}$ in the original formula is replaced by $\textbf{CD}$. As a result, in a single time step, the matrix exponential in \eqref{4_5} becomes $\me^{-\Delta t\textbf{CD}}$, which is approximated by a polynomial $\mathcal{P}(-\Delta t\textbf{CD})\equiv\sum_{m=0}^M a_m(-\Delta t\textbf{CD})^m$. Denoting the numerical approximation to $u(x_j,t_n)$ by $U^n_j$, where $t_n\equiv n\Delta t$ and letting $\textbf{U}^n\equiv(U^n_1,\cdots,U^n_N)^T$, we see that the numerical solution of the RK method for \eqref{4_4} is updated by 
\begin{equation}\label{4_6}
\textbf{U}^{n+1}=\mathcal{P}(-\Delta t\textbf{CD})\textbf{U}^n\,.
\end{equation}

Equation \eqref{4_6} would be useful for our error analysis if we had knowledge about the eigenvalues and the eigenvectors of the matrix $\textbf{CD}$. 
However, the availability of this information is strongly problem dependent. Another way to analyze the error is to look at the Discrete Fourier Transform of the semi-discretized system \eqref{4_4}. The DFT transforms \eqref{4_4} into the following finite-dimensional ODE system
\begin{equation}\label{4_7}
\dfrac{\dif}{\dif t}\widehat{\textbf{U}}(k',t) + (\widehat{\textbf{C}}*\widehat{\textbf{DU}})(k',t)=0\,,
\quad k'=-\dfrac{N}{2}+1,\cdots,\dfrac{N}{2}\,,
\end{equation}
where 
\begin{equation}\label{4_8}
\widehat{\textbf{C}}(k')=\dfrac{1}{N}\sum_{\alpha\equiv k'\pmod N}\widehat{c}(\alpha)\,,
\end{equation}
and
\begin{equation}\label{4_9}
\widehat{\textbf{DU}}(k',t) = \mi k^*\widehat{\textbf{U}}(k',t)\,.
\end{equation}
In this case the numerical wavenumber is
\begin{equation}\label{4_10}
k^*\equiv\dfrac{\kappa^*\left(\frac{2\pi k'\Delta x}{\ell}\right)}{\Delta x}\,,
\end{equation}
and 
\begin{equation}\label{4_11}
(\widehat{\textbf{C}}*\widehat{\textbf{DU}})(k',t)\equiv\sum_{\alpha=-N/2+1}^{N/2}\widehat{C}(k'-\alpha)\cdot \mi k^*\widehat{\textbf{U}}(k,t)\,.
\end{equation}
We note that \eqref{4_11} can be viewed as a linear transformation that maps the vector $\widehat{\textbf{U}}$ into the vector $(\widehat{\textbf{C}}*\widehat{\textbf{DU}})$. From this point of view, the ODE system is equivalent to 
\begin{equation}\label{4_12}
\dfrac{\dif}{\dif t}\widehat{\textbf{U}} + \widehat{\textbf{M}}\widehat{\textbf{U}}=0\,,
\end{equation}
which has a solution of the form 
\beq\label{4_12_1}
\widehat{\textbf{U}}(t)=\me^{-t\widehat{\textbf{M}}}\widehat{\textbf{U}}(0)\,. 
\eeq
But again, without a priori knowledge about $\widehat{\textbf{C}}$ and $\widehat{\textbf{M}}$, it is difficult to say more about the behavior of $\widehat{\textbf{U}}(t)$. Nevertheless, from \eqref{4_12_1}, we expect that the spectrum $\widehat{\textbf{U}}(t)$ migrates towards different wavenumbers exponentially quickly. This is, in fact, the primary source of numerical error (aliasing error), and the dispersion-and-dissipation error becomes secondary. 

\begin{example}[Non-constant coefficients]\label{ex:ncc}
We consider the initial value problem (\ref{4_1}), with $\ell=1$ and the coefficient 
\begin{equation}\label{4_13}
c(x)=A + B\sin^2(2\pi x)\,,
\end{equation}
and the initial condition 
\begin{equation}\label{4_14}
u_0(x)=\cos(30\pi x)e^{-80(x-0.5)^2}\,.
\end{equation}
If $A$ and $B$ are chosen so that $c(x)$ does not change sign in the periodic domain, we can solve the PDE analytically by the method of characteristics. The solution is a time-periodic function with period
\[
T_p\equiv\dfrac{1}{\vert A\vert\sqrt{1+B/A}}\,.
\]
Consequently, the corresponding infinite-dimensional ODE system in the Fourier space for this example admits only time-periodic solutions with the same period $T_p$.  Also, since $2\sin^2(2\pi x)=1-\cos(4\pi x)$, the transportation speed $c(x)$ is a periodic function with the period equal to 1/2. Thus, 
\beq\label{eq:uperiod}
u\left(x-\dfrac{1}{2},t\right) = u\left(x,t+\dfrac{T_p}{2}\right)\,.
\eeq
For our numerical example, we let $A = 1/5$ and $B =1$.
By virtue of the half-angle formula, we have 
\begin{equation}\label{4_16}
c(x)=\dfrac{1}{5}+\sin^2(2\pi x)=\dfrac{7}{10}-\dfrac{1}{4}\left(\me^{4\pi \mi x}+\me^{-4\pi \mi x} \right).
\end{equation}
This suggests that $\widehat{c}(k)$ is nonzero except $k=0$, $\pm2$. Therefore, after we introduce a spatial discretization to the PDE, the matrix $\widehat{\textbf{M}}$ in \eqref{4_12} is a finite dimensional cyclic pentadiagonal matrix with zero superdiagonal and subdiagonal. 
The waveform and spectrum of $u_0(x)$ are shown in Figure \ref{Fig:2_1}. The spectrum $\widehat{u}_0(k)$ decays exponentially fast as $\vert k\vert$ becomes large. Therefore, it is reasonable to assume that the spectrum $\widehat{u}_0(k)$ is nonzero only for the wavenumber $|k| \le 40\pi$.


From \eqref{eq:uperiod}, we know that the modulus  $\vert\widehat{u}(k,t)\vert$ is a time-periodic function with the period $T_p/2$. The spectrum $\vert\widehat{u}(k,t)\vert$, at $t=0$ and $t=T_p/4$, are shown in Figure \ref{Fig:4_1}(A). We see that at time $t=T_p/4$, the spectrum migrates to higher wavenumber. While the spectrum is confined within a finite range all the time, this range, approximately $\vert k\vert\leq200\pi$, is much larger than what we would expect from our knowledge about the initial condition ($|k| \le 40\pi$). Compared with the initial condition whose spectrum has a more well-defined cut-off wavenumber, the corresponding spectrum at $t=T_p/4$ decays much slower as the wavenumber increases, and the cut-off wavenumber is less well-defined.

Now we solve the problem numerically by three spatial discretization schemes: pseudospectral,  UNOPT10TH, and  KLL2ND. In each case the time integrator is RK8 with CFL number $r=0.1$. 
\begin{figure}[ht]
\begin{centering}
\includegraphics[scale=0.6]{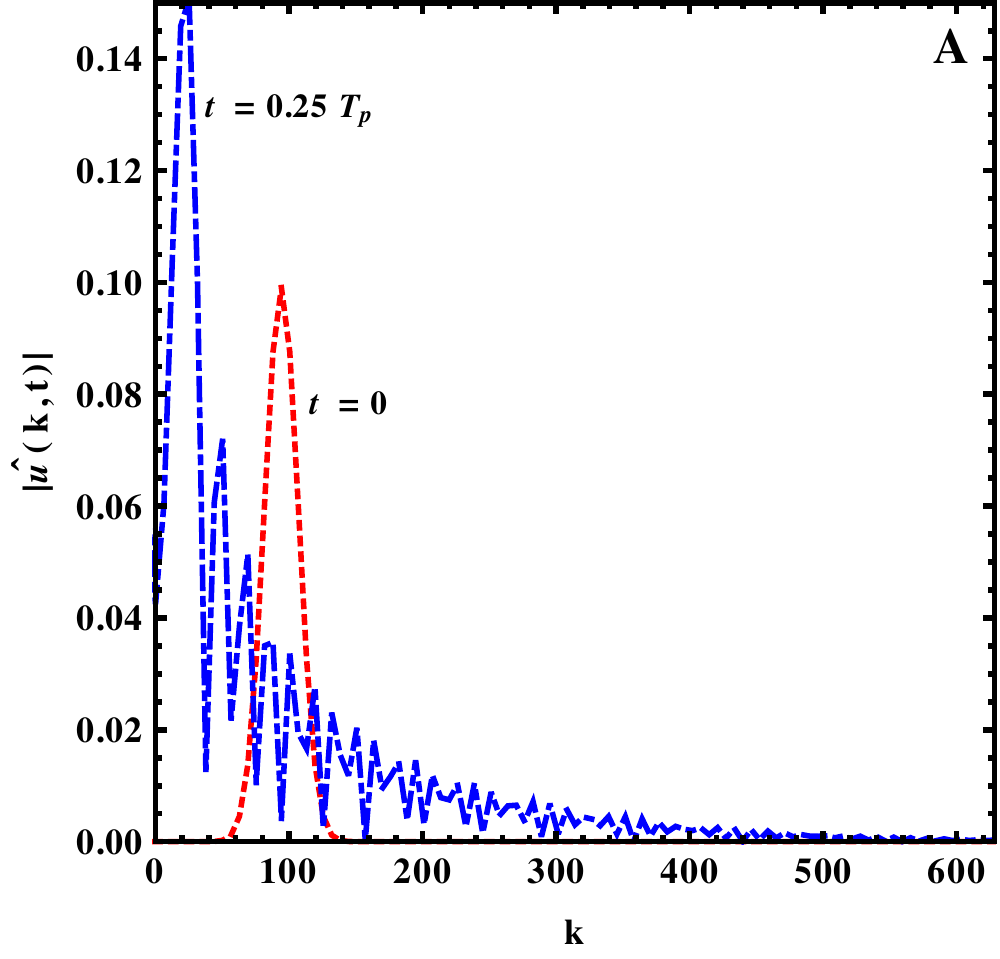}\hskip 0.5in
\includegraphics[scale=0.61]{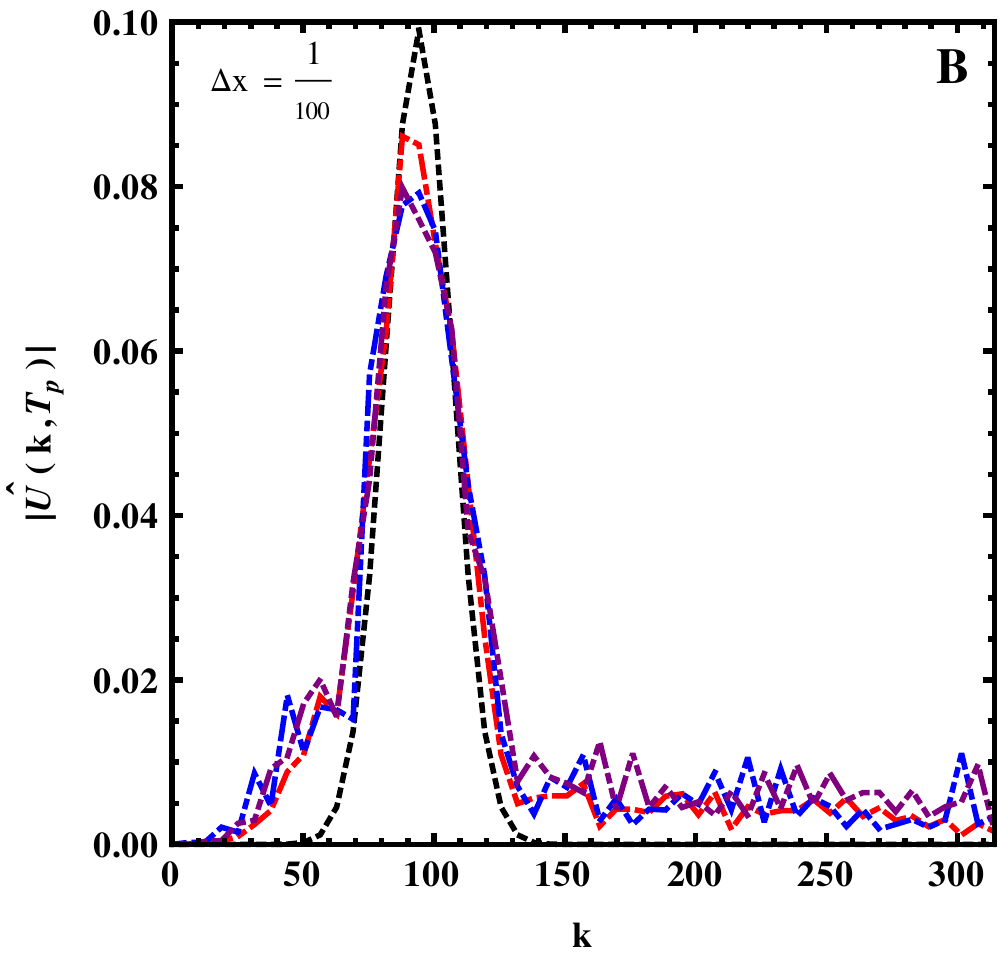}\\ \vskip 0.25in
\includegraphics[scale=0.6]{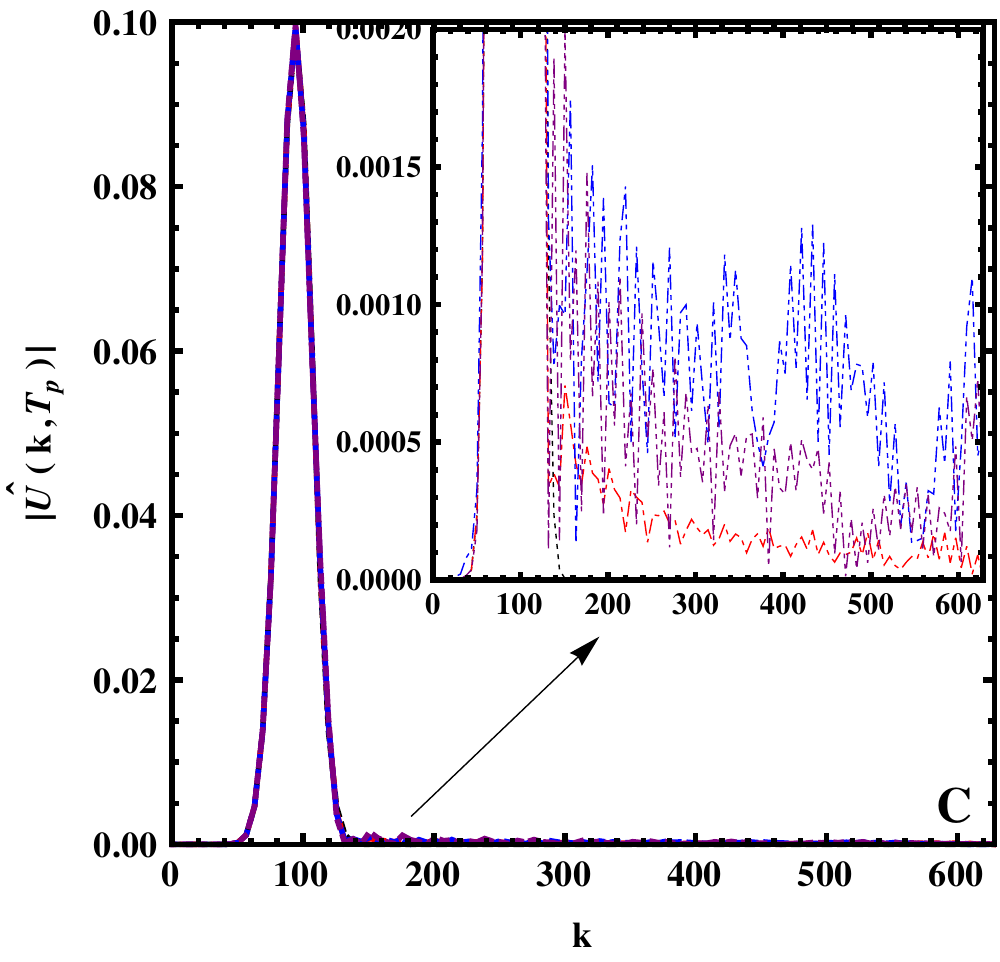}\hskip 0.5in
\includegraphics[scale=0.6]{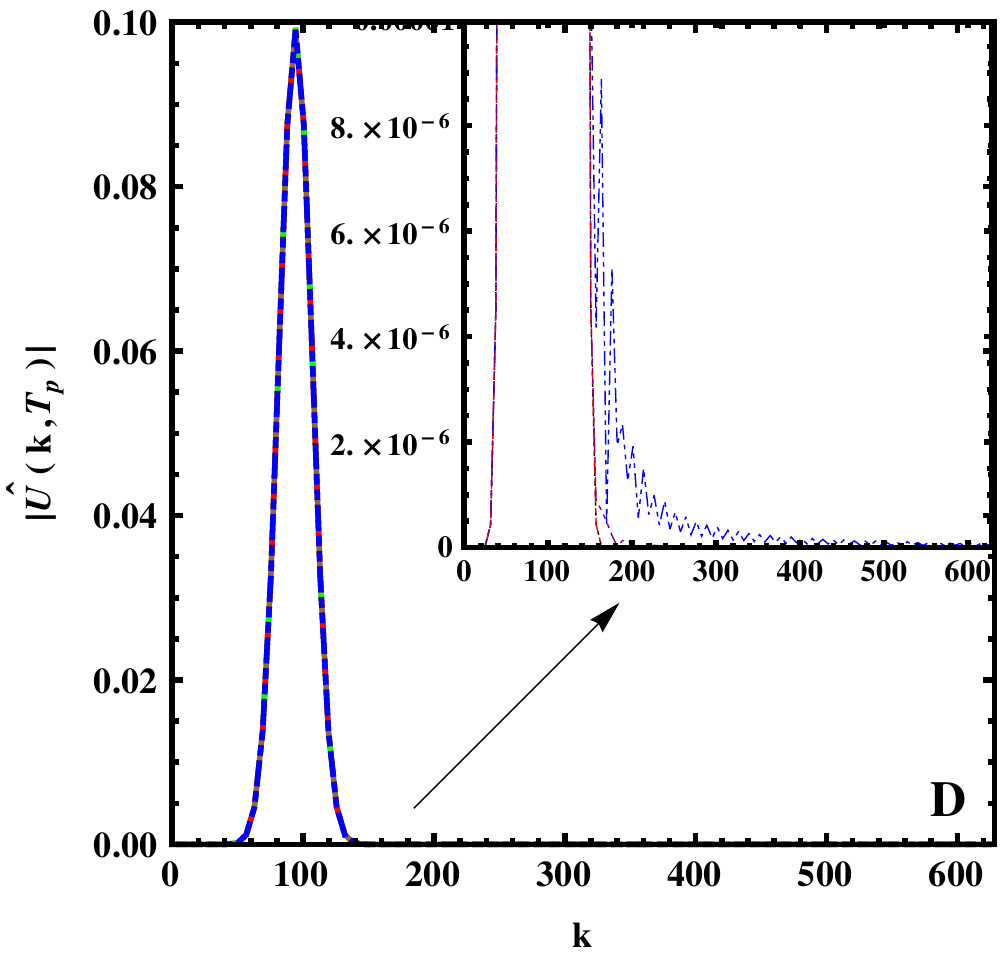}
\end{centering}
\caption{(A) The spectrum $\vert\widehat{u}(k,t)\vert$, at $t=0$ (dotted) and $t=T_p/4$ (dash-dotted. (B) Severe aliasing errors for numerical methods: the pseudostectral method,  KLL2ND, and UNOPT10TH (from top to bottom at $k=30\pi$. $\Delta x=1/00$. (C) The same as (B), but $\Delta x =1/200$. The dispersion-and-dissipation error starts to show it significance. (D) The same as (B),  but $\Delta x =1/500$. The aliasing error is almost negligible. The dispersion-and-dissipation error is the dominant error.} 
\label{Fig:4_1}
\end{figure}
Consider a mesh grid with $\Delta x=1/100$. The spectra of the numerical solutions for the three spatial discretizations at $t=T_p$ are shown in Figure \ref{Fig:4_1}(B), where around $k=30\pi$, from top to bottom are true spectrum,  the pseudostectral method,  KLL2ND, and UNOPT10TH (color version, black dotted: true, red dash-dotted: pseudostectral, purple dash-dotted: KLL2ND, blue dash-double-dotted: UNOPT10TH).  Note that by the Nyquist-Shannon sampling theorem, the highest resolvable wavenumber (Nyquist limit) provided by this mesh size is $2\pi/\ell\times 50=100\pi$, which is far smaller than what we need  (at least $200\pi$). As time increases from $0$ to $T_p/4$, the spectrum starts to migrate toward higher wavenumbers. When the migrating spectrum reaches $100\pi$, aliasing error occurs. At this point, the numerical solution has significant error, regardless of the spatial discretization scheme. 

We repeat the calculation for the mesh sizes $\Delta x=1/200$ and $\Delta x=1/500$. Figure \ref{Fig:4_1}(C) shows the spectra for $\Delta x=1/200$ and Figure \ref{Fig:4_1}(D) is for $\Delta x=1/500$. The Nyquist limit for $\Delta x=1/200$ is $200\pi$, while for $\Delta x=1/500$ is $500\pi$. Since $200\pi$ covers the widest range of spectrum for the solution (occurs at $t=T_p/4$),  the dispersion-and-dissipation error begins to show its significance at mesh size $\Delta x=1/200$. Figure  \ref{Fig:4_1}(C) shows that while the spectral method has the smallest error, KLL2ND produces smaller error than UPOPT10TH after $k >150$, for $\Delta x=1/200$. For $\Delta x=1/500$, the aliasing error is almost negligible, and the dispersion-and-dissipation error becomes the dominant error. Figure  \ref{Fig:4_1}(D) shows that in the frequency domain, for $\Delta x=1/500$, the true solution, KLL2ND, and the pseudostectral are visually indistinguishable, while UNOPT10TH produces considerable discrepancy. Similarly, while not immediately apparent here, our calculations indicate that in the physical domain, UNOPT10TH has $L^{\infty}$-error (maximum error) in the order of $O(10^{-6})$ for the solution of the PDE while the other two methods have much smaller errors.
\end{example}
\begin{note}
Finally, we remark that for certain operational simulations with periodic boundary conditions, such as the general circulation model (GCM) in atmospheric sciences, spectral methods are commonly used. 
Nevertheless,  based on the above numerical experiment, an optimized compact DRB FD scheme is capable of providing comparably accurate solutions within the same memory requirement, but with much less CPU time. Hence, optimized compact DRB FD schemes may be more economical for such simulations.
\end{note}


\section{Nonlinearity: a case-study}\label{sec:nonlinear}
We have analyzed the error dynamics of optimized compact DRB schemes through the linear transport equation with both constant and non-constant coefficients. It was shown by Hixon \cite{hixon00} that the benefits of the optimized schemes derived from the constant-coefficient linear transport equation have the potential to carry over to nonlinear applications. However, it is not clear that the error dynamics we observed in the previous examples continue to be valid for nonlinear cases. To address this question, in this section, we consider the nonlinear Hopf equation,
\begin{equation}\label{6_1}
u_t+(u^2/2)_x=0\,,\quad 0\leq x\leq\ell\,,
\end{equation}
in a periodic domain. For simplicity, we assume $\ell=1$ and restrict ourselves to the following initial condition suggested by Hixon \cite{hixon00}:
\begin{equation}\label{6_2}
u(x,0)=u_0(x)=\dfrac{3}{4}+\dfrac{1}{4}\sin(2\pi x)\,.
\end{equation}
The exact solution of this problem can be found by the method of characteristics. Shock formation (wave-breaking) for this initial condition occurs at $t=t_b=2/\pi\approx0.6366$. Our numerical simulations will be terminated right before the wave-breaking time $t_b$. It is known that the spectrum of the solution to the Hopf equation continues to extend towards higher wavenumber, in particular when time is close to $t_b$ (similar phenomenon was observed previously in our numerical experiment for the linear transport equation with non-constant coefficients in Example \ref{ex:ncc}). To see this behavior in Fourier space, we apply the Fourier transform so that \eqref{6_1} becomes 
\begin{equation}\label{6_4}
\dfrac{\dif}{\dif t}\widehat{u}(k,t) + (\widehat{u}*\widehat{u_x})(k,t)=0\,,
\quad k\in\mathbb{Z}\,.
\end{equation}
The presence of convolution allows interaction of different Fourier modes. Consequently, the spectrum of the solution varies with time. Furthermore, the nonlinearity may accelerate the migration of spectrum toward to the high wavenumber modes as time increases.

We now solve the nonlinear problem \eqref{6_1}, \eqref{6_2} by four spatial discretization schemes,  namely, KLL2ND, OPT2ND1P5, OPT2ND1P8, and UNOPT10TH. The time integrator is the MATLAB ODE solver, ODE45, with the parameters  RelTol = $10^{-9}$ and AbsTol = $10^{-10}$. The true solution of the Hopf equation at the grid points is obtained by solving a nonlinear equation resulting from the method of characteristics. 
\begin{figure}[ht]
\begin{centering}
\includegraphics[scale=0.45]{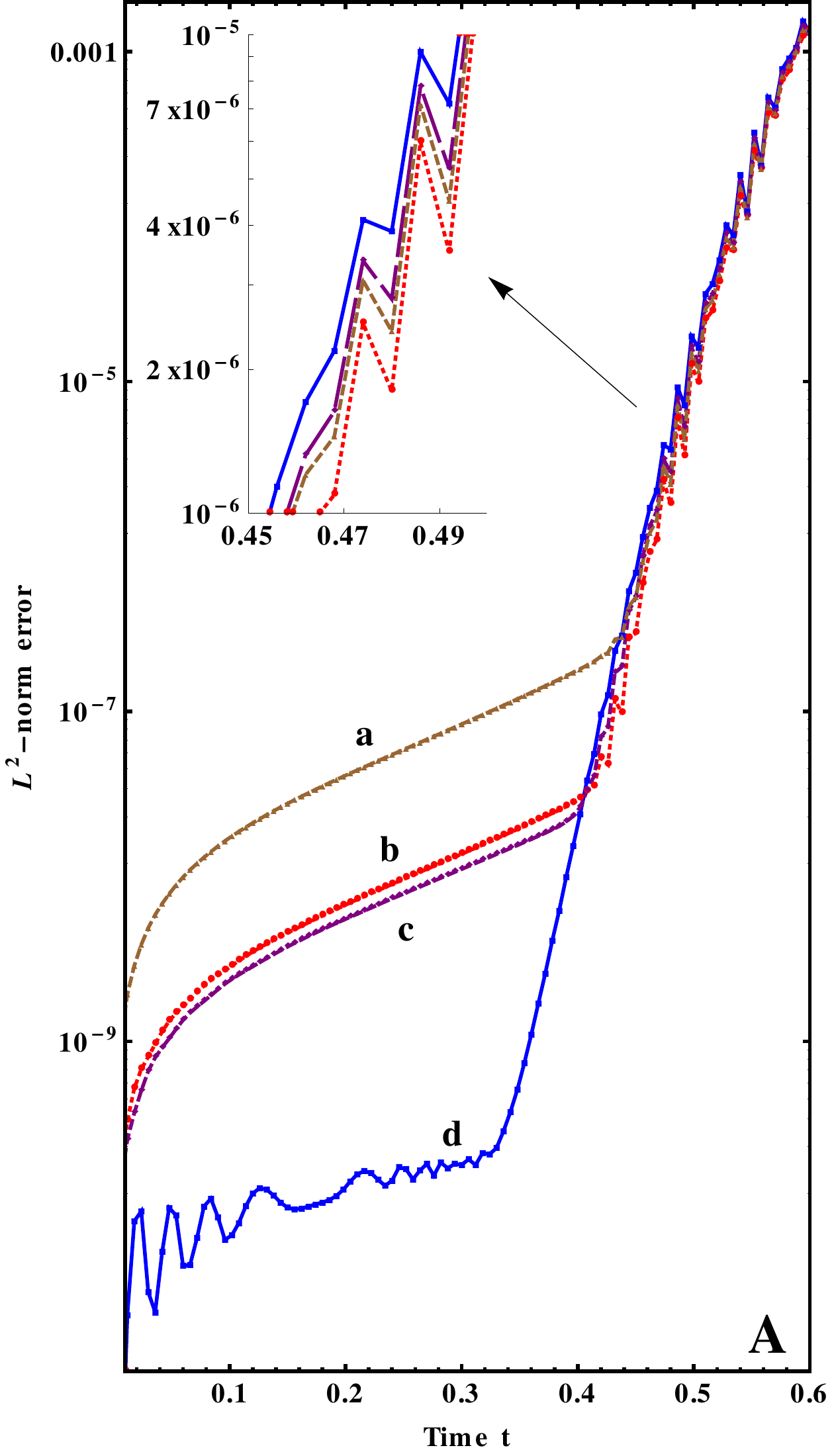}\hskip 0.5in
\includegraphics[scale=0.45]{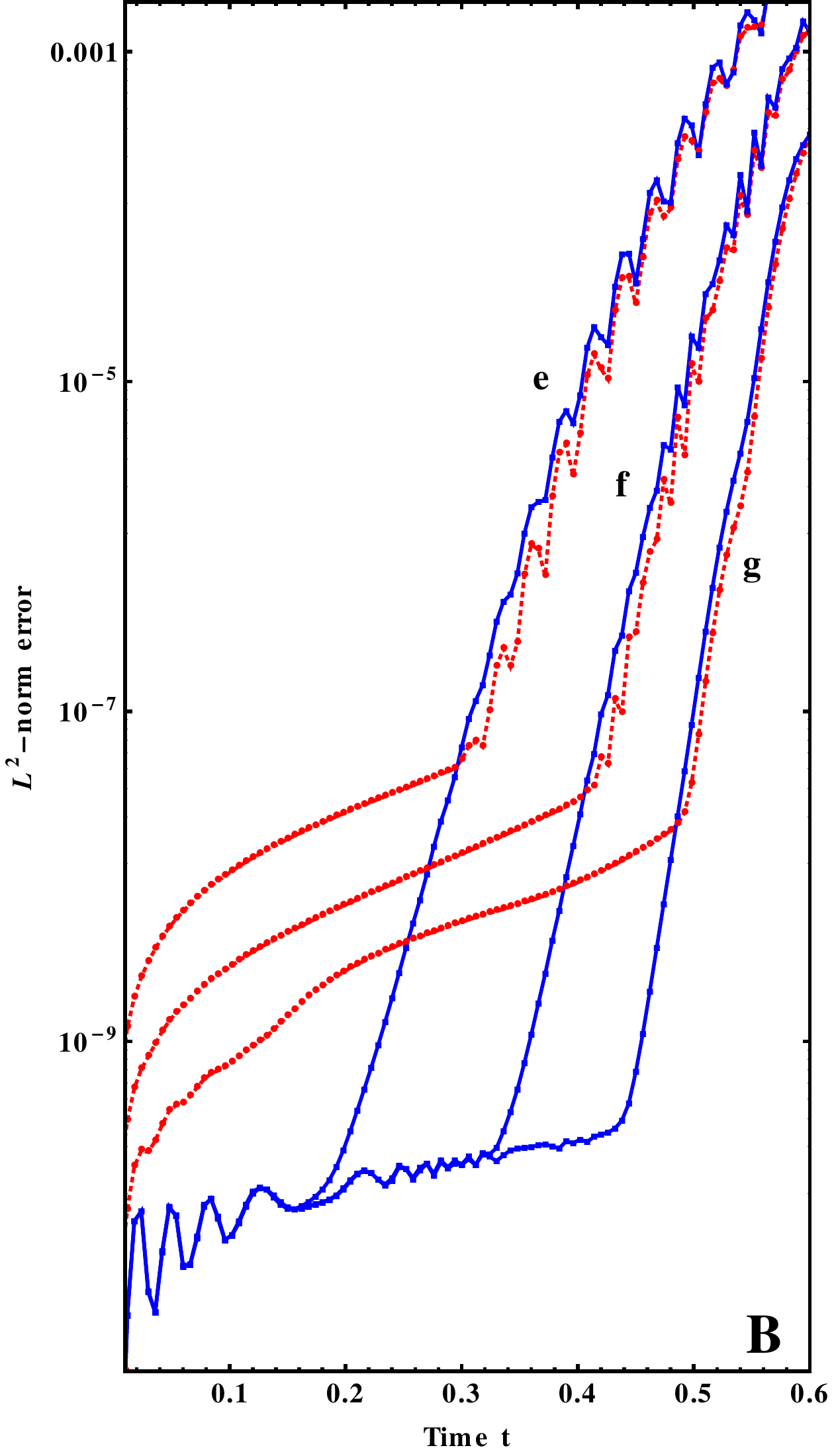}
\end{centering}
\caption{(A) The $L^2$-norm error against time $t$. The spatial discretization schemes are (a) OPT2ND1P8, (b) KLL2ND, (c) OPT2ND1P5 and (d) UNOPT10TH. The zoom-in figure shows the error at high wavenumber for $0.45 < t < 0.5$, from top to bottom, UNOPT10TH, OPT2ND1P5, OPT2ND1P8 and KLL2ND. The spatial grid size is $\Delta x=1/100$ in this example. (B) The grid refinement study for UNOPT10TH (solid line) and KLL2ND (dotted line). The spatial grid sizes are   (e) $\Delta x=1/50$, (f) $\Delta x=1/100$ and (g) $\Delta x=1/200$.}
\label{Fig:6_1}
\end{figure}
The spectrum of the initial condition (\ref{6_2}) is localized at $k=0$ and $k=1$. This is the smallest obtainable spectral bin number (wavenumber). Therefore, from our analysis we expect that---at least initially---the unoptimized high order-of-accuracy method, UNOPT10TH, will have the smallest $L^{2}$-norm error among the four tested methods. Figure \ref{Fig:6_1} (A) confirms that before $t=0.3$, (d) UNOPT10TH has the best performance. The other three optimized compact DRB schemes behave very much as in the linear case (see Figure \ref{Fig:3_2}(C) for comparison). For $0.3 < t < t_b$, the $L^{2}$-norm errors take a sharp turn for all four methods. This is the region where the spectrum accelerates into the high Fourier modes before the wave breaking. This is also the transient region where the dominant numerical dispersion-and-dissipation error moves from the low modified wavenumber region to the high modified wavenumber region. It is not too surprising that the optimized DRB scheme, KLL2ND, produces the smallest $L^2$-norm error, because this scheme was designed to reduce the numerical dispersion-and-dissipation error at high modified wavenumber in the linear case, as seen in Figure \ref{Fig:3_2}(C). The spatial grid size is $\Delta x=1/100$ in this example. 

Figure \ref{Fig:6_1}(B) shows the $L^2$-norm errors for UNOPT10TH and KLL2ND.  Three different grid sizes, (e) $\Delta x=1/50$, (f) $\Delta x=1/100$ and (g) $\Delta x=1/200$, are used. We see that before $t=0.2$, when refining grid, the error of KLL2ND is reduced. This truncation effect happens at very low modified wavenumber for KLL2ND.  For UNOPT10TH, the error stays the same, and this suggests that the error is dominated by the time integrator. After $t>0.2$, both methods start to take a turn to a steeper error. We observe that coarser grid takes an earlier turn. By definition, for a fixed wavenumber $k$, the modified wavenumber ($\kappa=k\Delta x$) is larger if $\Delta x$ is larger.  Hence the error of coarser grid taking a turn at an earlier time suggests that when the spectrum of numerical solution extends to the region of high wavenumber, coarser grid feels the effect of dispersion-and-dissipation error at high modified wavenumber earlier. This also explains why after taking the turn, the optimized DRB scheme has better performance than the unoptimized one (the improvement is about 10\% - 30\%), especially when a coarser grid is used. Basically, this numerical experiment confirms, consistent with the numerical experiments and observation previously reported by Hixon \cite{hixon00}, that the error dynamics and the benefits of the DRB compact schemes we observe from the linear case do carry over into nonlinear applications. We, however, provide the connection between the linear and nonlinear cases from the view point of numerical dispersion-and-dissipation error. 

\section{Higher-dimensional cases}\label{sec:2D}
In this section we indicate the extension of  our analysis to higher dimensional problems. We start by considering the two-dimensional transport equation 
\begin{equation}\label{7_1}
u_t+c_x u_x+c_y u_y=0\,,\quad 0\leq x,\:y\leq\ell\,,
\end{equation}
where $c_x$ and $c_y$ are constants. For simplicity, we assume periodic boundary conditions.  We take the Fourier transform with respect to the spatial variables $(x,y)$, and we find 
\begin{equation}\label{7_2}
\dfrac{\dif}{\dif t}\widehat{u}(k_x',k_y',t)=-\mi\dfrac{2\pi}{\ell}(c_x k_x' + c_y k_y')\widehat{u}(k_x',k_y',t)\,,
\quad k_x',\:k_y'\in\mathbb{Z}\,,
\end{equation}
which admits the solution
\begin{equation}\label{7_3}
u(x,y,t)=\sum_{k_x',\:k_y'\in\mathbb{Z}}\widehat{u}(k_x',k_y',0)\me^{-\mi(c_xk_x+c_yk_y)t}\me^{\mi k_x x}\me^{\mi k_y y}\,.
\end{equation}
Here, $k_x'$ and $k_y'$ are the wavenumbers in the $x$- and $y$-directions, respectively, and $k_x=2\pi k_x'/\ell$ and $k_y=2\pi k_y'/\ell$ are the corresponding angular wavenumbers. 

Let $x_i=i\Delta x$, $y_j=j\Delta y$, $\Delta x=\Delta y=\ell/N$, and $U_{ij}(t)\approx u(x_i,y_j,t)$, ${\bf U}(t)=\left[U_{ij}(t)\right]$. Define the matrix ${\bf D}_y={\bf L^{-1}R}$, where ${\bf L}$ and ${\bf R}$ are defined as in \eqref{2_16}. Let ${\bf D}_x={\bf D}_y^T$ (transpose of ${\bf D}_y)$. In matrix form, the partial derivatives of ${\bf U}$ with respect to $x$ and $y$ are approximated by ${\bf U}{\bf D}_x$ and ${\bf D}_y{\bf U}$. Therefore, the semi-discretized version of \eqref{7_1} may be written as
\begin{equation}\label{7_4}
\dfrac{\dif}{\dif t}{\bf U}(t) = - \left( c_x {\bf U}{\bf D}_x + c_y {\bf D}_y{\bf U} \right)\,.
\end{equation}
The right-hand-side of \eqref{7_4} is a linear transformation of ${\bf U}$, and the associated eigenvectors are given by 
\begin{equation}\label{7_5}
U_{ij}=\me^{\mi k_x x_i}\me^{\mi k_y y_j}\,,\quad 1\leq i,\:j\leq N\,,
\end{equation}
for which the corresponding eigenvalues are 
\begin{equation}\label{7_6}
\omega(k_x, k_y)=-\mi(c_xk_x^*+c_yk_y^*)\,,
\end{equation}
where $k_x^*$ and $k_y^*$  are $x$- and $y$- component of the numerical wavenumber, respectively.  Similarly as in the one-dimensional case, we define $k_x^*\equiv\kappa_x^*(k_x\Delta x)/\Delta x$, $k_y^*\equiv\kappa_y^*(k_y\Delta y)/\Delta y$ to be the components of effective wavenumber associated with $k_x$ and $k_y$, respectively. The effective modified wavenumber in two-dimension is a vector function $\kappa^*=(\kappa^*_x, \kappa^*_y)$. For simplicity, we assume that $\kappa_x^*=\kappa^*_y=\kappa^*$, where $\kappa^*$ is defined by Equation \eqref{2_20}. 

Taking DFT, we see that \eqref{7_4} is transformed into the decoupled ODE system
\begin{equation}\label{7_7}
\dfrac{\dif}{\dif t}\widehat{\bf U}(k_x',k_y',t)=-\mi(c_xk_x^*+c_yk_y^*)\widehat{\bf U}(k_x',k_y',t)\,,\quad -\dfrac{N}{2}+1\leq k_x',\:k_y'\leq\dfrac{N}{2}\,.
\end{equation}
Solving the above system of equations, we obtain the solution to the semi-discretized system
\begin{equation}\label{7_8}
U_{ij}(T)=\sum_{-N/2+1\leq k_x',\:k_y'\leq N/2}\widehat{\bf U}(k_x',k_y',0)\me^{-\mi(c_xk_x^*+c_yk_y^*)T}\me^{\mi k_x x_i}\me^{\mi k_y y_j}\,,
\end{equation}
for $1\leq i,\:j\leq N$.

Suppose we use the $M$-stage RK method  \eqref{2_30} as the time integrator for \eqref{7_7}. The numerical solution can be written down with the help of 
\begin{equation}\label{7_9}
\mathcal{P}\left(-\mi(c_xk_x^*+c_yk_y^*)\Delta t\right)\equiv\sum_{m=0}^M a_m\left(-\mi(c_xk_x^*+c_yk_y^*)\Delta t\right)^m,
\end{equation}
where the polynomial $\mathcal{P}(\cdot)$ is defined by \eqref{2_31}. Denoting $\widehat{\bf U}^n=\widehat{\bf U}(n\Delta t)$, we find that at time $t=(n+1)\Delta t$, the numerical solution of \eqref{7_7} may be written as
\begin{equation}\label{7_10}
\widehat{\bf U}^{n+1}(k_x',k_y')=\mathcal{P}\left(-\mi(c_xk_x^*+c_yk_y^*)\Delta t\right)\widehat{\bf U}^{n}(k_x',k_y'),\quad -\dfrac{N}{2}+1\leq k_x',\:k_y'\leq\dfrac{N}{2}\,.
\end{equation}
We define the effective temporal frequency to be
\begin{equation}\label{7_11}
\omega^*\equiv\dfrac{\mi}{\Delta t}\log\mathcal{P}\left(-\mi(c_xk_x^*+c_yk_y^*)\Delta t\right)\\
=(c_xk_x^*+c_yk_y^*)+O\left((c_xk_x^*+c_yk_y^*)^{p+1}\Delta t^p\right),
\end{equation}
assuming the RK integrator is $p$-th order accurate.  Finally, we can write down an explicit expression for the numerical solution of the two-dimensional transport equation as
\begin{equation}\label{7_12}
U_{ij}^n=\sum_{-N/2+1\leq k_x',\:k_y'\leq N/2}\widehat{\bf U}(k_x',k_y',0)\me^{-\mi\omega^*T}\me^{\mi k_x x_i}\me^{\mi k_y y_j}\,,
\end{equation}  
where $T=n\Delta t$, and the phase error associated with the wavenumber vector $(k_x,k_y)$ is 
\begin{equation}\label{7_13}
\theta_e=\left(c_x(k_x^*-k_x)+c_y(k_y^*-c_y)\right)T + O\left(T(c_xk_x^*+c_yk_y^*)^{p+1}\Delta t^p\right)\,.
\end{equation}
Evidently, the above equation is the two-dimensional extension of \eqref{2_41}.  In the limiting case $\Delta t\rightarrow0^+$, the phase error reduces to 
\begin{equation}\label{7_14}
\theta_e\approx\left(c_x(k_x^*-k_x)+c_y(k_y^*-k_y)\right)T\,.
\end{equation}


\begin{example}[A two-dimensional example]\label{ex:2d}
We consider a problem that is similar to the test problem of Zalesak's rotating disk \cite{zalesak}. The governing equation is
\begin{equation}\label{7_15} 
u_t+c_x(x,y)u_x+c_y(x,y)u_y=0\,,
\quad 
-\dfrac{1}{2}\leq x,\:y\leq\dfrac{1}{2}\,,
\end{equation} 
where 
\begin{subequations}\label{7_16}
\begin{align}
c_x(x,y)&=2\pi y\,,\\
c_y(x,y)&=-2\pi x\,.
\end{align}
\end{subequations}
We assume that the support of the initial condition is confined within the disk $D=\left\lbrace(x,y)\;:\;x^2+y^2<1/4\right\rbrace$. Then the true solution at time $t$ is given by rotating the initial profile by an angle of $2\pi t$ radians. That is, if the initial condition is $u_0(x,y)$, we may write the solution as
\begin{equation}\label{7_17}
u(x,y,t)=u_0(x\cos\theta-y\sin\theta,x\sin\theta+y\cos\theta)\,,
\quad 
(x,y)\in D\,,
\end{equation} 
where $\theta=2\pi t$. The initial condition for our numerical experiments is  
\begin{equation}\label{7_20}
u_0(x,y)=\left(1-u_1(x,y)\right)u_2(x,y)\,,
\end{equation}
where
\begin{subequations}\label{7_21}
\begin{align}
u_1(x,y)&\equiv H\left((x^2+\left(y-\dfrac{1}{12}\right)^2-\dfrac{1}{20},\dfrac{1}{20}\right),\\
u_2(x,y)&\equiv 1-\left(1-H\left(x^2-\dfrac{1}{100},\dfrac{1}{40}\right)\right)\left(1-H\left(\left(y+\dfrac{1}{8}\right)^2-\dfrac{1}{36},\dfrac{1}{40}\right)\right)\,,\\
\intertext{and}
H(z,\delta)&\equiv\dfrac{1}{2}\left(1+\erf\left(\frac{4z}{\delta}\right)\right)\,.
\end{align}
\end{subequations}
In \eqref{7_21}, $\erf(\cdot)$ denotes the error function. Here the function $H(z,\delta)$ is defined so that $H\approx0$ for $z<-\delta/2$ and $H\approx1$ for $z>\delta/2$. With the function $H$, the support of the initial condition is confined within the disk $D$. The initial condition is shown in Figure \ref{Fig:6_1} (left panel).

Before we perform our numerical experiments, we first examine the nature of \eqref{7_15} and the initial condition \eqref{7_20}. An application of the Fourier transform to \eqref{7_15} yields
\begin{equation}\label{7_18}
\dfrac{\dif}{\dif t}\widehat{u}(k_x',k_y',t) + (\widehat{c_x}*\widehat{u_x})(k_x',k_y',t) + (\widehat{c_y}*\widehat{u_y})(k_x',k_y',t) = 0\,,
\quad k_x',\:k_y'\in\mathbb{Z}\,.
\end{equation} 
In this case the two-dimensional convolution is defined by
\begin{equation}\label{7_19}
(\widehat{f}*\widehat{g})(k_x',k_y')\equiv\sum_{k_x'',\:k_y''\in\mathbb{Z}}\widehat{f}(k_x'-k_x'',k_y'-k_y'')\widehat{g}(k_x'',k_y'')\,.
\end{equation}
A feature of the exact solution \eqref{7_17} is that, although the convolution implies that Fourier modes may interact with each other, in this case the interaction is confined to modes with the same ``length" $\sqrt{k_x'^2+k_y'^2}$. In fact,  a Fourier mode in two dimensions corresponds to a wavenumber vector $(k_x',k_y')$, the wave propagates along a direction that is parallel to the vector. Along this direction, we can specify a single wavenumber $k=\sqrt{k_x'^2+k_y'^2}$, exactly like what we did for the one-dimensional case. Therefore, although the spectrum of the solution does migrate in the wavenumber space, it does not migrate towards to higher wavenumber. Instead, it rotates in the wavenumber space. Now we take the Fourier Transform for $u_0$.  The maximal value of $|\widehat{u}_0|$ is of $O(1)$, which occurs at $(k_x',k_y')=(0,0)$. The value of $|\widehat{u}_0(k_x',k_y')|$ is of $O(10^{-6})$ for $k=\sqrt{k_x'^2+k_y'^2}>50$. We remark that since the wavenumber is localized at 50, according to \eqref{7_13} or \eqref{7_14}, in principle the relation between the error and the modified wavenumbers (or the grid sizes $\Delta x$ and $\Delta y$) is predictable. In practice, however, because of the effect of the convolution in \eqref{7_18}, to achieve a desired accuracy we may need smaller grid sizes than the analysis suggests (even though the grid sizes should still be of the same order). 

\begin{figure}[ht]
\begin{centering}
\includegraphics[scale=0.6]{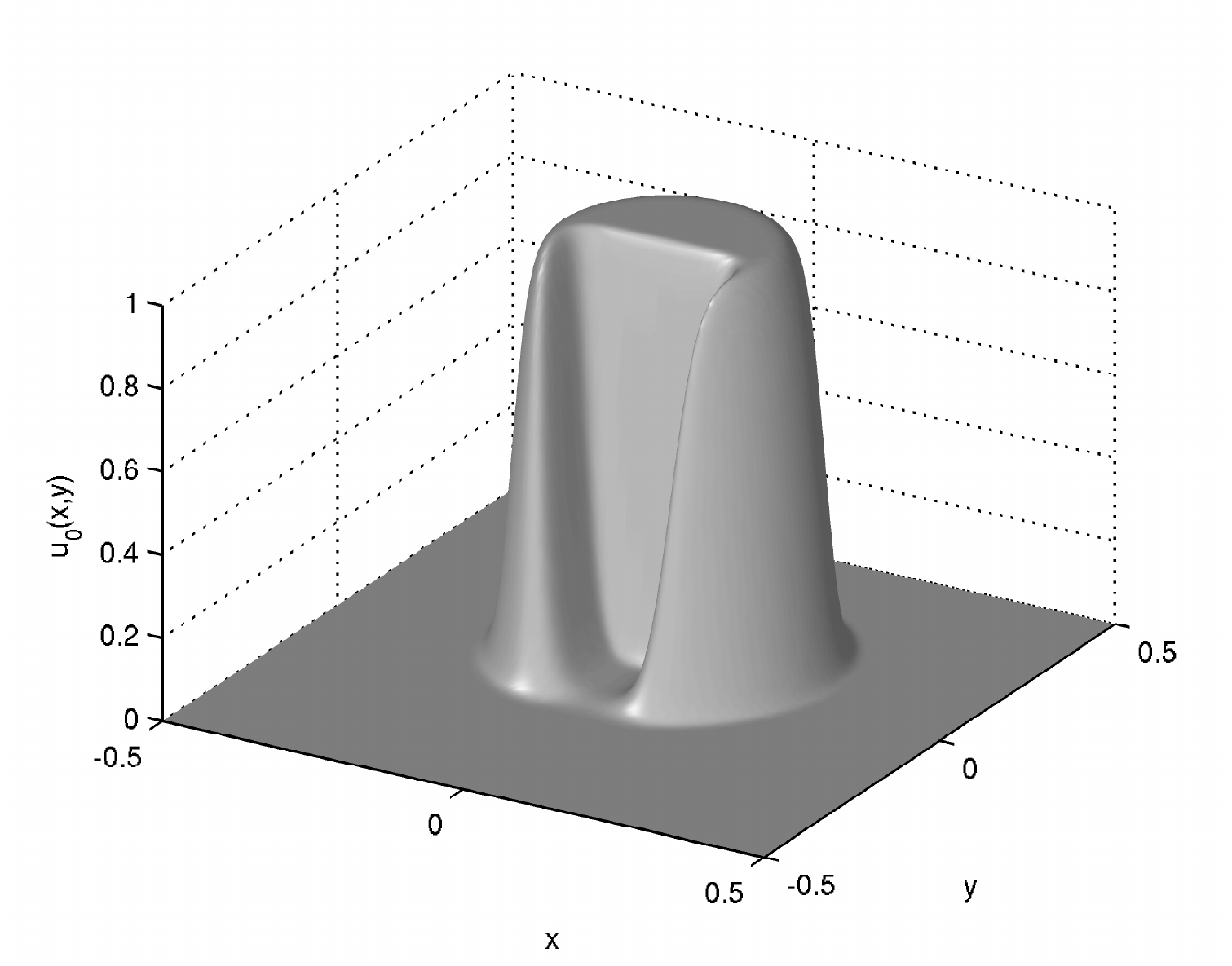}
\includegraphics[scale=0.6]{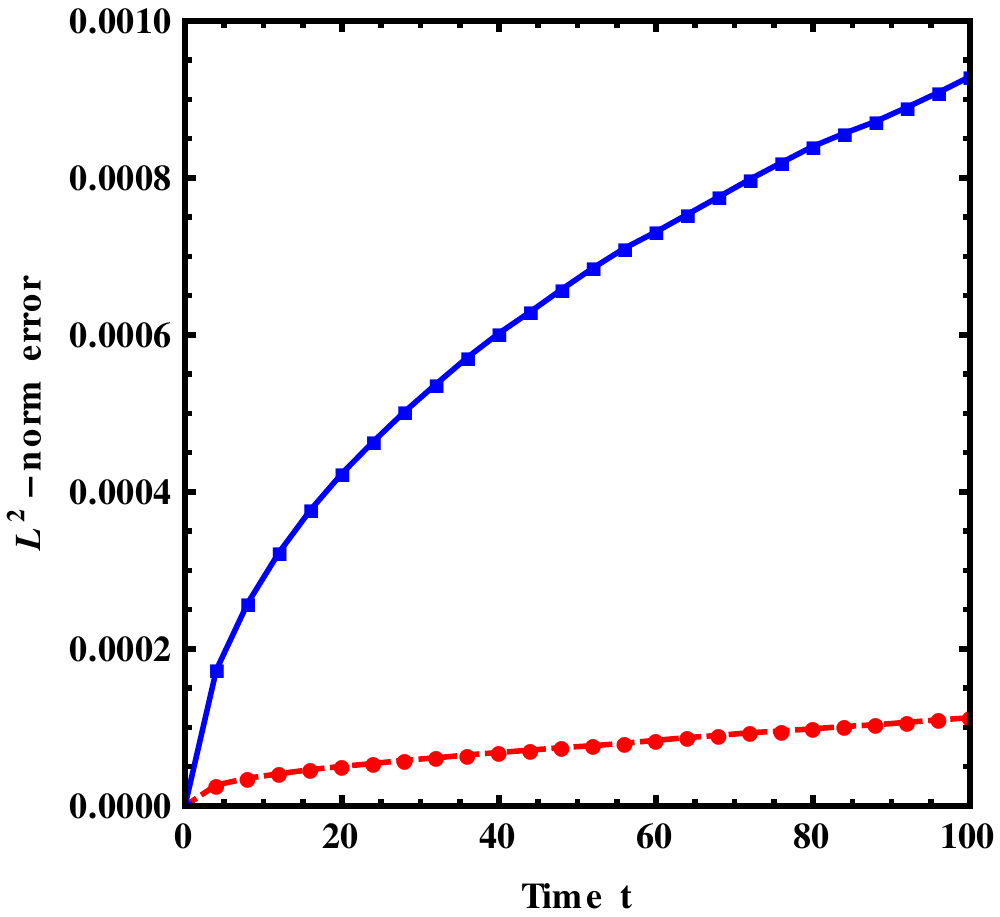}
\end{centering}
\caption{Left: The initial condition $u_0(x,y)$. Right: The $L^2$-norm error against time $t$ with $\Delta x=1/100$ for two spatial discretization schemes, KLL2ND (bottom) and UNOPT10TH (top). }
\label{Fig:7_1}
\end{figure}

We solve the two-dimensional problem by two spatial discretization schemes, KLL2ND and UNOPT10TH.  The right panel of Figure \ref{Fig:7_1} is the $L^2$-norm errors of a long-time simulation for the two methods. The top (solid) line is UNOPT10TH and the bottom (dashed ) line is KLL2ND. The figure shows that the error of KLL2ND almost grows linearly at this scale, and is an order smaller than the error of UNOPT10TH. The time integrator for this simulation is the RK8 method with the CFL number $r=0.1$. Figure \ref{Fig:7_2}(A) and (B) are the grid refinement study of the errors for the two methods, using two different time integrators. The time integrators for (A) and (B) are the RK8 method and the MATLAB ODE solver, ODE45, respectively. The CFL number for RK8 is $r=0.1$, while the options for ODE45 are RelTol = $10^{-5}$ and AbsTol = $10^{-7}$. In Figure \ref{Fig:7_2}(A), because the CFL number is small enough for RK8, the time integration error is negligible, the error dynamics is predicted by our analysis, as described in Fig \ref{Fig:2_5}. However, for Figure \ref{Fig:7_2}(B), we observe that when the grid size is small enough, the error is dominated by the time integrator. Nevertheless, before this happens, Figures \ref{Fig:7_2}(A) and (B) are visually indistinguishable. 

\begin{figure}[ht]
\begin{centering}
\begin{tabular}{c}
\includegraphics[scale=0.6]{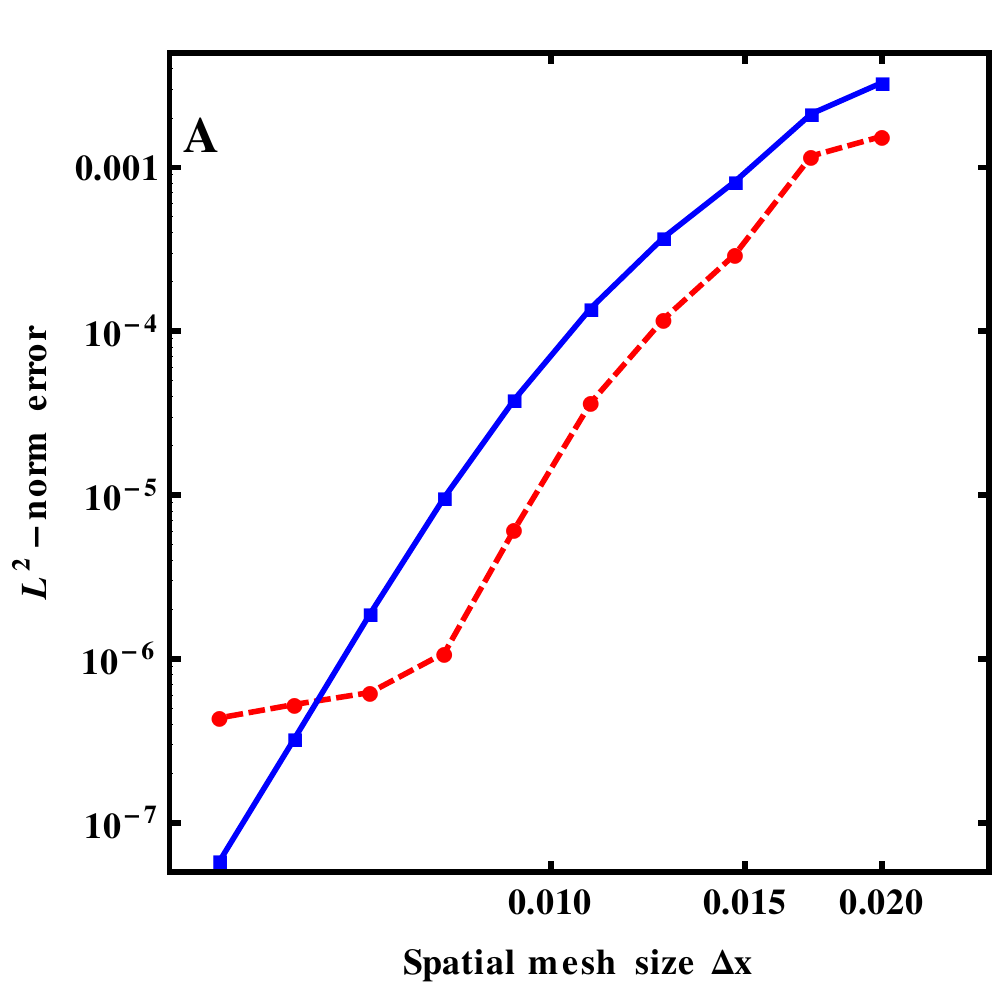}
\includegraphics[scale=0.6]{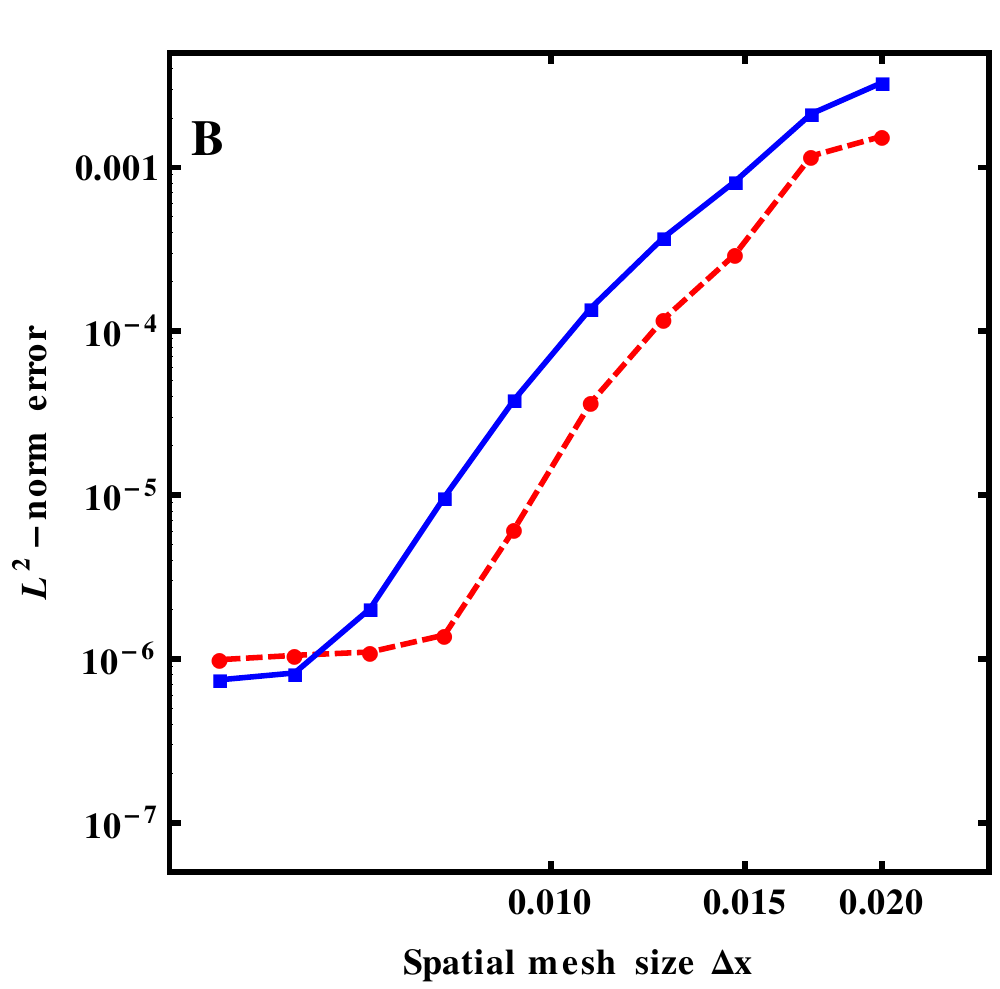}
\end{tabular}
\end{centering}
\caption{The $L^2$-norm error versus the spatial grid size $\Delta x$ at the final time $T=1$. The tested spatial discretization schemes are KLL2ND (dashed line) and UNOPT10TH (solid line). The time integrators are (A) RK8 and (B) MATLAB ODE solver, ODE45. The CFL number for RK8 is $r=0.1$, while the options for ODE45 are RelTol = $10^{-5}$ and AbsTol = $10^{-7}$.}
\label{Fig:7_2}
\end{figure}



\end{example}

\section{Non-periodic boundary conditions}\label{sec:non-periodic}
Our analysis in the previous sections is based on the assumption of periodicity. In this case, we are able to write down the explicit expression of the resolution characteristic $\kappa^*$ as in \eqref{2_20}, for the finite difference spatial discretization precisely because we have complete knowledge of the matrix ${\bf D=L^{-1}R}$ and its associated eigenvalue problem. Periodicity guarantees translational invariance, suggesting that the components of the eigenvectors of ${\bf D}$ are given by sinusoidal waves $\me^{\mi kx_j}$. Thus, the expression of $\kappa^*$ can be obtained by simply substituting the ansatz $U_j=\me^{\mi kx_j}=\me^{\mi kj\Delta x}$ into \eqref{2_10}. 
For problems with non-periodic boundary conditions, we expect that we will be forced to use an asymmetric spatial discretization for the boundary points. This asymmetry destroys the cyclic structure of the matrix ${\bf D}$, and consequently we do not know the exact form of the corresponding eigenvectors. Thus, we can no longer expect to obtain the explicit expression of the resolution characteristic. 
A common practice is to ignore this difficulty by substituting the same ansatz $U_j=\me^{\mi kx_j}=\me^{\mi kj\Delta x}$ into the asymmetric spatial discretization, regardless of the fact that the structure of the matrix ${\bf D}$ has changed, to obtain what we call a pseudo-resolution characteristic. 
This approach is not theoretically sound, and in general does not work. The reason again is that the cyclic structure of the matrix $\bf D$ is altered after an asymmetric discretization is introduced. Thus the sinusoidal wave  $U_j=\me^{\mi kx_j}$ lose their distinguished significance, and the previous analysis collapses.

\subsection{An example of error dynamics of composite schemes}\label{sec:ibvp}

In this section, we examine the transport equation with outflow boundary condition on one side of the domain. We use an example to illustrate the $L^{2}$-norm errors of several composite schemes (an interior scheme combined with a boundary scheme), and we show that the previous analysis does not predict the error when the wave travels into the boundary.

 
Thus, we consider Equation \eqref{2_5}, and we suppose that the wave speed $c>0$ is a constant. Since $c>0$, we are not permitted to prescribe the boundary condition at $x=\ell$. It is determined by the solution. We assume that the wave profile leaves the domain when it reaches the right boundary at $x=\ell$.
%
We let $\ell=1$, $c=1$, and we consider the initial condition suggested by Tam \cite{tam04},
\begin{equation}\label{5_2}
u_0(x)=\exp\left(-\dfrac{25}{9}\log2(x-\dfrac{1}{3})^2\right)\left(2+\cos\left(k_0(x-\frac{1}{3})\right)\right).
\end{equation}
The boundary condition at $x=0$ is $u(0,t)\equiv0$. The wavenumber $k_0=255$ is fixed throughout this section, unless stated otherwise. The leading edge of the wave packet is located at about $x=0.5$. This edge reaches the boundary ($x=1$) at about time $t=0.5$.


Consider a mesh grid with uniformly spaced grid points $x_j=j\Delta x$, $\Delta x=\ell/N$, $0\leq j\leq N$. The boundary points are $x_0$ and $x_N$.  We introduce three spatial discretization schemes, for which the spatial derivatives of the interior points, $1\leq j\leq N-3$ (assuming that $U_j=U_j'=0$ for $j\le 0$) are computed by the standard compact scheme (\ref{2_10}), while derivatives of the  points near the boundary are evaluated by the following boundary scheme reported in \cite{kim07}
\begin{multline}\label{5_3}
\beta_{-2} U'_{N-4}+  \alpha_{-2} U'_{N-3}+U'_{N-2}+\alpha_{2} U'_{N-1}+\beta_{2} U'_{N}\\
=\dfrac{1}{\Delta x}\left(a_2 U_{N}+b_2 U_{N-1}+c_2 U_{N-2}+d_2 U_{N-3}+e_2 U_{N-4}+f_2 U_{N-5}+g_2 U_{N-6}\right)\,,
\end{multline} 
\begin{multline}\label{5_4}
\beta_{-1} U'_{N-3}+  \alpha_{-1} U'_{N-2}+U'_{N-1}+\alpha_{1} U'_{N}\\
=\dfrac{1}{\Delta x}\left(a_1 U_{N}+b_1 U_{N-1}+c_1 U_{N-2}+d_1 U_{N-3}+e_1 U_{N-4}+f_1 U_{N-5}+g_1 U_{N-6}\right)\,,
\end{multline} 
\begin{multline}\label{5_5}
\beta_{0} U'_{N-2}+  \alpha_{0} U'_{N-1}+U'_{N}\\
=\dfrac{1}{\Delta x}\left(a_0 U_{N}+b_0 U_{N-1}+c_0 U_{N-2}+d_0 U_{N-3}+e_0 U_{N-4}+f_0 U_{N-5}+g_0 U_{N-6}\right)\,,
\end{multline} 
for $j=N-2,N-1,N$.  When matching the coefficients of the above equations, the highest order of accuracy equations (\ref{5_3}),  (\ref{5_4}), and (\ref{5_5}) can achieve are tenth-order, ninth-order, and eighth-order, respectively. Hence, the unoptimized boundary scheme overall is eighth-order accurate. We refer to this asymmetric boundary scheme as UNOPTBS

We test three spatial discretization schemes in this section. The first one is a composite scheme of the tenth-order unoptimized scheme (\ref{2_10}) for the interior points and the unoptimized boundary scheme UNOPTBS for the boundary points. We refer to this composite scheme as UNOPT8TH. The second one is the scheme reported in \cite{kim07}. This is a fourth-order composite optimized scheme. Both the interior and the boundary schemes are optimized. We refer to this scheme as KIM4THOPTBC. The last scheme is a fourth-order composite scheme that uses the fourth-order optimized interior scheme from KIM4TH and the unoptimized boundary scheme UNOPTBS. We refer to this scheme as 4THUNOPTBS. 
The computational costs for these three schemes are of the same order. 
\begin{figure}[ht]
\begin{centering}
\begin{tabular}{cc}
\includegraphics[scale=0.4]{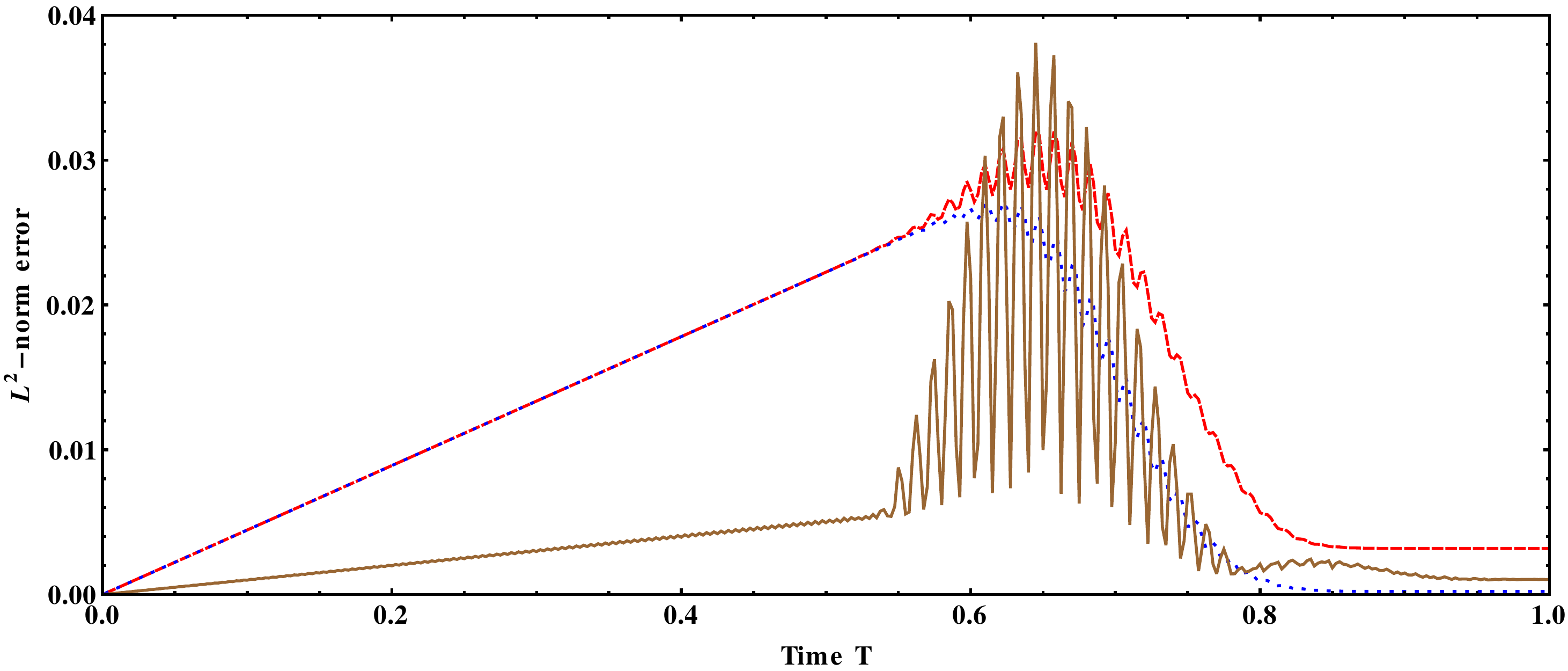}\\
\includegraphics[scale=0.405]{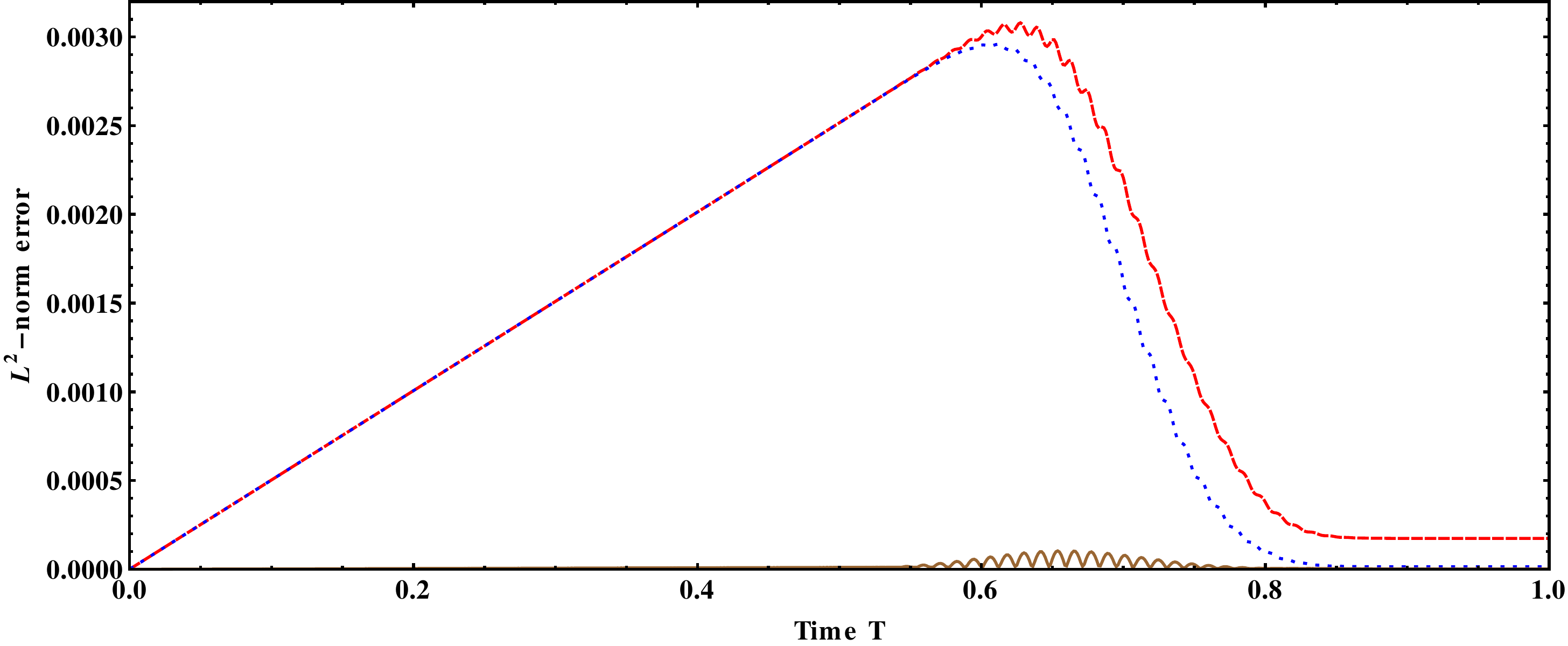}
\end{tabular}
\end{centering}
\caption{The resulting $L^2$-norm error against time with $\Delta x=1/200$ (top) and $\Delta x=1/400$ (bottom). Three schemes are used. Solid line: UNOP8TH. Dashed line: KIM4THOPTBC. Dotted line: 4THUNOPTBS.}
\label{Fig:5_1}
\end{figure}

Figure \ref{Fig:5_1} shows the $L^2$-norm errors of the three tested schemes against time for $\Delta x=1/200$ (top panel) and $\Delta x=1/400$ (bottom panel), respectively. We use RK8 as the time integrator with the CFL number $r=0.1$ for all calculations.  The solid line is UNOPT8TH, the dashed line is KIM4THOPTBC, and the dotted line is 4THUNOPTBS.
It is clear from Figure \ref{Fig:5_1} that before the wave front reaches the boundary $x=1$ (before $t=0.5$), the error corresponding to each scheme grows linearly with time. This is consistent with our analysis for the periodic case, despite the fact that the sinusoidal wave $U_j=\me^{\mi kx_j}$ is not an eigenvector of $\bf{D}$ (as defined in \S\ref{sec:compact}) anymore. Also, since KIM4THOPTBC and 4THUNOPTBS use the same interior scheme, the corresponding errors shown in Figure \ref{Fig:5_1} coincide with each other. The errors are no longer predictable by the analysis after the wave packet enters the boundary. It is interesting to see that in Figure \ref{Fig:5_1}, for the two schemes using the same interior optimized scheme, KIM4THOPTBC and 4THUNOPTBS, the one using unoptimized boundary scheme (4THUNOPTBS) produces smaller $L^2$-norm error than its counterpart for time between 0.58 and 0.8, when wave packet enters the boundary.  While all three schemes produce oscillatory errors, when the wave packet is traveling into the boundary, the error oscillation of the eight-order unoptimized scheme, UNOPT8TH is prominent for gird size $\Delta x=1/200$. Nevertheless, its error on average is still comparable to that of the optimized schemes. For the finer grid ($\Delta x=1/400$), UNOPT8TH produces considerably smaller error than the optimized schemes before and after the wave packet enters the boundary. Our aim for the rest of this section is to analyze and quantify the error dynamics near the boundary for composite schemes. Moreover, we propose a composite scheme based on this discussion. 

\subsection{Error analysis for composite schemes}
Suppose that at time $T=n\Delta t$, a sinusoidal wave on a mesh grid is ${\bf U}^n=(U_1^n,\cdots,U_N^n)^T$, where $U_j^n=\me^{\mi\alpha x_j}$ and $x_j=j\Delta x$ is the $j^{th}$ grid point. Similar to the case of periodic boundary conditions, for solutions at $T+\Delta t$, the spatial derivative at the grid points are evaluated by a composite compact scheme that can be viewed as a linear transformation $({\bf U}^n)' = {\bf D} {\bf U}^n$. Following Jordan and Sengupta et al.\ \cite{j07, sds07}, for any given $\alpha$ and $x_j$, $j=1,\cdots N$, we should be able find a diagonal matrix $\bf{K}$ so that 
\begin{equation}\label{5_8}
{\bf DU}^n = {\bf KU}^n\,,
\end{equation} 
where the diagonal entries of ${\bf K}$ are $k_1,\ldots,k_N$. In general $k_i\ne k_j$, if $i\ne j$. Thus, this matrix ${\bf K}$ signifies a phenomenon that does not arise in the case of periodic boundary conditions. (We recall that in that case ${\bf DU}^n = k^{*}{\bf I}{\bf U}^n$, for some scalar $k^{*}$, where ${\bf I}$ is the identity matrix.) Now we define
\begin{equation}\label{5_9}
\Delta{\bf K}={\bf K}-k^*{\bf I}\,.
\end{equation}
Thus, $\Delta{\bf K}$ provides a measurement of the the difference between the numerical wavenumbers that appear in composite schemes and those that appear in schemes for periodic boundary conditions. 
Suppose we use a two-stage second-order Runge-Kutta method to find the numerical solution at $T+\Delta t$.  From \eqref{2_30}, the time-stepping is given by
\begin{equation}\label{5_10}
\begin{split}
\textbf{K}_1&=-c\Delta t\textbf{D}\textbf{U}^n\,,\\
\textbf{K}_2&=-c\Delta t\textbf{D}(\textbf{U}^n + \dfrac{1}{2!}\textbf{K}_{1})\,,\\
\textbf{U}^{n+1}&=\textbf{U}^n + \textbf{K}_2\,.
\end{split}
\end{equation}
Substituting  ${\bf DU}^n = {\bf KU}^n$ and ${\bf K}=k^*{\bf I}+\Delta{\bf K}$ into \eqref{5_10}, we immediately obtain
\begin{equation}\label{5_11}
\textbf{U}^{n+1}=\left(1+(-ck^*\Delta t)+\dfrac{(-ck^*\Delta t)^2}{2!}\right)\textbf{U}^{n}+O(\Vert\Delta t\Delta\textbf{KU}^n\Vert)\,.
\end{equation}
In \eqref{5_11}, $\|\cdot\|$ denotes a suitable vector norm. If we compare \eqref{5_11} with \eqref{2_34}, we see immediately that the first term in \eqref{5_11} agrees with the right-hand side of \eqref{2_34} (which assumes periodic boundary conditions). The second term in Equation \eqref{5_11} thus accounts for the discrepancy between schemes with periodic and non-periodic boundary conditions. Furthermore, the form of the second term in \eqref{5_11} holds for general Runge-Kutta methods of any order.


\begin{figure}[ht]
\begin{centering}
\begin{tabular}{c}
\includegraphics[scale=0.35]{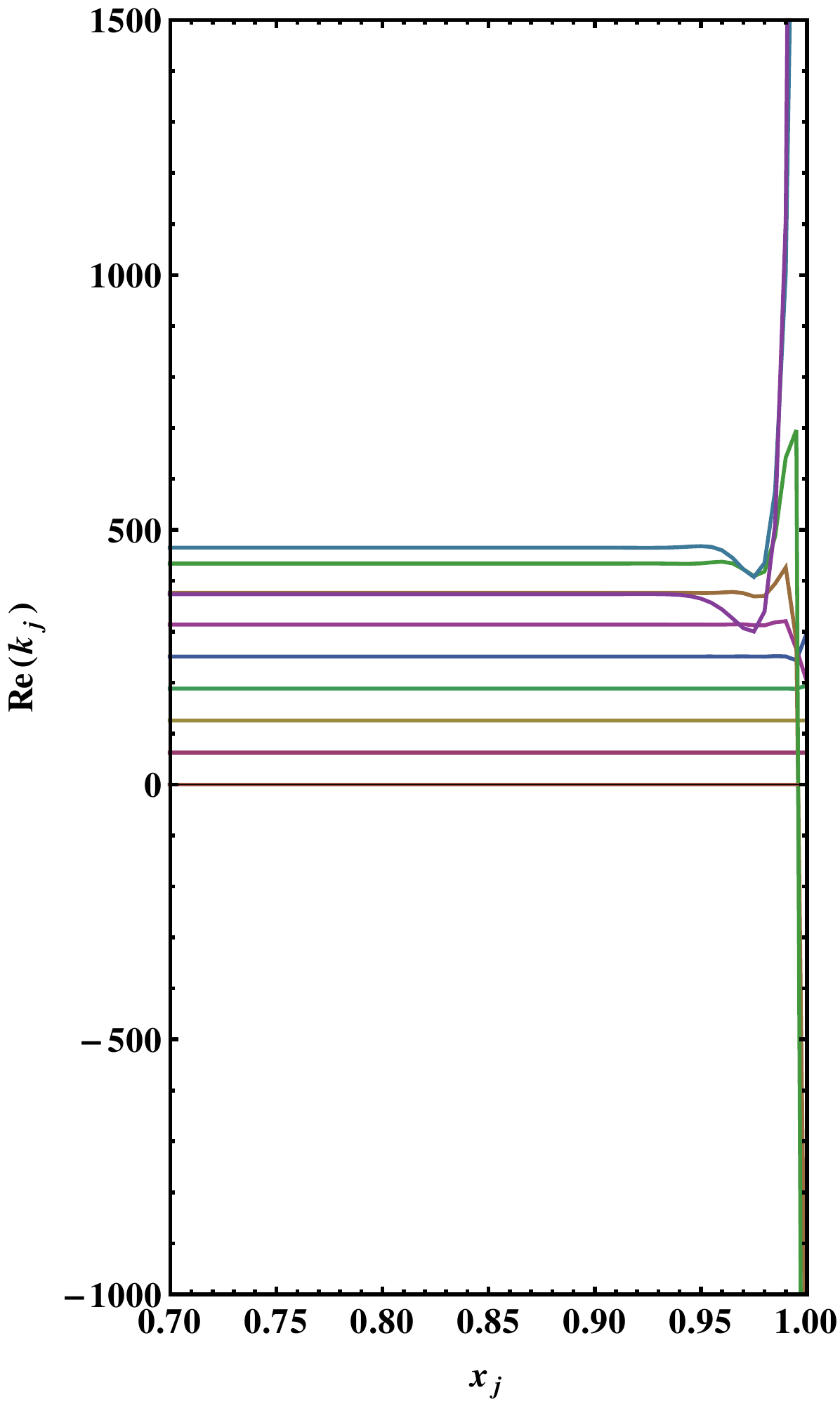}\hskip 0.5in
\includegraphics[scale=0.35]{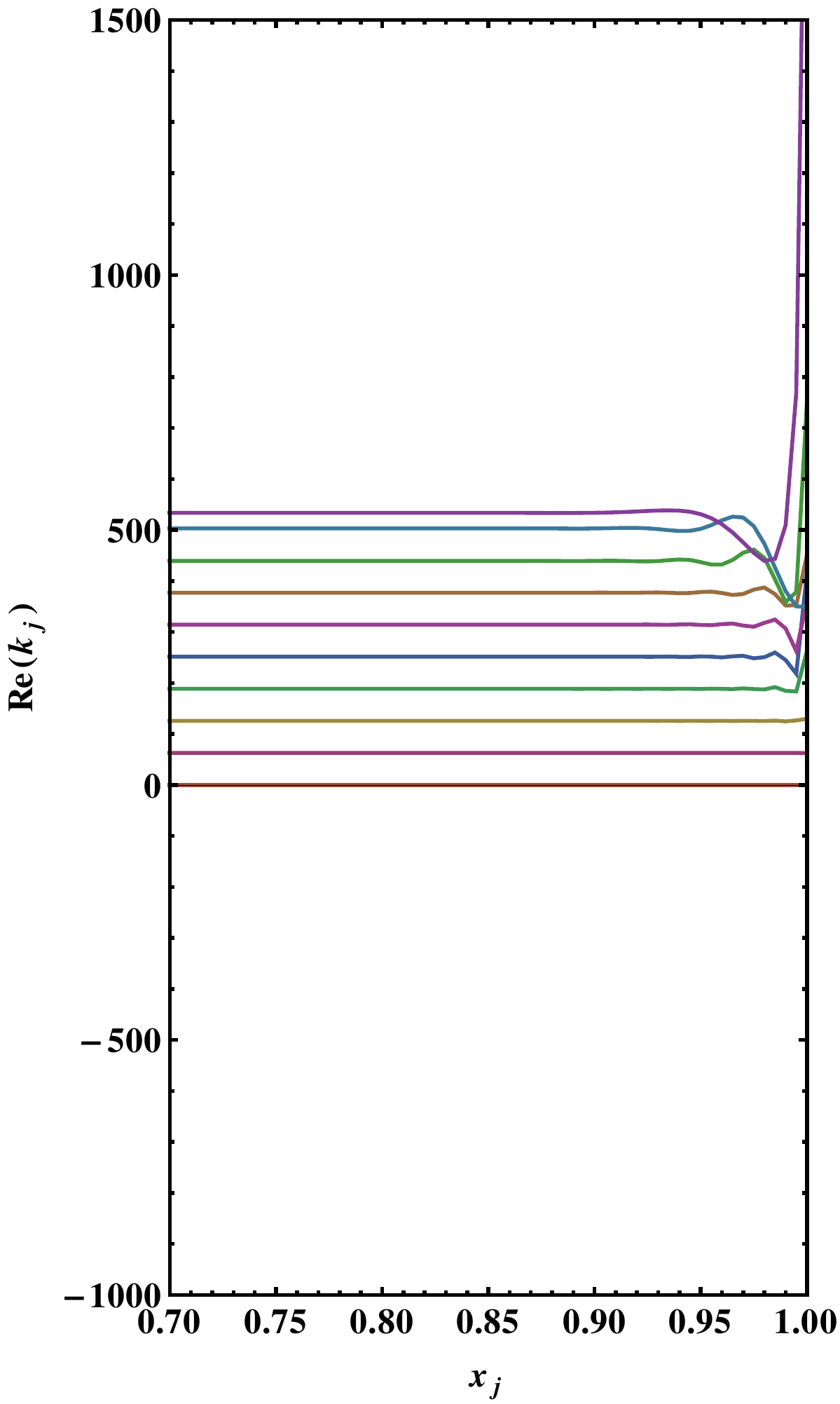}
\end{tabular}
\end{centering}
\caption{Left: The plot of the real part of $k_j$ (denoted $Re(k_j)$) against $x_j$ for the composite scheme UNOPT8TH. The $\alpha$ values used are $\alpha=2\pi k'$, $k'=0,10,20,\ldots,90$.  Right: the same as the left but for the composite scheme KIM4THOPTBC}
\label{Fig:5_2}
\end{figure}

Consider the composite schemes UNOPT8TH and KIM4THOPTBC. For each scheme, we compute $k_j$, $j=1,\ldots, 200$, and $\Delta x =1/200$ in \eqref{5_8}, for various $\alpha$, where $\alpha=2\pi k'$, $k'=0, 10, 20,\ldots,90$. Figure \ref{Fig:5_2} shows the real part of $k_j$ (denoted $\Re(k_j)$) versus $x_j$ for UNOPT8TH on the left and for KIM4THOPTBC on the right. From top to bottom, the $k'$ values are $90, 80, \ldots, 0$. It is clear from the figure that for $x_j\le0.9$, the values of $\Re(x_j)$ are constant for all $k$. Our calculations confirm that within this interval the constant value of $k_j$ and the scalar value of $k^{*}$ from the periodic case are approximately the same. This explains why the error dynamics of the composite schemes are predictable by the analysis, based on the assumption of periodic boundary conditions, before the wave packet travels into the boundary. The figure also show that for $x_j>0.9$, $k_j$ is no longer well predicted by $k^{*}$, and the analysis fails. 

\begin{note}
Based on Figure \ref{Fig:5_1}, Figure \ref{Fig:5_2}, equation \eqref{5_11}, and our numerical experiments, we make the following observations.
\begin{enumerate}
\item[(a)] Equation \eqref{5_11} is derived based on the assumption that the spectrum of the wave packet has one well-defined wavenumber throughout the simulation. In reality, the wave packet may have a finite range of spectrum, and different  wave components may interact with each other. The error dynamics may not  be as simple as what \eqref{5_11} describes.
\item[(b)] Optimized boundary schemes need not produce a smaller $\Delta{\bf K}$ (as in \eqref{5_11}). Evidence for this statement can be seen in the comparison of KIM4THOPTBC and 4THUNOPTBS near the boundary in Figure \ref{Fig:5_1}. 
\item[(c)] From Figure \ref{Fig:5_2} and Equation \eqref{5_11}, we expect that large errors occur near the boundary. This expectation is supported by Figure \ref{Fig:5_1}, where large oscillatory errors are observed among the tested composite schemes near the boundary. However, it is interesting to point out that from our numerical experiments, the frequency of the oscillations does not depend on $\Delta x$ or on what FD scheme is used. Rather, it depends on $k_0$ in the initial condition \eqref{5_2}. A large $k_0$ may result in a large frequency.
\item[(d)] Minimizing $\Delta {\bf K}$ in terms of an optimization problem is difficult. For example, if the objective function is defined by 
\begin{equation}\label{5_12}
\sum_{j=1}^{N}\int_0^{k_{\mathrm{max}}}\vert k_j-k^*\vert^2\,\dif k\,,
\end{equation}
the system giving rise the analytical expression of $k_j$'s is normally overdetermined, except for some special cases \cite{j07}. We explain the difficulty in Appendix \ref{appendix2}.
\end{enumerate}
\end{note}

\subsection{A hybrid composite scheme}

Since minimizing $\Delta {\bf K}$ in terms of an optimization problem is rather difficult, we propose a composite scheme that reduces the destruction of the structure of the interior scheme while coupling with the boundary scheme so that the size of $\|\Delta {\bf K}\|$ is reduced naturally.  

Let us consider again the transport equation with outflow boundary condition on the right. The incoming signal is from the left boundary $x=0$, so that the boundary $u(0, t)$ is given. To illustrate the idea,  as an example, we use the pentadiagonal discretization, \eqref{2_10}, for the interior scheme. Let the mesh grid be $x_j=j\Delta x$, $j=0,\cdots N$. Assuming that the derivatives at $x_j$, $j=1,\cdots N-3$ are evaluated by the interior scheme (provided we have knowledge about the ghost point values $u(-\Delta x, t)$ and $u(-2\Delta x, t)$). Since the outflow boundary condition does not provide information at the ghost points $x_{j}$, $j=N+1,\cdots, N+3$, to close the system and find the derivatives at $x_j$, $j=N-3,\ldots, N$, we have to introduce a boundary scheme, such as Equations \eqref{5_3}, \eqref{5_4}, and \eqref{5_5}.  However, the resulting $N$ equations are fully coupled. This coupling may disturb the original resolution characteristic of the pentadiagonal discretization; see Figure \ref{Fig:5_2}. Therefore, it is desirable to reduce this coupling. We remark that if the implicit (compact) discretization degenerates into an explicit scheme (e.g., $\alpha=\beta=0$, no coupling), the optimization procedure for $\Delta {\bf K}$ becomes feasible. 
 
We propose that instead of adding more equations to close the system, we close the system by reducing the number of unknowns. In other words, if we can get a suitable ansatz for $U'_{j}$, $j=N-2,\ldots, N$, then together with the equations for $1\leq j\leq N-3$, we close the system for the interior scheme on $N$ grid points with $N-3$ unknowns. Furthermore, if the ansatz is given by a method which is known to have an accurate resolution characteristic, we can expect that the resolution characteristic of the overall spatial scheme will benefit from that. 

\begin{figure}[ht]
\includegraphics[scale=0.6]{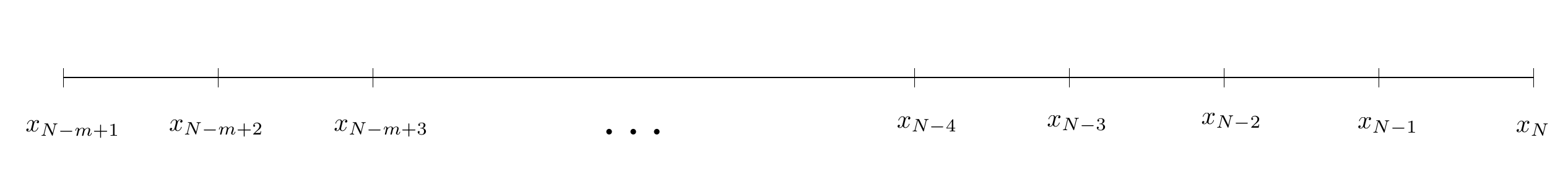}
\caption{Buffer zone's grid points for KLL2NDBC}
\label{Fig:7_4}
\end{figure}

To find a suitable ansatz, we recall that the resolution characteristics of unoptimized compact high-order methods behave nicely at low and moderate modified wavenumbers. Therefore, we propose an ansatz that approximates the derivatives of points near the boundary by an unoptimized high-order method. 
To describe our composite compact DRB scheme, we start with the DRB scheme KLL2ND described in \S\ref{sec:DRB-design}. When $x_N$ is the boundary point at $x=1$, formula \eqref{5_14} is well defined for $j\leq N-5$. In order to find the derivatives at $x_j$, $j=N-4,\ldots N$, we must make an ansatz for $U'_j$, $N-4\leq j\leq N$. 
To do this, we consider a buffer zone near the boundary, consisting of $m$ grid points, including the boundary points. The grid points are labeled as $x_{N-m+1}, x_{N-m+2}, \ldots, x_{N}$, as  shown in Figure \ref{Fig:7_4}.
For the left endpoints of the buffer zone, from $x_{N-m+1}$ to $x_{N-m+3}$, and the right endpoint of the buffer zone, from $x_{N-2}$ to $x_{N}$, we use the boundary scheme, that is, \eqref{5_3}, \eqref{5_4} and \eqref{5_5}, to find the derivatives. For each of the other points inside the buffer zone, from  $x_{N-m+4}$ to $x_{N-3}$, we use the tenth-order unoptimized pentadiagonal scheme, UNOPT10TH, to evaluate the derivatives. After this calculation, we have the derivatives  $U'_j$, $j=(N-m+1),\cdots,N$. Then, we take the values of $U'_j$, for $N-4\leq j\leq N$, as our ansatz, and we use them in \eqref{5_14} to solve for the derivatives for $j=1,\ldots,N-5$.     

\begin{note}[Bad Coupling]
The role of the buffer zone in the construction of the ansatz is to restrict the ``bad coupling'' to a small number, i.e., $m$, of grid points instead of $N$ grid points. In general $m\ll N$, and our numerical experiments suggest that $m\ge 9$. If $m<9$, the $m\times m$ system may be ill-posed.  We refer to our hybrid composite compact DRB scheme as KLL2NDBC
\end{note}

\begin{figure}[ht]
\begin{centering}
\begin{tabular}{cc}
\hskip -0.16in\includegraphics[scale=0.4]{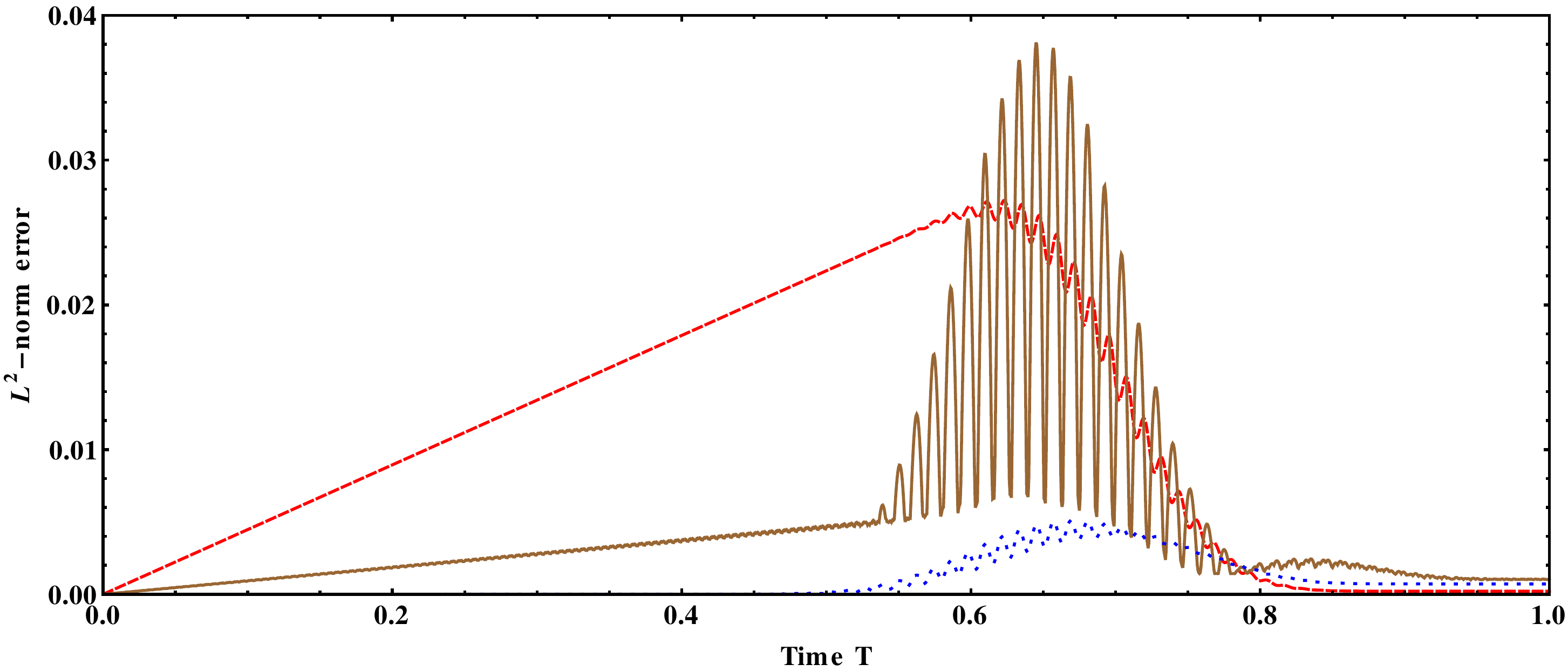}\\
\includegraphics[scale=0.4]{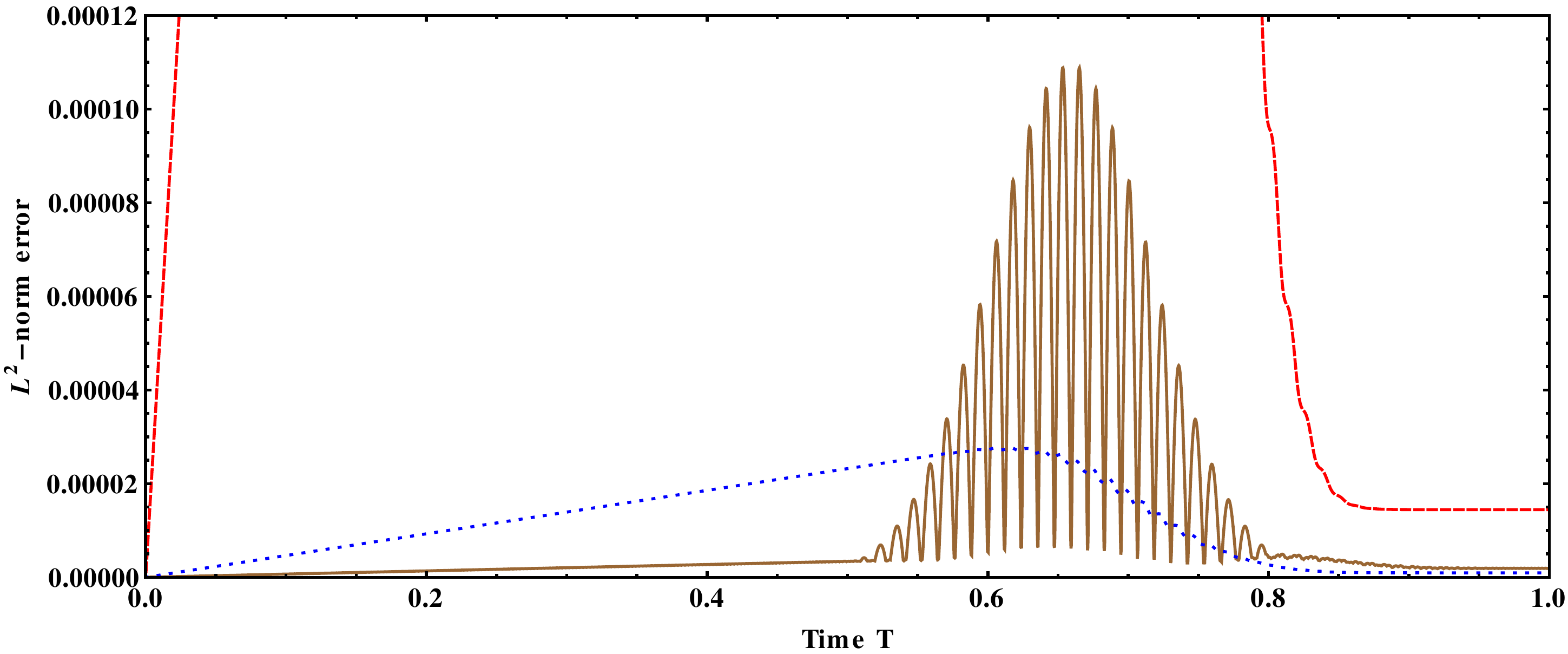}
\end{tabular}
\end{centering}
\caption{The $L^2$-norm error versus time with $\Delta x=1/200$ (top) and $\Delta x=1/400$ (bottom). The three tested spatial discretization schemes are UNOPT8TH (solid line), KIM4THOPTBC (dashed line), and KLL2NDBC (dotted line). The time integrator is RK8 with CFL number $r=0.1$.}
\label{Fig:5_3}
\end{figure}

Now we compute the IBVP from  \S\ref{sec:ibvp} using KLL2NDBC. Figure \ref{Fig:5_3} show the $L^2$-norm errors versus time for the composite schemes UNOPT8TH (solid line), KIM4THOPTBC (dashed line), and KLL2NDBC (dotted line). The time integrator is RK8 and the CFL number is $r=0.1$. The results in the top panel are for grid size $\Delta x =1/200$ while those in the bottom panel are for $\Delta x =1/400$. For both grids, the buffer zone is $m=10$ for KLL2NDBC. The figure shows that the proposed composite scheme produces the least error when the wave packet enters the right boundary.  We remark that using five computed $U_j'$, $j=N-4,\cdots N$, as our ansatz is suggested by our numerical experiments.  The proposed composite scheme may be unstable if more than five points are used for the ansatz.

\section{Concluding remarks}\label{sec:conclusion}

In this paper, we have given a detailed error analysis for optimized compact schemes that minimize the dispersion-and-dissipation errors for FD discretization. The analysis is based on the non-dispersive transport equation in a periodic domain. We derive the $L^2$-norm error of the solution in terms of the dispersion-and-dissipation errors. We have also carried out numerical experiments to compare various optimized schemes and some traditional unoptimized schemes. Our analysis and experiments characterize the interplay between Taylor truncation errors and the dispersion-and-dissipation errors. Based on the analysis and the experiments, we proposed a set of principles for designing a suitable optimized compact scheme for a given problem. Moreover, we used these principles to develop several optimized schemes, including a composite scheme for problems with non-trivial boundary conditions. Our ultimate goal, however, is not to claim that the proposed discretization schemes are superior than the traditional unoptimized schemes, or other optimized schemes. On the contrary, we hope to remind readers that there is no optimized FD scheme that is superior to its higher-order unoptimized counterpart or other optimized FD schemes on all meshes. When we select or design a FD algorithm, especially an optimized FD scheme, for solving a given PDE, we need to pay attention to the prior knowledge of the solution of the equation, as well as the following issues regarding the problem. They are as follows. (1) What is the error tolerance we are willing to allow for the solution? (2) What computer power do we have? This includes, in particular, the maximal memory capacity and the desired (or the limit of) CPU/GPU run-time. (3) What is the highest wavenumber we want to resolve? Because a good numerical algorithm often means a good algorithm for solving a specific problem under specific circumstances, these issues must be kept in mind, and we believe that our analysis will aid users as they seek good algorithms for their problems.

\section*{Acknowledgement}

Research of YHK and GL was supported in part by the National Science Foundation under grant number DMS-0845127.

\appendix

\section{Coefficients for spatial discretizations} 
We report the coefficients of the spatial discretization schemes developed in this paper. Whenever an optimization is performed on a definite integral, the integral is computed using Mathematica's built-in function {\bf NIntegrate} (this function finds numerical value of an integration for a given integrand over a specified domain), and the optimization is carried out by another built-in Mathematica function, {\bf FindMinimum} (this function solves the minimization problem for a given objective function). Note that, in general, the function {\bf NIntegrate} returns a numerical value instead of a function of undetermined independent variables. So, to define a function of certain undetermined variables by {\bf NIntegrate} and later substitute it into the command {\bf FindMinimum}, we use the key word ``{\bf \_?NumberQ}" to ensure the undetermined variables  are not treated as numerical values. The coefficients and the objective functions are listed below.

\begin{enumerate}

\item OPT2ND1P5: The discretization is given by \eqref{2_10}, and the coefficients are found by minimizing the objective function $\int_0^{\Gamma}|\kappa^*-\kappa|^2\,\dif\kappa$ subject to the constraint (the first equation in \eqref{eq:padeconstraints}) and $\Gamma=1.5$. The function $\kappa^*$ is specified by \eqref{2_20}.
\begin{equation*}\label{8_2}
\begin{split}
a&=1.382433766621239,\quad\alpha=0.5267276883269482,\\
b&=0.7782394862770302,\quad\beta=0.06166491031289359,\\
c&=0.01611194438141422.
\end{split}
\end{equation*}

\item OPT2ND1P8: The discretization is given by \eqref{2_10}, and the coefficients are found by minimizing the objective function $\int_0^{\Gamma}|\kappa^*-\kappa|^2\,\dif\kappa$ subject to the constraint (the first equation in \eqref{eq:padeconstraints}) and $\Gamma=1.8$. The function $\kappa^*$ is specified by \eqref{2_20}.
\begin{equation*}\label{8_3}
\begin{split}
a&= 1.361203029457461,\quad\alpha=0.5412576434842603,\\
b&=0.8382608566836231,\quad\beta=0.06893488676613146,\\
c&=0.02092117435969928.
\end{split}
\end{equation*}


\item KLL2ND: The discretization is given by \eqref{3_7}, and the coefficients are found by minimizing the objective function $\int_0^{\Gamma}|\kappa^*-\kappa|^2\,\dif\kappa$ subject to a constraint that guarantees a second-order of formal accuracy and $\Gamma=2.0$. The function $\kappa^*$ is specified by \eqref{3_6}. The second-order constraint can be derived by assuming the condition $\kappa^*=\kappa+O(\kappa^3)$. 
\begin{equation*}\label{8_5}
\begin{split}
a&= 1.271681048997683,\quad d=0.01957068852632900805,\\
b&=1.0324635517800887,\quad e=-0.024173863453322705888322,\\
c&=0.0710624871538077965,\quad f=0.001112355198712081607526,\\
\alpha&=0.585636483962933,\quad\beta=0.09783577170489735.
\end{split}
\end{equation*}

\item UNOPT8TH: The discretizations are given by \eqref{5_3}, \eqref{5_4}, and \eqref{5_5}. The coefficients are found so that each of the discretization can achieve its highest order of formal accuracy. 
\begin{equation*}\label{8_6}
\begin{split}
\alpha_{-2}&=\dfrac{8}{9},\quad\beta_{-2}=\dfrac{1}{6},\quad\alpha_{2}=\dfrac{4}{15},\quad\beta_{2}=\dfrac{1}{90},\\
a_2&=-\dfrac{68}{675},\quad b_2=-\dfrac{254}{225},\quad c_2=-\dfrac{7}{6},\quad d_2=\dfrac{40}{27},\\
e_2&=\dfrac{8}{9},\quad f_2=\dfrac{2}{75},\quad g_2=-\dfrac{1}{1350},\\
\alpha_{-1}&=\dfrac{5}{2},\quad\beta_{-1}=\dfrac{10}{9},\quad\alpha_{1}=\dfrac{1}{18},\\
a_1&=-\dfrac{257}{540},\quad b_1=-\dfrac{107}{30},\quad c_1=-\dfrac{5}{12},\quad d_1=\dfrac{110}{27},\\
e_1&=\dfrac{5}{12},\quad f_1=-\dfrac{1}{30},\quad g_1=\dfrac{1}{540},\\
\alpha_{0}&=12,\quad\beta_{0}=15,\\
a_0&=-\dfrac{79}{10},\quad b_0=-\dfrac{154}{5},\quad c_0=\dfrac{55}{2},\quad d_0=\dfrac{40}{3},\\
e_0&=-\dfrac{5}{2},\quad f_0=\dfrac{2}{5},\quad g_0=-\dfrac{1}{30}.
\end{split}
\end{equation*}

\item UNOPT10TH: Following Lele \cite{lele92}, the discretization is given by \eqref{2_10}, where the coefficients satisfy (\ref{eq:padeconstraints}).
\begin{equation}\label{8_1}
\begin{split}
a&=\dfrac{17}{12},\quad\alpha=\dfrac{1}{2},\\
b&=\dfrac{101}{150},\quad\beta=\dfrac{1}{20},\\
c&=\dfrac{1}{100}.
\end{split}
\end{equation}

\item KLL2NDBC: The coefficients for the boundary points in the buffer zone can be found in UNOPT8TH. The coefficients for the interior points of the buffer zone can be found in UNOPT10TH. The coefficients of the interior points outside the buffer zone can be found in KLL2ND.


\end{enumerate}

\section{Difficulty of minimizing $||\Delta{\bf K}||$ in terms of an optimization problem}\label{appendix2}

In order to minimize $\Vert\Delta{\bf K}\Vert$, it is necessary to establish an objective function for the optimization process. Equation \eqref{5_12} is an example. While other objective functions might be considered, a critical factor for a useful objective function is that analytical expressions for the $k_j$'s must be available. We claim that it is most unlikely to find such analytical expressions except in the case of explicit schemes and for a very special implicit scheme \cite{j07}. In this appendix we indicate why this is so.

To be concrete, consider the computational domain $0\leq x\leq1$, with $x_j=j/N=j\Delta x$, $0\leq j\leq N$. For simplicity, we use \eqref{2_10} to evaluate derivatives at the interior points, and we use \eqref{5_3}, \eqref{5_4}, and \eqref{5_5} for the points near the right-hand boundary $x=1$. The argument is the same for the boundary $x=0$. With the spatial discretization in hand, we construct the matrices ${\bf L}=[\ell_1\,\ell_2\,\cdots\,\ell_N]^T$ and ${\bf R}=[r_1\,r_2\,\cdots\,r_N]^T$, where the $\ell_j$'s and $r_j$'s are column vectors whose nonzero entries are determined by the coefficients of the spatial discretization. For a given function on $0\leq x\leq1$ whose numerical approximation is denoted by ${\bf U}$, the corresponding approximation for the first derivative is given by ${\bf U'}={\bf DU}={\bf L}^{-1}{\bf RU}.$

Now, fix a wavenumber $k$ and a sinusoidal wave ${\bf U}=[\me^{\mi kx_1}\,\me^{\mi kx_2}\,\cdots\, \me^{\mi kx_N}]^T$. We recall that the diagonal matrix ${\bf K}$ and its entries, the $k_j$'s, are defined by \eqref{5_8}. Equivalently, the diagonal entries of $\mathbf{K}$ must be chosen to satisfy
\begin{equation}\label{9_1}
{\bf LKU}={\bf RU}\,.
\end{equation}
However, this is impossible unless further assumptions are imposed to reduce the size of the system. From Figure \ref{Fig:5_2}, we see that $k_j\approx k^*$ for those $j$'s with $x_j$ being a distance away from the boundary. Hence a reasonable assumption that can help us to reduce the size of the linear system is the following: $k_j\equiv k^*$ for each $j\leq N/2$ (assuming $N$ is even). Here the threshold $N/2$ can be replaced by other numbers, as long as the assumption reflects the observation we made from Figure \ref{Fig:5_2}. With this assumption, let us look at the $(N/2-2)$-th equation in the system (\ref{9_1}). After the common factor that appears on both sides of the equation is cancelled out, the equation becomes 
\begin{multline}\label{9_2}
\beta k_{N/2-4}\me^{-2\mi\kappa}+\alpha k_{N/2-3}\me^{-\mi\kappa}+k_{N/2-2}+\alpha k_{N/2-1}\me^{\mi\kappa}+\beta k_{N/2}\me^{2\mi\kappa}\\
=\dfrac{\mi}{\Delta x}\left(a\sin(\kappa)+\dfrac{b}{2}\sin(2\kappa)+\dfrac{c}{3}\sin(3\kappa)\right).
\end{multline}
Here, we have $k_j=k^*$ for $N/2-4,\ldots,N/2$, due to the assumption we made. Therefore, we have 
\begin{equation}\label{9_3}
\beta k^*\me^{-2\mi\kappa}+\alpha k^*\me^{-\mi\kappa}+k^*+\alpha k^*\me^{\mi\kappa}+\beta k^*\me^{2\mi\kappa}=\dfrac{\mi}{\Delta x}\left(a\sin(\kappa)+\dfrac{b}{2}\sin(2\kappa)+\dfrac{c}{3}\sin(3\kappa)\right).
\end{equation}
For the $(N/2-1)$-th equation in the system \eqref{9_1}, after canceling out common terms and substituting the known $k_j$'s, we have
\begin{multline}\label{9_4}
\beta k^*\me^{-2\mi\kappa}+\alpha k^*\me^{-\mi\kappa}+k^*+\alpha k^*\me^{\mi\kappa}+\beta k_{N/2+1}\me^{2\mi\kappa}\\
=\dfrac{\mi}{\Delta x}\left(a\sin(\kappa)+\dfrac{b}{2}\sin(2\kappa)+\dfrac{c}{3}\sin(3\kappa)\right).
\end{multline}
Now we observe that \eqref{9_3} and \eqref{9_4} are identical, except the term $k_{N/2+1}$ in the latter equation. In fact, the similarity between the two equations implies $k_{N/2+1}\equiv k^*$, provided $\beta\neq0$. Inductively, one can show that $k_j\equiv k^*$ for $j\leq N-1$ by considering the $m$-th equation in the system (\ref{9_1}) for $m\leq N-3$. This leaves us only $k_N$ undetermined. Finally, $k_N$ together with $k^*$ must satisfy the last three equations in the system, which are likely to be overdetermined and have no solution.

In general, if the matrix ${\bf L}$ has $\phi_1$ super- and subdiagonals, and the matrix ${\bf R}$ has $\phi_2$ super- and subdiagonals. Assuming that $\phi_1\leq\phi_2$, the resulting system (\ref{9_1}) then leaves the last $(\phi_2-\phi_1)$ $k_j$'s undetermined and they must satisfy the last $\phi_2$ equations. This results in an overdetermined system for the undetermined $k_j$'s unless $\phi_1=0$. That is, unless the discretization is explicit.

Now let us consider the following two cases. 
\begin{description}
\item[Case (1)] The spatial discretization used for interior points has been specified, but the precise values for the coefficients have not yet been determined. For example, if we choose the pentadiagonal discretization \eqref{2_10} for interior points, but the coefficients $a$, $b$, $c$ and $\alpha$, $\beta$ have not been chosen. Then the system \eqref{9_1} can be reduced to $\phi_2$ equations with $k^*$ and $\phi_2-\phi_1$ $k_j$'s to be determined. If  $\phi_1=1$, we have exactly $\phi_2$ undetermined $k_j$'s, including $k^*$, which satisfy $\phi_2$ equations. Hence we can derive constraints satisfied by the coefficients of the spatial discretization. It is worth pointing out that all $k^*$ and $k_j$'s are functions of $k$ and $\Delta x$, and the system (\ref{9_1}) represents a system of functional relations instead of usual scalar equations. These functional relations impose strong constraints on the $k^*$ and $k_j$'s. Consequently, solving these constraints involving solving a system of polynomials of multivariables instead of a simple linear system as one might have expected. Larger $\phi_2$ implies larger degrees of the polynomials. When $\phi_2$ is too large, we expect the system to be unsolvable, and analytical expressions for $k_j$'s are not available in that case. Jordan \cite{j07} studied this special case $\phi_1=1$ with small $\phi_2=2,3$.
\item[Case (2)] Both the spatial discretization used for the interior points and the corresponding coefficients are specified. In this case, the reduced system from (\ref{9_1}) is overdetermined and does not admit any solution. Analytical expressions for $k_j$'s are unavailable as well. 
\end{description}

\end{document}